\documentclass[12pt]{amsart}
\usepackage{amssymb,amsfonts}
\usepackage{amscd}
\usepackage{amsmath}
\usepackage{amsfonts}
\usepackage{amssymb}
\usepackage{amsthm}
\usepackage{verbatim}
\usepackage{comment}
\usepackage{mathrsfs}
\usepackage{xcolor}
\usepackage{hyperref}
\hypersetup{
    colorlinks=true,
    citecolor=blue, 
     linkcolor=black,     
    urlcolor=magenta, 
}

\usepackage{tikz-cd}

\usepackage[top=1.4in,bottom=1.4in,left=1.2in,right=1.2in]{geometry}
\providecommand{\U}[1]{\protect\rule{.1in}{.1in}}
\providecommand{\U}[1]{\protect\rule{.1in}{.1in}}
\newtheorem{theorem}{Theorem}[section]
\newtheorem{lemma}[theorem]{Lemma}
\newtheorem{proposition}[theorem]{Proposition}

\newtheorem{example}[theorem]{Example}
\newtheorem{definition}[theorem]{Definition}
\newtheorem{definitions}[theorem]{Definitions}

\newtheorem{corollary}[theorem]{Corollary}
\newtheorem{remark}[theorem]{Remark}
\newtheorem{remarks}[theorem]{Remarks}

\numberwithin{equation}{section}

\newcommand{\be}{\begin{equation}}
\newcommand{\ee}{\end{equation}}
\newcommand{\ber}{\begin{eqnarray}}
\newcommand{\eer}{\end{eqnarray}}
\newcommand{\barr}{\begin{array}}
\newcommand{\earr}{\end{array}}

\newcommand{\h}[1]{\widehat{#1}}
\newcommand{\cA}{\mathcal{A}}
\newcommand{\cB}{\mathcal{B}}

\newcommand{\cE}{\mathcal{E}}
\newcommand{\cM}{\mathcal{M}}
\newcommand{\cF}{\mathcal{F}}
\newcommand{\cL}{\mathcal{L}}
\newcommand{\cI}{\mathcal{I}}
\newcommand{\cU}{\mathcal{U}}

\newcommand{\cH}{\mathcal{H}}
\newcommand{\cG}{\mathcal{G}}
\newcommand{\cJ}{\mathcal{J}}

\newcommand{\cN}{\mathcal{N}}

\newcommand{\cQ}{\mathcal{Q}}
\newcommand{\cT}{\mathcal{T}}
\newcommand{\cV}{\mathcal{V}}

\newcommand{\bC}{\mathbb{C}}
\newcommand{\bD}{\mathbb{D}}

\newcommand{\bB}{\mathbb{B}}

\newcommand{\bK}{\mathbb{K}}

\newcommand{\bN}{\mathbb{N}}

\newcommand{\bZ}{\mathbb{Z}}
\newcommand{\bX}{\mathbb{X}}
\newcommand{\bI}{\mathbb{I}}

\newcommand{\fz}{\mathfrak{z}}
\newcommand{\fw}{\mathfrak{w}}

\newcommand{\fI}{{\mathfrak{I}}}
\newcommand{\fJ}{{\mathfrak{J}}}

\newcommand{\fs}{{\mathfrak{s}}}

\newcommand{\fp}{{\mathfrak{p}}}

\newcommand{\sP}{{\mathscr{P}}} 
\newcommand{\sK}{{\mathscr{K}}}
\newcommand{\sL}{{\mathscr{L}}}

\newcommand{\Hs}{H_{\sigma}}
\newcommand{\Ht}{H_{\tau}}
\newcommand{\Hr}{H_{\rho}}

\newcommand{\zs}{\mathfrak{z}_{\sigma}}

\newcommand{\NC}{NC\nobreakdash-}

\newcommand{\oplusi}{\oplus_{\infty}}

\begin{document}
\title{Algebras, Kernels and Correspondences}
\author{Sibaprasad Barik}
\address{Statistics and Mathematics Unit, Indian Statistical Institute, 8th Mile, Mysore Road, Bengaluru, 560059, India}
\email{barik\_if@isibang.ac.in, sibaprasadbarik00@gmail.com}

\author{Paul S. Muhly}
\address{817 Hinman Ave., 4W, Evanston Il, 60202}
\email{paul-muhly@uiowa.edu}

\author{Baruch Solel}
\address{Department of Mathematics, Technion-Israel Institute of Technology, 32000, Haifa, Israel}
\email{mabaruch@technion.ac.il}

\begin{abstract}
    This is a nearly complete draft of our paper in which we attempt to do three things. First, we show that when formulated in the context of reproducing kernel correspondences over Rieffel's category, Normod-$M$, the Kolmogorov decomposition/representation theorem of Ball, Marx, and Vinnikov is functorial in $M$. Second, we investigate relations between our work and the work of Clou\^atre and Hartz concerning multipliers of reproducing kernels. This investigation leads to interesting questions about the geometry of the matricial unit disc associated with a correspondence over $M$. And third, we describe how some of the recent work of Salomon, Shalit and Shamovich may be formulated in our setting. This leads to the notion of subvarieties, their quotients, and their relation with the structure of reproducing correspondences.   
\end{abstract}

\maketitle

\section{Introduction}

While writing our 2013 paper, \emph{Tensorial function theory: from Berezin transforms to Taylor's Taylor series and back}, \cite{MS13}, we realized that Joseph Taylor's matricial approach to noncommutative function theory, described in \cite{Tay73}, could be generalized to certain spaces of functions defined over the collection of normal representations of a $W^*$-algebra. 
Subsequently, after the appearance of the paper by Joseph Ball, Gregory Marx, and Victor Vinnikov, \cite{BMV1}, we began to wonder about how reproducing kernel Hilbert spaces might also be generalized to that setting. 
Recently, we discovered that the main theorem of their paper, a generalization of Kolmogorov's decomposition/representation theorem, \cite[Theorem 3.1]{BMV1}, makes sense in the context of the representation theory of an arbitrary $W^*$-algebra provided the functions involved are properly defined. 
This assertion applies especially to the functions' domains and ranges. 
Moreover, once the proper setting has been established, and the theorem is articulated in such a way that reflects the setup, their proof of \cite[Theorem 3.1]{BMV1} applies almost word-for-word and reveals the surprising fact that their theorem really is independent of the $W^*$-algebra involved. 
So, in more contemporary language, it is fair to announce that Kolmogorov's decomposition theorem and its proof are functorial in the choice of the $W^*$-algebra involved.

For some indication of the details involved, let $M$ be a $W^*$-algebra, i.e., a $C^*$-algebra that is a dual space. Also, let $\widehat{M}$ be the collection of normal representations of $M$ on Hilbert space. 
Then $\widehat{M}$ is the object space of Rieffel's category Normod-$M$ \cite{Rie74m}.  The morphism space is made up from the intertwiners.  
That is, for representations $\sigma$ and $\tau$, $\cI(\sigma,\tau): =\{C\in B(H_{\sigma},H_{\tau})\vert C\sigma(m)=\tau(m)C,\, m\in M\}$\footnote{We find it notationally more convenient to write $\cI(\cdot,\cdot)$ instead of the traditional $\text{Hom}_M(\cdot,\cdot)$.} is the \emph{intertwiner space} determined by $\sigma$ and $\tau$. 
Also, let $E$ be a $W^*$-correspondence over $M$. 
This means, among several things, that $E$ is a bimodule over $M$ where the left action is given by a homomorphism $\varphi_M$ of $M$ into the bounded adjointable operators on $E$, $\cL(E)$.

The correspondence $E$ plays \emph{two} roles.  First, it generates the abstract operator algebra $H^{\infty}(E)$ that is studied here and in a number of other settings in our work. Second, it enables us to build a new representation of $M$, $\sigma^E$, from a given representation $\sigma$ of $M$ called the \emph{representation of $M$ induced from $\sigma$ by $E$}. The intertwiner space $\cI(\sigma^E,\sigma)$ contains the domains on which the functions we study are defined. These functions come from $H^{\infty}(E)$ through a process of evaluation that we shall define.  We do not study all the possibilities for the domains, but we focus on subsets of $\cI(\sigma^E,\sigma)$ called balls or discs that sometime carry weight sequences. The values of our functions take place in linear mapping spaces built from intertwining spaces. A bit more specifically, the values of our functions are taken in mapping spaces of the form $B_*(\cI(\sigma,\tau),B(H_{\sigma},H_{\tau}))$ - the weak-* continuous linear maps from $\cI(\sigma,\tau)$ to $B(H_{\sigma},H_{\tau})$.  The key point to keep in mind is that the entire theory we develop is \emph{fibred} over $\widehat{M}$.

Once our version of Kolmogorov's decomposition theorem is proved we will apply it in a number of settings.  
Perhaps the two most important applications appear in Sections \ref{sec:CPSzego} and \ref{sec:AppCompPick}. 
In Section \ref{sec:CPSzego} we show that the Szego kernel, which is naturally defined on the fibred unit disc $\bD = \{\bD(\sigma)\}_{\sigma\in \widehat{M}}$, where $\bD(\sigma)$ is the open unit disc in $\cI(\sigma^E,\sigma)$, is a complete Pick kernel.
Section \ref{sec:AppCompPick} is devoted to a generalization of the work of Guy Salomon, Orr Shalit and Eli Shamovich
on ``subvarieties of the disc'' and the corresponding quotients of $H^{\infty}(E)$ they introduced in \cite{SSS}.  We will also discuss the relationship between certain nonmatricial kernels that appear in the literature and our current work in Section \ref{nonmatricial}. The geometry of the unit balls and discs will be discussed in Sections 5 and 6 and will be connected to  the work of Raphael Cluoatre and Michael Hartz in \cite{CH}, which also was a source of inspiration for this paper.

\section{Our Setting}\label{sect:Setup}
\subsection{Intertwiner spaces}\label{ssect:Soil}
We follow the standard convention in which a $W^*$-algebra is a $C^*$-algebra that is a dual space. The term von Neumann algebra is reserved for the image of a $W^*$-algebra under a prescribed \emph{normal} representation (cf. \cite{Sak71}). All representations of $W^*$-algebras that we consider are normal.

Throughout, $M$ will denote a $W^*$-algebra with \emph{separable} predual.
We write $\widehat{M}$ for the set of representations of $M$ on a prescribed separable Hilbert space. 
When it is helpful to keep track of the Hilbert space on which a representation $\sigma$ takes  place, we will index the Hilbert space  with $\sigma$.

We emphasize: \emph{Unless stated otherwise, elements $\sigma$ of $\widehat{M}$ need not be faithful (i.e., $\ker(\sigma)$ can be nonzero). They will, however, be unital: $\sigma(1_M)=I_{H_{\sigma}}$.}

\subsubsection{Locally countable sets of representations}It will be necessary to consider some representations of $M$ on certain nonseparable Hilbert spaces.
Specifically, we consider nonempty sets, $\Omega$, which may be uncountable, that index families, $\{\sigma_{\omega}\}_{\omega\in \Omega}$, of normal representations of $M$, each of which acts on a given separable  Hilbert space $H_{\sigma_{\omega}}$. 
The only general requirement that we place on the families is that  multiplicity is countable, i.e., for each $\omega\in \Omega$, $\{\upsilon\in \Omega\vert \sigma_{\upsilon} = \sigma_{\omega}\}$ is a countable set. We call such sets \emph{locally countable}. By default, countable sets are locally countable, as is $\widehat{M}$.

We say that a family of representations of $M$, $\{\sigma_{\omega}\}_{\omega \in \Omega}$, is \emph{fibred over $\Omega$} or \emph{graded by $\Omega$}. We use this terminology to emphasize one of the parallels between \NC function theory and this presentation. 

\subsubsection{Sums of representations, and intertwiner spaces}\label{sss:SumsRepInt}
Given a set of representations of $M$ that is indexed by a locally countable set, $\{\sigma_{\omega}\}_{\omega \in \Omega}$, we form the direct sum, $H_{\Omega}$, of the Hilbert spaces $H_{\sigma_{\omega}}$ and we form the representation, $\bigoplus_{\omega\in \Omega}\sigma_{\omega}$, of $M$ on $H_{\Omega}$ defined by the formula
\begin{equation}\label{eq:sigmaomega}
    \bigoplus_{\omega\in \Omega}\sigma_{\omega}(a)\{h_{\upsilon}\}_{\upsilon \in \Omega}:= \{\sigma_{\omega}(a)h_{\omega}\}_{\omega\in \Omega}\qquad a\in M, \{h_{\upsilon}\}_{\upsilon\in \Omega}\in H_{\Omega}.
\end{equation}
Evidently, $\bigoplus_{\omega\in \Omega}\sigma_{\omega}$ is a normal representation of $M$ on $H_{\Omega}$ whose image,
\begin{equation*}
   (\bigoplus_{\omega\in \Omega}\sigma_{\omega})(M)=\bigoplus_{\omega\in \Omega}\sigma_{\omega}(M), 
\end{equation*}
is a von Neumann algebra in $B(H_{\Omega})$. It is fruitful to think of the image of $\bigoplus_{\omega\in \Omega}\sigma_{\omega}$ as block diagonal matrices where the blocks are parametrized by $\Omega$.

The \emph{commutant} of $(\bigoplus_{\omega\in \Omega}\sigma_{\omega})(M)$, then, is a \emph{matrix} of intertwining operator spaces indexed by $\Omega \times \Omega$:
\begin{equation}\label{eq:keyint}
    (\bigoplus_{\omega\in \Omega}\sigma_{\omega}(M))'= [\cI(\sigma_{\upsilon},\sigma_{\omega})]_{(\omega,\upsilon)\in \Omega\times \Omega},
\end{equation}
where, for $\sigma$ and $\tau$ in $\Omega$,
\begin{equation*}
    \cI(\sigma,\tau):= \{C\in B(H_{\sigma},H_\tau):\; C\sigma(a)=\tau(a)C,\; a\in M\}.
\end{equation*}
So, in particular, $\cI(\sigma,\sigma)= \sigma(M)'$ acting on $H_{\sigma}$.

Evidently, $\cI(\sigma,\tau)$ is an ultraweakly closed linear subspace of $B(H_{\sigma},H_\tau)$ and the containment\footnote{$\cI(\tau, \omega)\cdot \cI(\sigma, \tau)$ stands for the set of all products $TS$ where $T\in \cI(\tau, \omega)$ and $S\in \cI(\sigma, \tau)$.}
\begin{equation}\label{eq:maps}
    \cI(\tau, \omega)\cdot \cI(\sigma, \tau) \subseteq \cI(\sigma, \omega), \qquad \sigma, \tau, \omega \in \Omega,
\end{equation}
is easily seen to be satisfied by direct calculation. However, a bit more revealing calculation proves the same inclusions and with no extra effort: Form the commutant of the matrix $\text{diag}(\sigma,\tau,\omega)$, viz.,
$$\begin{bmatrix}
    \cI(\sigma,\sigma) & \cI(\tau,\sigma) & \cI(\omega,\sigma)\\
    \cI(\sigma,\tau) & \cI(\tau,\tau) & \cI(\omega,\tau)\\
    \cI(\sigma,\omega) & \cI(\tau,\omega) & \cI(\omega,\omega)
\end{bmatrix}$$
acting on $H_{\sigma}\oplus H_{\tau}\oplus H_{\omega}$, and then calculate its square.

This calculation also reveals 
\begin{enumerate}
    \item that $\cI(\sigma,\tau)$ is a right $\sigma(M)'$ -- left $\tau(M)'$ bimodule and
    \item that the map $C\to C^*$ is an antilinear bimodule isomorphism from $\cI(\sigma,\tau)$ onto $\cI(\tau,\sigma)$.
\end{enumerate}
We shall therefore write $\cI(\sigma,\tau)^*= \cI(\tau,\sigma)$.

These computations show that $\widehat{M}$ is a $W^*$-category in the sense of \cite{GLR1985}. Indeed, it is Rieffel's category \emph{Normod-M} \cite{Rie74m}. We won't presume that the reader is familiar with $W^*$-category theory. Instead, we shall focus on Rieffel's paper for what we need, leaving $W^*$-categories for another time.

\subsubsection{An important Lemma}
The following lemma is one of the pillars of our analysis. It allows us to "move around" among intertwiner spaces and it underwrites use of use of Morita equivalence. Of most immediate importance it enables us to define the domains of our noncommutative functions.

\begin{lemma}\label{intertwiners}
Let $\sigma, \tau$, and $\rho$ be three representations in $\widehat{M}$. If the kernel of $\rho$ is contained in the kernels of both $\tau$ and $\sigma$, then
\begin{enumerate}
    \item the containment \eqref{eq:maps}
can be promoted to the equation
\begin{equation}\label{eq:key} 
\cI(\tau,\sigma)=\overline{\text{span}(\cI(\rho,\sigma)\cdot\cI(\tau,\rho))}^{w^*}\footnote{Here $\text{span}(\cI(\rho,\sigma)\cdot\cI(\tau,\rho))$ denotes the linear span of 
$\cI(\rho,\sigma)\cdot\cI(\tau,\rho)$ 
and $\overline{\text{span}(\cI(\rho,\sigma)\cdot\cI(\tau,\rho))}^{w^*}$ denotes its $weak^*$ closure.
},
\end{equation} and 
\item
a positive $R\in \sigma(M)'$  may be written
\be
R= \sum_{i=1}^{\infty} X_i^*X_i
\ee where each $X_i$ lies in $\cI(\sigma,\rho)$ and the series converges in the ultraweak topology on $B(H_{\sigma})$.
\end{enumerate}
\end{lemma}

It is easy to give examples where $1$ or $2$ fails if $\ker\rho$ is not contained in $\ker\sigma \cap \ker\tau$.

\begin{proof}

Form the representation 
\begin{equation*}
    \pi:= \sigma\oplus \rho \oplus \tau
\end{equation*}
of $M$ acting on $H_{\pi}:= H_{\sigma} \oplus H_{\rho} \oplus H_{\tau}$ and let $\sP_{\sigma}$, $\sP_{\rho}$ and $\sP_{\tau}$ denote the projections of $H_{\pi}$ onto its first, second, and third summands, respectively. 
These projections are in $\pi(M)'$ and their sum is the identity in $\pi(M)'$. 
By the definition of an intertwiner space, it is immediate that $\cI(\rho,\sigma)= \sP_{\sigma}\pi(M)'\sP_{\rho}$, $\cI(\tau,\rho)= \sP_{\rho}\pi(M)'\sP_{\tau}$, and $\cI(\tau,\sigma)= \sP_{\sigma}\pi(M)'\sP_{\tau}$. 
Let $V$ be the linear span in $\pi(M)'$ of the set $\sP_{\sigma}\pi(M)'\sP_{\rho}\pi(M)'\sP_{\tau}$. 
The argument for \eqref{eq:maps} shows that $V\subseteq \sP_{\sigma}\pi(M)'\sP_{\tau}$. 
But also, $V$ may be written as $\sP_{\sigma}(\text{span}(\pi(M)'\sP_{\rho}\pi(M)'))\sP_{\tau}$. 
However, $\fI:=\text{span}(\pi(M)'\sP_{\rho}\pi(M)')$ is the two-sided ideal in $\pi(M)'$ generated by $\sP_{\rho}$.  
So, $\overline{V}^{w^*} = \sP_{\sigma} \overline{\fI}^{w^*} \sP_{\tau}$ and $\overline{\fI}^{w^*}$ is the $w^*$-closed ideal in $\pi(M)'$ generated by $\sP_{\rho}$. 
By Corollary 3 in \cite[I.3]{Dix} a $w^*$-closed ideal in a von Neumann algebra is generated by a central projection. 
So, 
\begin{equation}
    \overline{V}^{w^*}=\sP_{\sigma} \pi(M)'\fp \sP_{\tau}
\end{equation} 
for a projection $\fp$ in the center of $\pi(M)'$. 
Since the center of $\pi(M)'$ is the center of $\pi(M)$ which, in turn, is the image of the center of $M$ under $\pi$, we may write $\fp = \pi(p_{\pi})$ for a projection $p_{\pi}$ in the center of $M$. 
This yields the equation
\begin{equation}
\fp = \pi(p_{\pi})= \sigma(p_{\pi})\oplus \rho(p_{\pi})\oplus \tau(p_{\pi}),
\end{equation}
which, in turn, implies that the projections $\sigma(p_{\pi})$, $\rho(p_{\pi})$, and $\tau(p_{\pi})$ are pairwise orthogonal in $\pi(M)'$. 
Since $\sP_{\rho} \le \pi(p_{\pi})$,  $I_{H_{\rho}} \le \rho(p_{\pi})$, which implies $\rho(1_M - p_{\pi}) = 0$. 
But then, by hypothesis, $\sigma(1_M - p_{\pi})=\tau(1_M - p_{\pi})=0$. Hence $\pi(1_M-p_{\pi})=0$
, so $\fp = \pi(1_M)=I_{H_{\pi}}$, which proves equation \eqref{eq:key}.

For the proof of the second assertion of the lemma, first note that if we let $\sigma = \tau$, the first part shows that 
\begin{equation}
\sigma(M)'=\cI(\sigma,\sigma)=\overline{\cI(\rho,\sigma)\cdot\cI(\sigma,\rho)}^{w^*}=\overline{\cI(\sigma,\rho)^*\cdot\cI(\sigma,\rho)}^{w^*},
\end{equation}
which shows that \emph{any} element $\sigma(M)'$ can be written as an ultraweakly convergent \emph{net of sums} of pairs of products taken from $\cI(\sigma,\rho)^*\cdot \cI(\sigma,\rho)$. 
Our problem, then, is to show that if $R$ is a positive element in $\sigma(M)'$, then it can be written an ultraweakly convergent \emph{countable series} of the indicated form.

Consider the direct sum representation $\pi: = \sigma \oplus \rho$ of $M$ on the Hilbert space $H_{\pi}: = H_{\sigma} \oplus H_{\rho}$. Let $\sP_{\sigma}$ and $\sP_{\rho}$ in $B(H_{\pi})$ be the projections onto the first and second summands, respectively. Since $\sP_{\sigma}$ and $\sP_{\rho}$ commute with $\pi(M)$, they belong to the commutant $\pi(M)'$. As established previously, we can identify the intertwiner space $\cI(\sigma, \rho)$ with the subspace $\sP_{\rho} \pi(M)' \sP_{\sigma}$ of the commutant. 
Let $c(\sP_{\rho})$ denote the central support of $\sP_{\rho}$ in $\pi(M)'$. Since the center of $\pi(M)'$ is the image under $\pi$ of the center of $M$ it follows that
$c(\sP_{\rho}) = \pi(p)$ for some central projection $p \in M$. The inequality $\sP_{\rho} \le \pi(p)$ restricts to $I_{H_{\rho}} \le \rho(p)$ on the second summand, implying that $c(\sP_{\rho}) = I_{H_{\pi}}$. This, in turn, shows that the central support of $\sP_{\rho}$ dominates the central support of $\sP_{\sigma}$.
Now we may apply the standard comparison theory of projections in von Neumann algeras found in Sections 1 and 2 of \cite[
I.2]{Dix}. In particular, \cite[I.2.4. Theorem 3]{Dix} implies $\sP_{\sigma}$ is equivalent to a projection $F$ such that $F \le \sP_{\rho} \otimes I_{\ell^2}$ within the amplification $\pi(M)' \otimes B(\ell^2)$.
Consequently, there exists a sequence of partial isometries $\{V_i\}_{i=1}^\infty$ in $\pi(M)'$ such that:
\begin{enumerate}
    \item $V_i^* V_i \le \sP_{\sigma}$ for all $i$, and the sum, $\sum_{i=1}^\infty V_i^* V_i = \sP_{\sigma}$, converges in the ultraweak topology,
while    
\item $V_i V_i^* \le \sP_{\rho}$ for all $i$, and the projections $V_i V_i^*$ are mutually orthogonal.
\end{enumerate}
Because $V_i = \sP_{\rho} V_i \sP_{\sigma} \in \sP_{\rho} \pi(M)' \sP_{\sigma}$, each $V_i$ naturally identifies with an intertwiner in $\cI(\sigma, \rho)$. 

Now, let $R \in \sigma(M)'$ be a positive operator. There exists a unique positive square root $R^{1/2} \in \sigma(M)'$. 
Define $X_i = V_i R^{1/2}$. 
Since $V_i \in \cI(\sigma, \rho)$ and $R^{1/2} \in \cI(\sigma, \sigma) = \sigma(M)'$, their composition satisfies:
\[
X_i \sigma(a) = V_i R^{1/2} \sigma(a) = V_i \sigma(a) R^{1/2} = \rho(a) V_i R^{1/2} = \rho(a) X_i \quad \text{for all } a \in M.
\]
Thus, $X_i \in \cI(\sigma, \rho)$ for every $i$.
Finally, we may compute the ultraweakly convergent sum,
\begin{multline*}
\sum_{i=1}^\infty X_i^* X_i = \sum_{i=1}^\infty (V_i R^{1/2})^* (V_i R^{1/2}) 
= \sum_{i=1}^\infty R^{1/2} V_i^* V_i R^{1/2} 
= R^{1/2} \left( \sum_{i=1}^\infty V_i^* V_i \right) R^{1/2} \\
= R^{1/2} \sP_{\sigma} R^{1/2} 
= R^{1/2} I_H R^{1/2} 
= R,
\end{multline*}
to complete the proof.
\end{proof}

\subsection{Induced Modules: The soil in which our functions grow}\label{ssect:Indmod}

Throughout, $E$ will denote a $W^*$-correspondence over $M$ in the sense of \cite[Definition 2.2]{MS02}. 
This means, in particular, that $E$ is a bimodule over $M$ on which we denote the right action of $M$ simply with a $\verb|\cdot|$ while we write the left action using a normal homomorphism $\varphi_E$ from $M$ into the $W^*$-algebra of ultraweakly continuous $M$-module maps of $E$, $\cL(E)$.

We emphasize that $E$ need not be related to any particular Hilbert space representations of $M$. 
In fact, we find it preferable to view the pair $E$ and $\varphi_E$ as forming a map from $\widehat{M}$ into $\widehat{M}$ in the following fashion.

A $\sigma\in \widehat{M}$, acting on the Hilbert space $H_\sigma$, may be promoted to a representation $\sigma^E$ of $\cL(E)$ acting on the Hilbert space $E\otimes_{M} H_{\sigma}$ -- the tensor product of $E$ with $H_{\sigma}$ that is balanced over the actions of $M$; i.e., 
\begin{equation*}
    (\xi\cdot a)\otimes h=\xi\otimes \sigma(a)h\qquad \xi\in E, h\in H_{\sigma}, a\in M.
\end{equation*}The promoted representation, $\sigma^E$ of $\cL(E)$, is defined by the formula
\begin{equation}\label{eq:inducedrep}
    \sigma^E(T):= T\otimes I_{H_{\sigma}}, \qquad T\in \cL(E).
\end{equation}
It is evident that $\sigma^E$ is a normal representation of $\cL(E)$ acting on $E\otimes_M H_{\sigma}$ and so one may compose it with $\varphi_E$ to obtain a normal representation of $M$ on $E\otimes_M H_{\sigma}$. 
The composition $\sigma^E\circ\varphi_E $ is known in the literature as the representation of $M$ that is \emph{induced} by $\sigma^E$ and $\varphi_E$, and is often written $Ind_{\varphi_E}^{\sigma_E}$, cf. \cite{RW1998}.

It is not the representation $Ind_{\varphi_E}^{\sigma_E}$ by itself that plays a role in our theory. Rather, it is the \emph{intertwiner space},
\begin{equation*}
    \cI(Ind_{\varphi_E}^{\sigma_E},\sigma),
\end{equation*}
that is of central importance. When it is necessary to include $E$ and $\varphi_E$ in addition to $\sigma$, we shall write
    $\sigma^E\circ \varphi$
for $\sigma^E\circ \varphi_E$.
Usually in our discussions, $E$ and $\varphi$ are fixed in advance and $\sigma$ is the only variable. In that event, we write $\h{\sigma}$ for $Ind_{\varphi_E}^{\sigma_E}$, which leads to 
$\cI(\h{\sigma},\sigma)$. Thus we use multiple notation for the same space:
\begin{equation}\label{eq:domnotation}
   \cI(Ind_{\varphi_E}^{\sigma_E},\sigma)= \cI(\sigma^E\circ\varphi,\sigma)=\cI(\h{\sigma},\sigma).
\end{equation}

\begin{definition}\label{EsIntspac}
    The space $\cI(\h{\sigma},\sigma)$ is called the $(E,\sigma)$-\emph{intertwiner space} determined by $E$ and $\sigma$.
\end{definition}

This space is an \emph{operator space} that is contained in
\begin{equation*}
B(H_{\h{\sigma}},H_{\sigma}),
\end{equation*}
where $H_{\h{\sigma}}=E\otimes_M H_{\sigma}$. It is closed in all the operator topologies on $B(H_{\h{\sigma}},H_{\sigma})$ and it is a right $\h{\sigma}(M)'$ - left $\sigma(M)'$ operator bimodule. 

The domains of the matricial functions we study are contained in $\cI(\h{\sigma},\sigma)$. We are inspired by the poetic assertion of Herman Weyl from the introduction to his book on Riemann surfaces \cite{Weyl64} that Riemann surfaces are the ``native land'' for analytic functions, they are the ``soil
in which the functions grow and
thrive''. The soil in which our functions grow is $\cI(\h{\sigma},\sigma)$ and to emphasize this, we write its elements using lower case Fraktur letters from the end of the alphabet: 
\begin{equation}
    \cI(\h{\sigma},\sigma)=\{\fz\in B(H_{\h{\sigma}},H_{\sigma}):\fz\h{\sigma}(a)=\sigma(a)\fz, a\in M\}.
\end{equation}
(In alternative, extended notation:
\begin{equation*}
\cI(\h{\sigma},
    \sigma)=\mathcal{I}(\sigma^E \circ \varphi, \sigma)=\{\fz:E\otimes_M H_{\sigma} \rightarrow H_{\sigma} : \; \fz(\sigma^E \circ \varphi_E)(a)=\sigma(a) \fz , \; a\in M\}.)
    \end{equation*}

It may be helpful to note that in listed references the domains of our functions are viewed as contained in $\cI(\sigma,\h{\sigma})$ as opposed to $\cI(\h{\sigma},\sigma)$. The change to the current convention found in subsequent listed references was made to avoid the uncontrollable shower of asterisks that were the result of the first convention.  

The next lemma is a corollary of Lemma \ref{intertwiners} and shows  how the domains change when one changes representations. It also shows that the intertwiner spaces act on themselves in an an analogue of the quadratic representation in the theory of Jordan algebras.

\begin{lemma}\label{corintertwiners}
Suppose $\sigma, \rho\in \widehat{M}$ with $\ker(\rho)\subseteq \ker(\sigma)$. Then
the $w^*$-closure of 
$$\text{span}\{aT | a\in \cI(\rho,\sigma), \; T\in B(H_{\rho}) \}: 
= \overline{\cI(\rho,\sigma)B(H_{\rho})}^{w^*}$$
is $B(H_{\rho},H_{\sigma})$, and 
\begin{multline}\label{eq:transform1}
\overline{(I_E\otimes\cI(\rho,\sigma))\cI(\rho,\h{\rho})\cI(\sigma,\rho)}^{w^*}
=\overline{(I_E\otimes\cI(\rho,\sigma))\cI(\rho,\rho\circ\varphi_E)  \cI(\sigma,\rho)}^{w^*}\\
=\cI(\sigma, \sigma^E\circ\varphi_E)=\cI(\sigma,\h{\sigma}).
\end{multline} 
while 
\begin{multline}\label{eq:transform2}
\overline{\cI(\rho,\sigma)\cI(\h{\rho},\rho)(I_E\otimes\cI(\sigma,\rho))}^{w^*}=\overline{(I_E\otimes \cI(\rho,\sigma))\cI(\rho,\rho^E \circ \varphi_E )\cI(\sigma,\rho)}^{w^*}\\=\cI( \sigma^E \circ \varphi_E,\sigma )=\cI(\h{\sigma},\sigma).\end{multline}
\end{lemma}
\begin{proof} 
Use Lemma~\ref{intertwiners} to write $I_{H_{\sigma}}=w^*-\lim \sum_{i=1}^n a_i a_i^*$,  $a_i\in \cI(\rho,\sigma)$). Then,
given $A\in B(H_{\rho},H_{\sigma})$, $A=w^*-\lim \sum_{i=1}^n a_i (a_i^*A)$. So $A\in  \overline{\cI(\rho,\sigma)B(H_{\rho})}^{w^*}$. A $\fw\in \cI( \sigma, \sigma^E \circ \varphi_E)$ may be written as the $w^*$ limit: $w^*-lim_{n,k} (I_E\otimes \sum_{i=1}^n a_ia_i^*)\fw (\sum_{j=1}^k a_ja_j^*)$. Since $(I_E\otimes a_i^*)\fw a_j \in \cI( \rho, \rho^E \circ \varphi_E)$, $a_j\in \cI(\rho,\sigma)$ and $a_i^*\in \cI(\sigma,\rho)$, the proof of \eqref{eq:transform1} is complete. The the proof of \eqref{eq:transform2} follows by taking adjoints as described in Subsection \ref{sss:SumsRepInt}.
\end{proof}

\subsection{Matricial Calculus}\label{ssect:matcalc}
We turn now to an operator space version of the matricial structures studied in the NC\nobreakdash-function theory literature. We have tried to formulate terminology which highlights the similarities between our work and NC\nobreakdash-function theory and reveals the significant differences. The differences result from the fact that in NC\nobreakdash-function theory the structures involved are fibred over the nonnegative integers while in our theory the structures are fibred over the normal representations of a $W^*$-algebra. As we intend to show, with careful specializations, many NC\nobreakdash-function theoretic results can be derived from our theory.

\subsubsection{Sums of intertwiner spaces and matricial families}
It is important for this work to understand how to sum intertwiner spaces.  
Algebraically, there is no issue, we just forms the usual external direct sum.  
But an intertwiner space is also an operator space and the question arises about how to give the direct sum an operator space structure.  
It turns out that there are many possibilities, but for our purposes there are only two: the $\ell^{\infty}$-norm and the column norm.  
The $\ell^{\infty}$-norm is the norm of choice when studying sums of intertwiner spaces. 
The column norm arises when studying Hilbert space representations of our algebras.
We begin with some of the main points about the $\ell^{\infty}$ here citing Gilles Pisier's book \cite{Pis2003} for the basics of operator space theory.

To understand what the $\ell^{\infty}$-norm means for intertwiner spaces, we begin with four Hilbert spaces, $H_1$, $H_2$, $H_3$ and $H_4$. We then partition $B(H_1\oplus H_2 \oplus H_3 \oplus H_4)$ into blocks to reveal the formula
\begin{multline*}
    B(H_1\oplus H_2 \oplus H_3 \oplus H_4)\\\simeq\begin{bmatrix}
        B(H_1\oplus H_2,H_1\oplus H_2) & 
        B(H_1\oplus H_2,H_3\oplus H_4) 
        \\B(H_3\oplus H_4, H_1\oplus H_2) & B(H_3\oplus H_4, H_3\oplus H_4)
    \end{bmatrix}\\
    \simeq \begin{bmatrix}
        B(H_1,H_1)&  B(H_1,H_2)&  B(H_1,H_3)&  B(H_1,H_4)\\B(H_2,H_1)& B(H_2,H_2)& B(H_2,H_3)& B(H_2,H_4)\\
        B(H_3,H_1)& B(H_3,H_2)& B(H_3,H_3)& B(H_3,H_4)\\
        B(H_4,H_1)& B(H_4,H_2)& B(H_4,H_3)& B(H_4,H_4)
    \end{bmatrix},
\end{multline*}
indicating that the resulting matrix operator spaces are completely isometrically isomorphic. Indeed, the isomorphisms are all implemented by easily calculated unitary operators on the usual Hilbert space direct sum $H_1\oplus H_2 \oplus H_3\oplus H_4$. When we apply this observation to the Hilbert spaces $H_{\h{\sigma}}$, $H_{\sigma}$, $H_{\h{\tau}}$, and $H_{\tau}$, and focus on the intertwiner subspaces, we immediately arrive at the completely isometric isomorphism,

\begin{equation}\label{eq:odirectsum}
    \cI(\widehat{\sigma \oplus \tau},\sigma \oplus \tau) \simeq \begin{bmatrix}
        \cI(\widehat{\sigma},\sigma)& \cI(\widehat{\sigma},\tau)\\
        \cI(\widehat{\tau},\sigma)& \cI(\widehat{\tau},\tau)\end{bmatrix}.
\end{equation}

The operator space $B(H_1,H_1)\oplus_{\infty} B(H_2,H_2)$ is completely isometrically isomorphic to $\begin{bmatrix}
    B(H_1,H_1)&0\\
    0& B(H_2,H_2)
\end{bmatrix}$. So, the inclusion, $$\begin{bmatrix}
    B(H_1,H_1)&0\\
    0& B(H_2,H_2)
\end{bmatrix}\\
\hookrightarrow
\begin{bmatrix}
    B(H_1,H_1)&B(H_1,H_2)\\
    B(H_2,H_1)& B(H_2,H_2)
\end{bmatrix},$$
leads to the completely isometric embedding
\begin{equation}\label{eq:ospaceident}
B(H_1,H_1)\oplus_{\infty} B(H_2,H_2)
\hookrightarrow
\begin{bmatrix}
    B(H_1,H_1)&B(H_1,H_2)\\
    B(H_2,H_1)& B(H_2,H_2)
\end{bmatrix}.
\end{equation}
The identification, \eqref{eq:odirectsum}, and the formula, \eqref{eq:ospaceident}, then give precise meaning to the assertion that the containment
\begin{equation}\label{eq:containment}
    \cI(\widehat{\sigma},\sigma)\oplus_{\infty} \cI(\widehat{\tau},\tau) \subseteq
    \cI(\widehat{\sigma\oplus\tau},\sigma\oplus\tau)
\end{equation}
identifies $\cI(\widehat{\sigma},\sigma)\oplus_{\infty} \cI(\widehat{\tau},\tau)$ completely isometrically to a subspace of $\cI(\widehat{\sigma\oplus\tau},\sigma\oplus\tau)$, and, in turn, make the following definitions well-defined:
\begin{definitions}\label{def:dirsumsets} \begin{enumerate}
    \item If $\cU(\sigma) \subseteq \cI(\widehat{\sigma},\sigma)$ and $\cU(\tau) \subseteq
    \cI(\widehat{\tau},\tau)$ are \emph{subsets}\footnote{These are are not-necessarily
\emph{linear subspaces}.}, then the \emph{set} $\cU(\sigma)\oplusi \cU(\tau)$ is defined to be 
$\{S \oplus T\in \cI(\widehat{\sigma}, \sigma)\oplus_{\infty} \cI(\widehat{\tau}, \tau): 
S\in \cU(\sigma),T\in \cU(\tau)\}$.
    \item A family $\cU:=\{\cU(\sigma)\}_{\sigma\in \widehat{M}}$ 
where each $\cU(\sigma)$ is a subset of $\cI(\widehat{\sigma},\sigma)$ is said to be \emph{closed under direct sums} in case the containment \eqref{eq:containment}
carries $\cU(\sigma)\oplusi \cU(\tau)$ into $\cU(\sigma\oplus \tau)$ for every $\sigma$ and $\tau$. 
\item If $\cU$ closed under direct sums, then $\cU$ is called a \emph{matricial family of sets} indexed by $\widehat{M}$.
\item The full family $\{\cI(\widehat{\sigma},\sigma)\}_{\sigma\in \widehat{M}}$ is denoted simply by an unalloyed $\cI$.
\end{enumerate}
\end{definitions}

Perhaps the simplest example of a matricial family of sets for our purposes is a matricial disc.
\begin{definition}\label{def:mudisc}
    The \emph{matricial unit disc centered at the origin} is $\bD:=\{\bD(\sigma)\}_{\sigma\in \widehat{M}}$ where $\bD(\sigma):=\{C\in \cI(\widehat{\sigma},\sigma): \Vert C \Vert < 1\}$.
\end{definition}
The fact that $\bD(\sigma)\oplusi \bD(\tau)\subseteq \bD(\sigma \oplus \tau)$ is an immediate reflection of the definition of the $\ell_{\infty}$ direct sum norm. Clearly, one may change the radius to any positive number to get other matricial families.  Changing the center is more complicated and will not be pursued here.

However, we want to point out some variations of the matricial unit disc that play prominent roles in our theory. They will be discussed further in Section\,\ref{example}. 

For this purpose, note that one may form the tensor powers of $E$ that are balanced over $M$ through the representation $\varphi$ of $M$ into $\cL(E)$, denoted $E^{\otimes n}$. The representation $\varphi$, then, may be promoted to a representation of $M$ on $E^{\otimes n}$, $\varphi_n$, via the formula
\begin{equation*}
    \varphi_n(a)(\xi_1\otimes \xi_2 \otimes \cdots \otimes \xi_n): = (\varphi(a)\xi_1)\otimes \xi_2 \otimes \cdots \otimes \xi_n,\qquad a\in M.
\end{equation*}
On the other hand, given $\fz\in \cI(\widehat{\sigma},\sigma)$ we may form the map $\fz^{(n)}$ from $E^{\otimes n}\otimes H_{\sigma} $ to $H_{\sigma}$ defined inductively by the formula
\begin{equation}\label{eq:power}
\fz^{(n)}:= \fz(I_E \otimes \fz)(I_{E^{\otimes 2}}\otimes \fz)\cdots (I_{E^{\otimes (n-1)}}\otimes \fz),
    \end{equation}
which we call the $n^{\text{th}}$ \emph{tensorial power} of $\fz$, $n\geq 1$. Clearly, $\Vert \fz^{(n)}\Vert \leq \Vert \fz \Vert^n$.

If we sum the tensor powers of $E$ we obtain the \emph{Fock space}, 
$\cF(E):=\sum_{n=0}^{\infty}E^{\otimes n}$, together with the representation of $M$,
$\varphi_{\infty}:= \text{diag}\,{(\varphi_0, \varphi_1, \varphi_2,  \cdots)}.$ Here, $E^{\otimes 0}$ is $M$ and $\varphi_0$ is left multiplication.

Now let $\{X_k\}_{k=1}^{\infty}$ be a sequence of positive semidefinite operators, where $X_k$  lies in the \emph{commutant} of $\varphi_k(M)$ in $\cL(E^{\otimes k})$. We assume, too, that $X_1$ is invertible and that $\limsup{\Vert X_k \Vert^{\frac{1}{k}}}<\infty$. In that event one can show \cite[Lemma 4.2]{MS16} that for $\fz\in \fI(\widehat{\sigma},\sigma)$ with $\Vert \fz \Vert<(\limsup{\Vert X_k \Vert^{\frac{1}{k}}})^{-\frac{1}{2}}$ the series
\begin{equation*}
    \sum_{k=1}^{\infty} \Vert \fz^{(k)}(X_k\otimes I_{H_k})\fz^{(k)*} \Vert
\end{equation*}
converges. For such sequences $X$ we define:

\begin{definition}\label{def:wmudisc}
    The \emph{weighted matricial unit disc determined by $X$ and  centered at the origin} is $\bD_X:=\{\bD(X,\sigma)\}_{\sigma\in \widehat{M}}$ where 
    \begin{equation*}
        \bD(X,\sigma):=\{\fz\in \cI(\widehat{\sigma},\sigma): \Vert \sum_{k=1}^{\infty}\fz^{(k)}(X_k\otimes I_{H_{\sigma}})\fz^{(k)*} \Vert < 1\}.
    \end{equation*}
\end{definition}

Observe that the hypotheses that $X_1$ is invertible guarantees that each $\bD(X,\sigma)$ contains a neighborhood of the origin in $\cI(\widehat{\sigma},\sigma)$. Observe, too, that if $X_1=I$ and all the other $X_k$ are zero, then $\bD_X = \bD$.

Another source of examples comes from the spaces of ``\emph{absolutely continuous points}'' contained in a $\bD_X$. These are defined in \cite{MS11} but do not play a role in this paper.

\begin{definition}\label{def:GlobMat}
    Suppose $\cA=\{\cA(\sigma)\}_{\sigma \in \widehat{M}}$ and $\mathcal{V}=\{\cV(\sigma)\}_{\sigma \in \widehat{M}}$ are two matricial families of sets indexed by $\widehat{M}$  and suppose that for each $\sigma \in \widehat{M}$, $f_{\sigma}:\cA(\sigma) \rightarrow \cV(\sigma)$ is a function. Then the family $f:=\{f_{\sigma}\}_{\sigma\in \widehat{M}}$ is called
    \begin{enumerate}
        \item \emph{a global family of functions} in case $f$ \emph{respects direct sums} in the sense that\begin{equation}
        f_{\sigma\oplus \tau}(\fz\oplus \fw)=f_{\sigma}(\fz)\oplus f_{\tau}(\fw),\qquad \fz\oplus \fw \in \cA(\sigma\oplus \tau),\end{equation} and
        \item a \emph{matricial family} in case $f$ \emph{respects intertwiners} in the sense that 
    \end{enumerate}

    \[ Cf_{\sigma}(\mathfrak{z})=f_{\tau}(\mathfrak{w})C  \]
    for every $\mathfrak{z}\in \cA(\sigma)$, every $\mathfrak{w}\in \cA(\tau)$ and every $C\in \cI(\sigma,\tau)$ such that \[ C\mathfrak{z}=\mathfrak{w}C.\]
\end{definition}
The reason we mention these notions here is that the two conditions together lead to our tensorial function theory which is a natural generalization of Joe Taylor's approach to noncommutative holomorphic functions \cite{MS13}. So, when our kernels (to be defined) have these properties it will follow that the spaces of functions they generate consist of (tensorial) holomorphic functions. We will track these conditions in our definitions to follow.

\subsection{Kernels}
Throughout this subsection, $\cU=\{\cU(\sigma)\}_{\sigma \in \widehat{M}}$ will denote a matricial family of subsets of $\cI=\{\cI(\widehat{\sigma},\sigma)\}_{\sigma\in \widehat{M}}$ in the sense of Definition \ref{def:dirsumsets}, i.e. $\cU(\sigma)\subseteq \cI(\widehat{\sigma},\sigma)$. Also, if $\sK$ and $\sL$ are dual operator spaces, then $B_{*}(\sK,\sL)$ denotes the space of bounded, ultraweakly continuous, linear maps from $\sK$ to $\sL$.

Recall that the matricial unit disc centered at the origin \ref{def:wmudisc} is the prototypical example of a matricial family $\cU$.

The following definitions and terminology are modeled on definitions found in Section 2.3 of \cite{BMV1}.

\begin{definition}\label{kernel}
Let $\cU=\{\cU(\sigma)\}_{\sigma \in \widehat{M}}$ be a matricial family of sets. 
A map $K$ from $\cU\times\cU$ to  $\coprod_{\sigma,\tau \in \widehat{M}} B_{*}(\cI(\tau,\sigma),B(H_{\tau},H_{\sigma}))
$,
\begin{equation}
    K:\cU\times \cU \to \coprod_{\sigma,\tau \in \widehat{M}} B_{*}(\cI(\tau,\sigma),B(H_{\tau},H_{\sigma})),
\end{equation}
is called a $\coprod_{\sigma,\tau \in \widehat{M}}B_{*}(\cI(\tau,\sigma),B(H_{\tau},H_{\sigma}))$-\emph{valued kernel} on $\cU\times \cU$. 
Such a $K$ is called \emph{graded} with respect to $\cU$ (or by $\cU$) in case 
\begin{equation}\label{eq:switch}
    K(\fz,\fw )\in B_*(\cI(\tau,\sigma),B(H_{\tau},H_{\sigma})), \qquad (\fz,\fw)\in \cU(\sigma)\times \cU(\tau).
\end{equation}
\end{definition}

It is important to note how the variables switch in \eqref{eq:switch}: $(\fz,\fw)\in \cU(\sigma)\times \cU(\tau)$, but $K(\fz,\fw)$ is a map from $\cI(\tau,\sigma)$ to $B(H_{\tau},H_{\sigma})$. Further note that we may evaluate $K(\fz,\fw)$ on an operator $P\in \cI(\tau,\sigma)$ to obtain an operator $K(\fz,\fw)(P)$ in $B(H_{\tau},H_{\sigma})$. Now in future formulas, we will want to form expressions
\begin{enumerate}
    \item $K(\fz,\fw)(PA)\in B(H_{\tau},H_{\sigma})$ and
    \item $K(\fz,\fw)(PA)C \in B(H_{\rho},H_{\sigma})$.
\end{enumerate}
The first forces $PA$ to be in $\cI(\tau,\sigma)$. The second equation forces $C$ to have range in $H_{\tau}$. 

\begin{definition}\label{def:findirectsums}
A graded kernel $K$ is said to \emph{respect (finite) direct sums} in case for 
$\fz_i \in \cU(\sigma_i)$, $1\leq i\leq n$ , $\fw_j \in \cU(\tau_j)$ , $1\leq j\leq m$,

\begin{multline}
    K(\fz_1\oplusi \fz_2 \oplusi \cdots \fz_n, \fw_1\oplusi \fw_2\oplusi \cdots \fw_m)\left(\begin{bmatrix}
        C_{1,1} & C_{1,2} & \cdots & C_{1,m}\\
        C_{2,1} & C_{2,2} & \cdots & C_{2,m}\\
        \vdots  & \vdots  & \vdots & \vdots \\
        C_{n,1} & C_{n,2} & \cdots & C_{n,m}
    \end{bmatrix}\right)\\
    =
    \begin{bmatrix}
        K(\fz_1,\fw_1)(C_{1,1}) & K(\fz_1,\fw_2)(C_{1,2}) & \cdots & K(\fz_1,\fw_m)(C_{1,m})\\
        K(\fz_2,\fw_1)(C_{2,1}) & K(\fz_2,\fw_2)(C_{2,2}) & \cdots & K(\fz_2,\fw_m)(C_{2,m})\\
        \vdots  & \vdots  & \vdots & \vdots \\
        K(\fz_n,\fw_1)(C_{n,1}) & K(\fz_n,\fw_2)(C_{n,2}) & \cdots & K(\fz_n,\fw_m)(C_{n,m})
    \end{bmatrix}. 
\end{multline}
for all $C_{i,j}\in \cI(\tau_i,\sigma_j)$, $1\leq i \leq n, 1\leq j \leq m$. 
A graded kernel that respects finite direct sums will be called a \emph{global kernel}.
\end{definition}
It may be helpful to the reader who wants to compare our notation with that of \cite{BMV1} to note that our use of the $\ell^{\infty}$-direct sum, $\fz_1\oplusi \fz_2\oplusi \cdots \oplusi \fz_n$, would be written there as the block diagonal matrix
$\begin{bmatrix}
    \fz_1 & 0 & \dots & 0\\
    0 & \fz_2 & \dots & 0\\
    \vdots & \vdots & \ddots & \vdots\\
    0 & \dots & 0 & \fz_n
\end{bmatrix}$.
See \cite[pages 1848, 1849]{BMV1}.

The following is an example of a global kernel that is defined on a weighted matricial unit disc. It comes from \cite{MS16} and will be discussed further in Section~\ref{example}.

\begin{example}\label{WS}
As we shall see, it is possible to begin with a sequence $R:=\{R_k\}_{k\geq 1}$ of positive operators, where $R_k$ is in $\cL(E^{\otimes k})$ and commutes with $\varphi_k(M)$. The precise properties of $R$ are not important here except to say that from it one may construct a weight sequence $X$ that leads to the weighted matricial unit disc $\bD_X$ that is defined in Definition \ref{def:wmudisc}. Further, as we shall see, 
for $(\fz,\fw)\in \bD(X,\sigma)\times\bD(X,\tau)$ and $A\in \cI(\tau,\sigma)$, the formula
\begin{equation} 
K(\fz,\fw)(A)=\sum_{k=0}^{\infty} \fz^{(k)}(R_k^2 \otimes A)\fw^{(k) *}\in \cI(\tau,\sigma) \subseteq B(H_\tau, H_\sigma) 
\end{equation}
defines a global kernel that preserves intertwiner spaces.
This example is studied in detail in Section~\ref{example}.
\end{example}

\begin{definition}\label{def:matricialker}
A global kernel $K:\cU\times \cU\to \coprod_{\sigma,\tau\in \widehat{M}} B_{*}(\cI(\tau,\sigma),B(H_{\tau},H_{\sigma}))$ is called a \emph{matricial} kernel provided it respects intertwiners in the sense that, for $\fz_1 \in \cU(\sigma_1)$, $\fz_2 \in \cU(\sigma_2)$, $\fw_1 \in \cU(\tau_1)$ and $\fw_2 \in \cU(\tau_2)$, every $C\in \cI(\sigma_1,\sigma_2)$ such that $ C\fz_1=\fz_2(I_E\otimes C) $, every $D\in \cI(\tau_1,\tau_2)$ such that $ D\fw_1=\fw_2(I_E\otimes D)$ and every $P\in \cI(\tau_1,\sigma_1)$, we have 
\begin{equation}\label{kernelintertwining} 
CK(\fz_1,\fw_1)(P)D^*=K(\fz_2,\fw_2)(CPD^*).
\end{equation}
\end{definition}

\begin{definition}\label{cpkernel}

A kernel $K$ is called a \emph{completely positive} kernel (cp kernel) if for every $\sigma \in \widehat{M}$ and every $\fz \in \cU(\sigma)$, the map $K(\fz,\fz)$ is a positive map from $\sigma(M)'  (=\cI(\sigma,\sigma))$ to $B(H_{\sigma})(=B(H_\sigma,H_\sigma))$.
\end{definition}

\begin{remarks}\label{sum}
\begin{enumerate}
\item
 It follows from the definition of a global kernel (Definition \ref{kernel}) that if $K$ is cp, then for every $\sigma \in \widehat{M}$ and every $\fz \in \cU(\sigma)$, the map $K(\fz,\fz)$ is in fact a completely positive map that maps $\sigma(M)'$ into $B(H_{\sigma})$.
\item 
 A cp kernel \emph{intertwines adjoints}  in the sense that for every pair $(\fz,\fw)\in \cU(\sigma)\times \cU(\tau)$ and every $C\in \cI(\fz,\fw)$
 \begin{equation}
     K(\fz,\fw)(C^*)=(K(\fw,\fz)(C))^*.
 \end{equation}
 To see this, let $m>0$ be such that $\left( \begin{array}{cc} mI & C \\ C^* & mI \end{array}\right)$ is positive. Then using the assumption that $K$ is cp, conclude that 
 \begin{equation}0\leq K(\fz \oplusi \fw,\fz \oplusi \fw)(\left( \begin{array}{cc} mI & C \\ C^* & mI \end{array}\right)) = \left( \begin{array}{cc} mK(\fz,\fz)(I) & K(\fz,\fw)(C) \\K(\fw,\fz)( C^*) & mK(\fw,\fw)(I) \end{array}\right). 
 \end{equation}
 In particular, $( K(\fz,\fw)(C))^*=K(\fw,\fz)( C^*)$.
 \vspace{0.1cm}
 \item A cp (global) kernel is matricial (that is, it respects intertwiners  in the sense of Definition ~\ref{def:matricialker}) if and only if, for every $\fw\in \cU(\rho)$, the function
 $$k_{\fw}(\cdot):=K(\cdot,\fw),$$  satisfies the ``intertwining condition'': For all $C\in \cI(\sigma_1, \sigma_2)$ such that $C\fz_1 = \fz_2(I_E\otimes C)$, the following equality holds:
 $$Ck_\fw(\fz_1)=k_\fw(\fz_2)C.$$

 To see this, assume that $\{k_{\fw}\}$ all are matricial and compute, for $C,D,P$ as in (\ref{kernelintertwining}),
 \begingroup
\allowdisplaybreaks
  \begin{align*}
  	CK(\fz_1,\fw_1)(P)D^*&=Ck_{\fw_1}(\fz_1)(P)D^*=k_{\fw_1}(\fz_2)(CP)D^*=(D(k_{\fw_1}(\fz_2)(CP))^*)^*\\&=(D(K(\fz_2,\fw_1)(CP))^*)^*=(DK(\fw_1,\fz_2)(P^*C^*))^*\\&=(Dk_{\fz_2}(\fw_1)(P^*C^*))^*=(k_{\fz_2}(\fw_2)(DP^*C^*))^*\\&=(K(\fw_2,\fz_2)(DP^*C^*))^*=K(\fz_2,\fw_2)(CPD^*).
  \end{align*}
  \endgroup
  For the other direction, if $K$ satisfies (\ref{kernelintertwining}) and $C\in \mathcal{I}(\sigma_1,\sigma_2)$ is such that $C\fz_1=\fz_2(I_E\otimes C)$ then, for $B\in \mathcal{I}(\tau,\sigma)$,
  $$Ck_{\fw}(\fz_1)(B)=CK(\fz_1,\fw)(B)=K(\fz_2,\fw)(CB)=k_{\fw}(\fz_2)(CB).$$ Thus $Ck_{\fw}(\fz_1)=k_{\fw}(\fz_2)C$.
\end{enumerate}
\end{remarks}

\section{The Kolmogorov correspondence}\label{sec:Correspondence}
Suppose $\cU=\{\cU(\sigma)\}_{\sigma \in \widehat{M}}$ is a matricial family of sets for $(M,E)$ and $K:\cU \times \cU \rightarrow \coprod_{\sigma,\tau \in \widehat{M}} B_{*}(\cI(\tau,\sigma),B(H_{\tau},H_{\sigma}))$ is a cp kernel.
We fix a faithful $\rho\in \widehat{M}$; $\rho$ will be fixed for the entire discussion. We write $\cB(\sigma)$ for $B_*(\cI(\rho,\sigma),B(H_{\rho},H_{\sigma}))$ and $\cB$ for the family $\{\cB(\sigma)\}_{\sigma \in \widehat{M}}$. In order to conform with the notation and terminology of \cite[Theorem 3.1]{BMV1}, we create \emph{kernel elements} from $K$ as follows. For $(\fw,A,C) \in \cU(\tau)\times \cI(\tau,\rho)\times B(H_{\rho},H_{\tau})$ we define a \emph{kernel element} $K_{(\fw,A,C)}$ that maps $\cU(\sigma)$ into $\cB(\sigma)$  i.e., 
\begin{equation*}
K_{(\fw,A,C)}:\cU(\sigma)\to B_*(\cI(\rho,\sigma),B(H_{\rho},H_{\sigma})),
\end{equation*}
    by
\be\label{Kwac}  K_{(\fw,A,C)}(\fz)(B)=K(\fz,\fw)(BA)C  \ee for $\fz\in \cU(\sigma)$ and $B\in \cI(\rho,\sigma)$ .

Since $K$ is a global kernel, it respects finite direct sums in each variable by definition. From this it follows that the function elements $K_{(\fw,A,C)}$ are global functions in the sense of Definition \ref{def:GlobMat}. 
 
 The function $K_{(\fw,A,C)}$ is graded by definition and to see that $K$ respects intertwiners when $K$ is matricial we show that whenever 
 $D\in \cI(\sigma_1,\sigma_2)$ and satisfies $D\fz_1=\fz_2 (I_E\otimes D)$, then
 \be  DK_{(\fw,A,C)}(\fz_1)(B) = K_{(\fw,A,C)}(\fz_2)(DB) .\ee
 Indeed, using Definition~\ref{kernel} (3), we have $$K_{(\fw,A,C)}(\fz_2)(DB)=K(\fz_2,\fw)(DBA)C=DK(\fz_1,\fw)(BA)C=DK_{(\fw,A,C)}(\fz_1)(B).$$

\vspace{6mm}

We now define
\be E(K)^{(0)}:=\text{span} \{K_{(\fw, A,C)} : \fw \in \cU(\tau), A\in \cI(\tau,\rho), C\in B(H_{\rho},H_{\tau})\} .\ee

The elements of $E^{(0)}(K)$ are global functions from $\cU$ to $\cB$, that are also matricial functions if $K$  is a matricial kernel. Clearly $E^{(0)}(K)$ is a vector space under pointwise operations. We define a left action of $\rho(M)'$  and a right action of $B(H_{\rho})$ on this space, and a $B(H_{\rho})$-valued inner product so that the self-dual completion of $E^{(0)}(K)$, $E(K)$,  becomes a $\rho(M)',B(H_{\rho})-W^*$-correspondence.

The left action of $\rho(M)'$ will be written with $\varphi$ and given $f\in E^{(0)}(K)$, $\fz \in \cU(\sigma)$ and $B\in \cI(\rho,\sigma)$, we set
\be\label{leftaction}(\varphi(D)f)(\fz)(B)=f(\fz)(BD)
\ee for $D\in \rho(M)'$.
For $f=K_{(\fw,A,C)}$ we get
$$\varphi(D) K_{(\fw,A,C)}= K_{(\fw,DA,C)}.$$

The right action of $R \in B(H_{\rho})$ on $f\in E(K)^{(0)}$ is defined by
\be\label{rightaction}(f\cdot R)(\fz)(B)=f(\fz)(B)R \ee
and, thus, $$K_{(\fw,A,C)}\cdot R=K_{(\fw,A,CR)}.$$

The inner product of two generators is defined to be

\be\label{innerproduct}
 \langle K_{(\fw,A,C)},K_{(\fw',A',C')} \rangle =C^*K(\fw,\fw')(A^*A')C' \in B(H_{\rho})\ee and then extend to any pair of elements of $E(K)^{(0)}$ by sesquilinearity. Since $K$ is a cp kernel, we have $\langle f,f \rangle \geq 0$ for all $f\in E(K)^{(0)}$.
We also have \begin{align*}
	\langle K_{(\fw,A,C)},K_{(\fw',A',C')}\cdot R \rangle=\langle K_{(\fw,A,C)},K_{(\fw',A',C'R)} \rangle&=C^*K(\fw,\fw')(A^*A')C'R\\&=\langle K_{(\fw,A,C)},K_{(\fw',A',C')} \rangle R
\end{align*}
and
\begin{align*}
	\langle &K_{(\fw,A,C)},\varphi(D)K_{(\fw',A',C')} \rangle =\langle K_{(\fw,A,C)},K_{(\fw',DA',C')} \rangle =C^*K(\fw,\fw')(A^*DA')C'\\
	&=C^*K(\fw,\fw')((D^*A)^*A')C'=\langle K_{(\fw,D^*A,C)},K_{(\fw',A',C')} \rangle=\langle \varphi(D^*)K_{(\fw,A,C)},K_{(\fw',A',C')} \rangle
\end{align*}

We write $E(K)$ for the self dual completion of $E(K)^{(0)}$.

Observe that if $f=K_{(\fw',A',C')}$, then $$\langle K_{(\fw,A,C)},f \rangle =C^*K(\fw,\fw')(A^*A')C'=C^*f(\fw)(A^*).$$ 
Consequently, for every $f\in E(K)$,
\be\label{reproducing}
 \langle K_{(\fw,A,C)},f \rangle =C^*f(\fw)(A^*)
 \ee which is the reproducing property.
 
Suppose that $K$ is a matricial kernel. We want to prove that every function $f\in E(K)$ is matricial. We show that whenever $\fz\in \cU(\sigma)$, $\tilde{\fz}\in \cU(\tilde{\sigma})$, $B\in \cI(\rho,\sigma)$, $f\in E(K)$ and $C$ intertwines $\fz$ and $\tilde{\fz}$ (that is, $C\in \cI(\sigma,\tilde{\sigma})$ and $C\fz=\tilde{\fz}(I_E\otimes C)$), we have
\be\label{fmatricial} f(\tilde{\fz})(CB)=Cf(\fz)(B) .\ee
Note that since both sides of (\ref{fmatricial}) lie in $B(H_{\rho},H_{\sigma})$, it suffices to show that for every $D\in B(H_{\rho},H_{\sigma})$ we have \be\label{df}  D^*f(\tilde{\fz})(CB)=D^*Cf(\fz)(B) .\ee
For this, note first that it follows from the definition of a kernel (see (\ref{kernelintertwining})) that for every $E\in \cI(\rho,\tau)$ and $\fw \in \cU(\tau)$,
\[ K(\fw,\tilde{\fz})(EB^*)D=K(\fw,\fz)(EB^*)C^*D .\]
It follows that $K_{(\tilde{\fz},B^*C^*,D)}=K_{(\fz,B^*, C^*D)}$ and thus, using (\ref{reproducing}), $$D^*f(\tilde{\fz})(CB)=\langle K_{(\tilde{\fz},B^*C^*D)},f \rangle = \langle K_{(\fz,B^*, C^*D)}, f\rangle =D^*Cf(\fz)(B)$$ completing the proof of (\ref{fmatricial}).

Our discussion to this point proves the first implication in the statement of the following theorem which is analogous to Theorem 3.1 of \cite{BMV1}.  

\begin{theorem}\label{cpkerneltheorem}
Suppose $\cU=\{\cU(\sigma)\}_{\sigma \in \widehat{M}}$ is a matricial family of sets for $(M,E)$ and $K:\cU \times \cU \rightarrow \coprod_{\sigma,\tau \in \widehat{M}} B_{*}(\cI(\tau,\sigma),B(H_{\tau},H_{\sigma}))$ is a graded function. Then the following assertions are equivalent.
\begin{enumerate}
\item[(1)] $K$ is a [matricial ] cp kernel.
\vspace{0.1cm}
\item[(2)] Given a faithful $\rho \in \widehat{M}$, there is a $\rho(M)'$,$B(H_{\rho})
$-$W^*$-correspondence $E(K)$ whose elements are global [matricial] functions $f:\cU \rightarrow \cB$ such that
    \begin{enumerate}
    \item[(a)] $E(K)$ contains the functions $K_{(\fw,A,C)}$, for $\fw\in \cU(\tau)$, $A\in \cI(\tau,\rho)$ and $C\in B(H_{\rho},H_{\tau})$ defined in (\ref{Kwac}).
    \item[(b)] The functions $K_{(\fw,A,C)}$ have the reproducing property (\ref{reproducing}).
    \item[(c)] The left and right actions are as in (\ref{leftaction}) and (\ref{rightaction}).
\end{enumerate}
\vspace{0.1cm}
\item[(3)] (Kolmogorov decomposition)  Given a faithful $\rho \in \widehat{M}$, there is a $\rho(M)'$,$B(H_{\rho})
$ -$W^*$-correspondence $\mathcal{X}$ together with a matricial function $H : \cU \rightarrow \cL(\cI(\rho,\cdot)\otimes_{\rho(M)'}\mathcal{X}, B(H_{\rho},H_{\cdot}))$, where $\cI(\rho,\sigma)$ is a $\sigma(M)',\rho(M)'$ -$W^*$-correspondence, and $H(\fz):\cI(\rho,\sigma)\otimes_{\rho(M)'}\mathcal{X} \rightarrow B(H_{\rho},H_{\sigma})$ is a bounded adjointable operator for $\fz\in \cU(\sigma)$, such that for every $P\in \cI(\tau,\sigma)$, $C\in B(H_{\rho},H_{\tau})$, $\fz\in \cU(\sigma)$ and $\fw\in \cU(\tau)$,
\be\label{kolmogorov}
H(\fz)(P\otimes I_{\mathcal{X}})H(\fw)^*(C)=K(\fz,\fw)(P)C.
\ee
If $K$ is matricial, so is $H$ in the sense that, whenever $\fz\in \cU(\sigma)$, $\fw\in \cU(\tau)$ and $C\in \cI(\sigma,\tau)$ satisfies $C\fz=\fw(I\otimes C)$, then $CH(\fz)=H(\fw)(C\otimes I_{\mathcal{X}})$.

\end{enumerate}

\end{theorem}
\begin{proof}
It was shown above that (1) implies (2). To prove that (2) implies (3), we set $\mathcal{X}=E(K)$ and, for $\fz\in \cU(\sigma)$, we define
$$ H(\fz): \cI(\rho,\sigma)\otimes_{\rho(M)'}E(K) \rightarrow B(H_{\rho},H_{\sigma}) $$ by
\be
H(\fz)(A\otimes f):=f(\fz)(A)
\ee 
for $A\in \cI(\rho,\sigma)$ and $f\in E(K)$.
We first claim that $H$ is well defined and bounded. For this, fix an element $\sum_{i=1}^n A_i \otimes f_i \in \cI(\rho,\sigma)\otimes_{\rho(M)'}E(K)$ and write $A=(A_1, A_2, \ldots , A_n)\in B(H_{\rho}\oplus H_{\rho} \oplus \cdots \oplus H_{\rho}, H_{\sigma})$. Write $A=V|A|$ for the polar decomposition of $A$. Since each $A_i$ belongs to $\cI(\rho,\sigma)$, $|A|\in M_n(\rho(M)')$ and $V$ is the unique partial isometry in $B(H_{\rho}\oplus H_{\rho} \oplus \cdots \oplus H_{\rho}, H_{\sigma})$ with $\ker(V)=\ker(A)$. We write $V=(V_1,V_2,\ldots V_n)$ where $V_i\in B(H_{\rho},H_{\sigma})$. Note that $\ker(A)\subseteq H_{\rho}\oplus H_{\rho} \oplus \cdots \oplus H_{\rho}$ is invariant under $B\oplus B \oplus \cdots \oplus B$ for $B\in \rho(M)$ and, thus, if $U$ is a unitary operator in $M$ and if we set $V_U=(\sigma(U^*)V_1\rho(U),\ldots, \sigma(U^*)V_n\rho(U))$ then $V_U$ is a partial isometry satisfying $A=V_U|A|$ and $\ker(V_U)=\ker(V)=\ker(A)$. It follows from uniqueness that $V=V_U$ and, therefore, $V_i=\sigma(U^*)V_i\rho(U)$ for every $i$. Since this holds for every unitary $U\in M$, $V_i\in \cI(\rho,\sigma)$ for all $i$.

For every $C\in B(H_{\rho},H_{\sigma})$ we have, using (\ref{reproducing}), $$\sum_i C^*f_i(\fz)(A_i)=\sum_i \langle K_{(\fz,A_i^*,C)}, f_i\rangle$$ and \begin{align*}
\big\|\sum_i f_i(\fz)(A_i)\big\|&=\sup_{||C||=1} \big\|\sum_i \langle K_{(\fz,A_i^*,C)}, f_i\rangle\big\|\\&=\sup_{||C||=1}||\langle (K_{(\fz,A_1^*,C)},\ldots, K_{\fz,A_n^*,C})^t,(f_1,\ldots,f_n)^t\rangle||
\end{align*} where the last inner product is in $E(K)\oplus E(K) \oplus \cdots \oplus E(K)$.

Since $A=V|A|$ we can write $|A|=(A_{i,j})\in M_n(\rho(M)')$ and get $A_i^*=\sum_j A_{i,j}V_j^*$. Thus $K_{(\fz,A_i^*,C)}=\sum_j K_{(\fz,A_{i,j}V_j^*,C)}=\sum_j \varphi(A_{i,j})K_{(\fz,V_i^*,C)}$ (using (\ref{leftaction})). Writing $\varphi_n(|A|)=(\varphi(A_{i,j}))$, we have $(K_{(\fz,A_1^*,C)},\ldots, K_{(\fz,A_n^*,C)})^t=\varphi_n(|A|)(K_{(\fz,V_1^*,C)},\ldots, K_{\fz,V_n^*,C)})^t$. Since $\varphi_n(|A|)$ is selfadjoint, we get
\begingroup
\allowdisplaybreaks
\begin{align*}
	\big\|\sum_i f_i(\fz)(A_i)\big\|&=\sup_{\|C\|=1}\|\langle (K_{(\fz,V_1^*,C)},\ldots, K_{(\fz,V_n^*,C)})^t,\varphi_n(|A|)(f_1,\ldots,f_n)^t\rangle\|\\
	&\leq \sup_{\|C\|=1}\|(K_{(\fz,V_1^*,C)},\ldots, K_{(\fz,V_n^*,C)})^t\| \|\varphi_n(|A|)(f_1,\ldots,f_n)^t\|.
\end{align*}
\endgroup
Now, 
\begingroup
\allowdisplaybreaks
\begin{align*}
	\|\varphi_n(|A|)&(f_1,\ldots,f_n)^t\|^2=\|\langle \varphi_n(|A|)(f_1,\ldots,f_n)^t,\varphi_n(|A|)(f_1,\ldots,f_n)^t\rangle\|\\&=\|\langle (f_1,\ldots,f_n)^t,\varphi_n(A^*A)(f_1,\ldots,f_n)^t\rangle\|=\|\sum_{i,j}\langle f_i, \varphi(A_i^*A_j)f_j \rangle\|\\
	&=\|\sum_{i,j}\langle A_i \otimes f_i, A_j \otimes f_j \rangle\|=\|\langle \sum_i A_i\otimes f_i,\sum_j A_j \otimes f_j \rangle\|\\
	&=\|\sum_i A_i\otimes f_i\|^2
\end{align*}
\endgroup
Also,
\begin{align*}
	\|(K_{(\fz,V_1^*,C)},\ldots, K_{(\fz,V_n^*,C)}^t\|^2 &=\|\langle (K_{(\fz,V_1^*,C)},\ldots, K_{(\fz,V_n^*,C)}^t,(K_{(\fz,V_1^*,C)},\ldots, K_{(\fz,V_n^*,C)})^t\rangle\|\\&=||\sum_i \langle K_{(\fz,V_i^*,C)},K_{\fz,V_i^*,C}\rangle||=||C^*K(\fz,\fz)(\sum_i V_iV_i^*)C\|\end{align*} 
	(using the definition of the inner product in $E(K)$). Since $\sum_i V_iV_i^*=VV^*\leq I$, we have
$$||(K_{(\fz,V_1^*,C)},\ldots, K_{(\fz,V_n^*,C}))^t||\leq ||C||||K(\fz,\fz)||^{1/2}.$$

Combining the computations above we get
$$||\sum_i f_i(\fz)(A_i)||\leq ||\sum_i A_i\otimes f_i||||K(\fz,\fz)||^{1/2}.$$
Thus $H(\fz) $ is bounded. In fact $||H(\fz)||\leq ||K(\fz,\fz)||^{1/2}$.

Since, for $D\in \rho(M)'$, $\sum_i f_i(\fz)(A_iD)=\sum_i \varphi(D)f_i(A_i)$ it follows that $H(\fz)$ is well defined. This completes the proof of the claim above.

We still have to prove that $H(\fz)$ is adjointable and, for that, we  compute $H(\fz)^*$.

To compute $H(\fz)^*$, we first use Lemma~\ref{intertwiners} and fix a sequence $\{A_i\}$ in $\cI(\rho,\sigma)$ such that $I=w^*-lim \sum_{i=1}^n A_iA_i^*$. We claim that, for $C\in B(H_{\rho},H_{\sigma})$ and $\fz\in \cU(\sigma)$,
$$ H(\fz)^*(C)=\sum_i A_i \otimes K_{(\fz,A_i^*,C)}.$$
To verify this, we fix $A\otimes K_{(\fw,B,D)} \in \cI(\rho,\sigma)\otimes_{\rho(M)'} E(K)$ and compute:
\begingroup
\allowdisplaybreaks
\begin{align*}
	\sum_i& \langle A_i \otimes K_{(\fz,A_i^*,C)},A\otimes K_{(\fw,B,D)} \rangle =\sum_i \langle  K_{(\fz,A_i^*,C)}, \varphi(A_i^*A)K_{\fw,B,D}\rangle\\& = \sum_i \langle  K_{(\fz,A_i^*,C)},K_{(\fw,A_i^*AB,D)}\rangle= \sum_i C^*K(\fz,\fw)(A_iA_i^*AB)D\\&=C^*K(\fz,\fw)(\sum_i A_iA_i^*AB)D=C^*K(\fz,\fw)(AB)D\\&=\langle C,K(\fz,\fw)(AB)D\rangle=\langle C, K_{(\fw,B,D)}(\fz)(A)\rangle=\langle C, H(\fz)(A \otimes  K_{(\fw,B,D)}) \rangle. 
\end{align*}
\endgroup
verifying the claim. (Note that we used the $w^*$-continuity of $K(\fz,\fw)$.)
It is left to prove the decomposition (\ref{kolmogorov}). For this, fix $P\in \cI(\tau,\sigma)$, $C\in B(H_{\rho},H_{\tau})$, $\fz\in \cU(\sigma)$ and $\fw\in \cU(\tau)$, and compute
\begingroup
\allowdisplaybreaks
\begin{align*}
	H(\fz)(P\otimes I)H(\fw)^*(C)&=\sum_i H(\fz)(P\otimes I)( A_i \otimes K_{(\fw,A_i^*,C)})=\sum_i H(\fz)(PA_i \otimes K_{(\fw,A_i^*,C)})\\&=\sum_i K_{(\fw,A_i^*,C)}(\fz)(PA_i)=\sum_i K(\fz,\fw)(PA_iA_i^*)C=K(\fz,\fw)(P)C.
\end{align*}
\endgroup
The last implication that $H$ is matricial if $K$ is follows immediately from (\ref{fmatricial}).

Finally, the fact that (3) implies (1) is straightforward from the definitions.
\end{proof}

\begin{remark}\label{subcategory}
We began this narrative talking about $\widehat{M}$ as a collection of representations of the $W^*-algebra$, $M$, but soon we turned to a discussion of the intertwiner spaces.  Had we been more meticulous we might have written $\widehat{M}^0$ in place of $\widehat{M}$ at the outset and declared it to be the space of \emph{objects} of Rieffel's category Normod-$M$. 
The morphism space for Normod-$M$, which consists of all the intertwiners, would then have been written ${\widehat{M}}^1$. 
This change in emphasis and notation makes it easy to talk about subcategories in a fashion suitable for our needs. 
So, by a subcategory $\Sigma = (\Sigma^0, \Sigma^1)$ of $\widehat{M}$, we mean that $\Sigma^0$ is a set of representations contained in $\widehat{M}^0$ and $\Sigma^1$ is a collection of intertwiner spaces between the objects in $\Sigma^0$.  The subcategory is called \emph{full} in case $\Sigma^1$ contains all the intertwiner spaces in $\widehat{M}^1$ that exist among the objects in $\Sigma^0$. The subcategory is called \emph{additive} in case $\Sigma^0$ is closed under direct sums.

Close inspection of the proof of Theorem~\ref{cpkerneltheorem} reveals that the proof is valid for any subcategory $\Sigma$ that is full, additive and contains a representation whose kernel is minimal among all the kernels of the representations in $\Sigma^0$. In particular, the result is true for any full, additive subcategory that contains a faithful representation. In this setting we write $E(K,\Sigma)$ to indicate the dependence of $E(K)$ on the subcategory $\Sigma$.
\end{remark}

The following example indicates how known results may fall in the umbra of our Theorem~\ref{cpkerneltheorem}. 

\begin{example}\label{DruryArveson}
Let $M=\mathbb{C}$, $E=\mathbb{C}^d$ and $\rho$ be the one dimensional representation of $\mathbb{C}$.
Clearly $\widehat{M}=\{n\rho :\; n\in \mathbb{N} \}$. For every $n\in \mathbb{N}$ we set $\cU(n\rho)$ to be the set of all contractive tuples $W=(W_1,\ldots,W_d)$ with $W_j\in M_n$ and $W_iW_j=W_jW_i$ for all $i,j$.  It is easy to see that $\cU=\{\cU(n\rho)\}$ is a matricial set. We now define $K:\cU \times \cU \rightarrow \coprod_{n,m} B_*(M_{n,m},M_{n,m})$ by
\be
K(Z,W)(A)=\Sigma_{\alpha \in \mathbb{F}^d_+} Z^{\alpha}A(W^{\alpha})^*
\ee
where $Z^{\alpha}=Z_{\alpha_1}Z_{\alpha_2}\cdots Z_{\alpha_k}$ if $|\alpha|=k$.

Then $K$ is a matricial kernel.

We claim that the correspondence $E(K)$ (which is a Hilbert space since $\rho(M)'=B(H_{\rho})=\mathbb{C}$) is the symmetric Fock space over $\mathbb{C}^d$ (that is, the Drury-Arveson space $H_d^2$).

The space $E(K)$ is generated by the functions $K_{(W,A,C)}$ where $W\in \cU(m)$ (for some $m\in \mathbb{N}$), $A\in M_{1,m}$ and $C\in M_{m,1}$. For every such function, we define $VK_{(W,A,C)}=\Sigma_{\alpha \in \mathbb{F}^d_+} A(W^{\alpha})^*C E_{\alpha}$ where $E_{\alpha}=E_{\alpha_1}\otimes \cdots \otimes E_{\alpha_{|\alpha|}}$ and $\{E_i\}$ is the standard basis in $\mathbb{C}^d$. Since $\{W_j\}$ commute, the coefficient of $E_{\alpha}$ is equal to the coefficient of $E_{\beta}$ if $\beta$ is a permutation of $\alpha$. It follows that $VK_{(W,A,C)}$ lies in the symmetric Fock space.

Also, we have,
$$\langle K_{(W,A,C)},K_{(W',A',C')} \rangle =C^*K(W,W')(A^*A')C'=\Sigma_{\alpha \in \mathbb{F}^d_+} C^*W^{\alpha}A^*A'(W^{'\alpha})^*C'$$   $$=\langle VK_{(W,A,C)},VK_{(W',A',C'} \rangle.$$
Thus, $V$ extends linearly to an isometry from $E(K)$ to the closed subspace generated by the vectors $\{VK_{(W,A,C)}=\Sigma_{\alpha \in \mathbb{F}^d_+} A(W^{\alpha})^*C E_{\alpha}\}$ and it is left to show that this space is the symmetric Fock space.

Since this space is invariant under the gauge group action (by considering the scalar multiple map, $W\mapsto \lambda W$), it suffices to show that for a fixed $k$, the set 
\begin{equation}\{\Sigma_{\alpha \in \mathbb{F}^d_+, |\alpha|=k} A(W^{\alpha})^*C 
E_{\alpha}\},\qquad A\in M_{1,m}, C\in M_{m,1}\end{equation} generates the symmetric tensor product $\mathbb{C}^d \circledS \cdots \circledS\mathbb{C}^d$ ($k$ 
copies).
It is enough, in fact, to let $W$ run over all scalar tuples ($W_j=H_jI
$ for some $H$ in the unit ball of $\mathbb{C}^m$), $A=(1,0,\cdots,0)$ and $C=A^t$ to get a generating set for the symmetric tensor product because what we get is $H\otimes \cdots \otimes H$ ($k$ times). This completes the proof of the claim.

\end{example}

The following concept is  minor modification of the notion of Morita equivalence of $C^*$-correspondences developed in \cite{MS00}.

\begin{definition}\label{ME}
Suppose $M,N, M_1$ and $N_1$ are $W^*$-algebras, $X$ is an $M-N$ $W^*$-correspondence and $Y$ is an $M_1-N_1$ $W^*$-correspondence. $X$ and $Y$ are called (strongly) Morita equivalent in case the algebras $M$ and $M_1$ are strongly Morita equivalent via an $M-M_1$ equivalence $W^*$-correspondence $Z$, the algebras $N$ and $N_1$ are strongly Morita equivalent via an $N-N_1$ equivalence $W^*$-correspondence $W$ and there is an isomorphism of correspondences
\be
 X\otimes_N W \cong Z \otimes_{M_1} Y.
\ee
If this holds, we write $X\sim^{SME} Y$ or $X\sim^{SME}_{W,Z} Y$.
\end{definition}

\begin{proposition}\label{propME}
Suppose $\Sigma\subseteq \widehat{M}$ is a full, additive subcategory and $\rho$ and $\pi$ are two  normal representations in $\Sigma$ with minimal kernels.
Suppose $K$ is a matricial cp kernel as in Theorem~\ref{cpkerneltheorem}, and $E^{\rho}(K,\Sigma)$, $E^{\pi}(K,\Sigma)$ are the $W^*$-correspondences constructed in the proof of Theorem~\ref{cpkerneltheorem}. Then these correspondences are strongly Morita equivalent.
\end{proposition}

Thus, with every such kernel one associates a Morita equivalence class of correspondences.

\begin{proof}
Note that $\cI(\rho,\pi)$ is a $\pi(M)'-\rho(M)'$ equivalence correspondence inducing the strong Morita equivalence of $\rho(M)'$ and $\pi(M)'$ because $\cI(\rho,\pi)\otimes_{\rho(M)'}\cI(\pi,\rho)\cong \pi(M)'$ and $\cI(\pi,\rho)\otimes_{\pi(M)'}\cI(\rho,\pi)\cong \rho(M)'$ (Lemma~\ref{intertwiners}). Similarly $B(H_{\rho},H_{\pi})$ is a $B(H_{\pi})-B(H_{\rho})$ equivalence correspondence inducing the strong Morita equivalence of $B(H_{\rho})$ and $B(H_{\pi})$ because $B(H_{\rho},H_{\pi})\otimes_{B(H_{\rho})}B(H_{\pi},H_{\rho})\cong B(H_{\pi})$ and $B(H_{\pi},H_{\rho})\otimes_{B(H_{\pi})}B(H_{\rho},H_{\pi})\cong B(H_{\rho})$.

We shall show that the correspondences $\cI(\rho,\pi)\otimes_{\rho(M)'}E^{\rho}(K,\Sigma)$ and $E^{\pi}(K,\Sigma)\otimes_{B(H_{\pi})}B(H_{\rho},H_{\pi}) $ are isomorphic.

It is enough to prove that the linear map that is defined on generators by
\be\label{map}
 D \otimes K^{\rho}_{(\fw,A,CB)} \mapsto K^{\pi}_{(\fw,DA,C)}\otimes B,
 \ee
for $\tau\in \Sigma$, $D\in \cI(\rho,\pi)$, $\fw\in \cU(\tau)$, $A\in \cI(\tau,\rho)$, $C\in B(H_{\pi},H_{\tau})$ and $B\in B(H_{\rho},H_{\pi}))$, is a well defined bimodule isometry.

Write $L=\text{span}\{ D \otimes K^{\rho}_{(\fw,A,CB)}: \tau\in \Sigma, D\in \cI(\rho,\pi), \fw\in \cU(\tau), A\in \cI(\tau,\rho), C\in B(H_{\pi},H_{\tau}) , B\in B(H_{\rho},H_{\pi})\}$ and $R=\text{span}\{K^{\pi}_{(\fw,DA,C)}\otimes B :\tau\in \Sigma, D\in \cI(\rho,\pi), \fw\in \cU(\tau), A\in \cI(\tau,\rho), C\in B(H_{\pi},H_{\tau}) , B\in B(H_{\rho},H_{\pi})\}$.

To show that the map in (\ref{map}) preserves the inner product we compute
\begin{multline}
\langle D \otimes K^{\rho}_{(\fw,A,CB)},D_1 \otimes K^{\rho}_{\fw_1,A_1,C_1B_1}\rangle=\langle K^{\rho}_{(\fw,A,CB)},\varphi(D^*D_1)K^{\rho}_{(\fw_1,A_1,C_1B_1)}\rangle\\   =\langle K^{\rho}_{(\fw,A,CB)},K^{\rho}_{(\fw_1,D^*D_1A_1,C_1B_1)}\rangle =B^*C^*K(\fw,\fw_1)(A^*D^*D_1A_1)C_1B_1 \\  =B^*\langle K^{\pi}_{(\fw,DA,C)},K^{\pi}_{(\fw_1,D_1A_1,C_1)}\rangle B_1\\=\langle K^{\pi}_{\fw,DA,C}\otimes B,K^{\pi}_{(\fw_1,D_1A_1,C_1)}\otimes B_1\rangle.
\end{multline}

Thus the map preserves the inner product. It is straightforward to check that the map is a bimodule map from $L$ onto $R$ . Since  $\cI(\rho,\pi)\otimes_{\rho(M)'}E(K,\rho)$ is the self-dual completion of $L$, $L$ is $w^*$-dense there and, similarly, $R$ is $w^*$-dense in $E(K,\pi)\otimes_{B(H_{\pi})}B(H_{\rho},H_{\pi}) $ and the map defined in (~\ref{map}) can be extended uniquely to an isomorphism of the $W^*$-correspondences.
\end{proof}

In the following proposition we clarify the relation of the construction of $E(K,\Sigma)$ in Theorem~\ref{cpkerneltheorem},  applied to a full additive subcategory $\Sigma$ to the construction of $\mathcal{H}(K)$ in \cite[Theorem 3.1]{BMV1}.

\begin{proposition}\label{BMV}
\begin{enumerate}
\item If in Theorem~\ref{cpkerneltheorem}, we set $M=\mathbb{C}$ and $\Sigma$ is the set of all finite dimensional representations of $\mathbb{C}$, identified with $\mathbb{N}$, and if $\rho$ is the one dimensional representation of $\mathbb{C}$, then we obtain the equivalence of (1) and (2) in \cite[Theorem 3.1]{BMV1} for the case where $\cA=\mathcal{Y}=\mathbb{C}$.
\item In Theorem~\ref{cpkerneltheorem}, fix normal representation $\rho$ of $M$, set $\Sigma$ to be all finite multiples of $\rho$, identified with $\mathbb{N}$, let $\cU=\{\cU(n\rho)\}$ be a matricial family of sets for $(M,E)$, and write $E(K,\rho)$ for the resulting $\rho(M)'- B(H_{\rho})$ correspondence. In the context of  \cite[Theorem 3.1]{BMV1}, let $\Omega=\cU$, $\mathcal{A}=\rho(M)'$ and $\mathcal{Y}=H_{\rho}$ . Then the two kernels may be identified and we have
    \be \label{HKEK}
    \mathcal{H}(K)\otimes (H_{\rho})^*\cong E(K,\rho)
    \ee
    where the tensor product is the exterior tensor product, $(H_{\rho})^*$ is viewed as a $W^*$-module over $B(H_{\rho})$,  and the isomorphism is a $\rho(M)'-B(H_{\rho})$ correspondence isomorphism, where $\mathcal{H}(K)$ is viewed as a $\rho(M)'-\mathbb{C}$ correspondence with the left action given by the map $\sigma$ of \cite[Theorem 3.1]{BMV1}.

\end{enumerate}

\end{proposition}
\begin{proof}
If, in (2), we set $M=\mathbb{C}$ and let $\rho$ to be the one dimensional representation of $\mathbb{C}$ we get (1). Thus, it is enough to prove part (2).

The kernel $K$ in \cite{BMV1} is a global map and, once we set $\Omega=\cU$, $\cA=\rho(M)'$ and $\mathcal{Y}=H_{\rho}$, we get a global map $K:\cU \times \cU \rightarrow \coprod_{n,m} B(M_{m,n}(\rho(M)'),M_{m,n}(B(H_{\rho})))$. But $M_{m,n}(\rho(M)')$ can be identified with $\cI(n\rho,m\rho)$ and $M_{m,n}(B(H_{\rho}))$ can be identified with $B(H_{n\rho},H_{m\rho})$. Thus, $K$ can be thought of as a global kernel in our setting (restricted to $\Sigma$ as above). These identifications preserve the property of being cp.

Now we turn to prove (\ref{HKEK}). For this, note first that we can identify $M_{1,m}(\cA)$ (where $\cA=\rho(M)'$ ) with $\cI(m\rho,\rho)$ and $\mathcal{Y}^m$ (where $\mathcal{Y}=H_{\rho}$ ) with $B(\mathbb{C}^m,H_{\rho})$. Also note that $H_{\rho}^m \otimes H_{\rho}^*$ is isomorphic (as a $C^*$-module over $B(H_{\rho})$) to $B(H_{\rho})^m$ via the obvious map). Using these identifications, it is enough to check that mapping $K_{(\fw,V,Y)}\otimes h^*\in \mathcal{H}(K)\otimes (H_{\rho})^*$ (with $h\in H_{\rho}$) to $K_{(\fw,V,Y\otimes h^*)}\in E(K,\rho) $ (with $y\otimes h^*$ viewed as an element of $B(H_{\rho})^m\cong B(H_{\rho},H_{m\rho})$) induces, by linearity, the required isomorphism.

First, we compute the inner product of generators.
\begin{multline}
    \langle K_{(\fw,V,Y)}\otimes h^*,K_{(\fw',V',Y')}\otimes k^*\rangle=\langle K_{\fw,V,Y},K_{\fw',V',Y'}\rangle \langle h^*,k^*\rangle\\=\langle Y, K(\fw,\fw')(V^{*}V')Y'\rangle_{H_{\rho}^m} (h\otimes k^*).
    \end{multline}
Note that the formula we use here for the inner product of two generators in $\mathcal{H}(K)$ is different from the one used in Equation (3.6) of \cite{BMV1} but we had to change it so that linearity will be in the second slot to conform with the inner product on $W^*$-modules.

On the other hand,
$$\langle K_{(\fw,V,Y\otimes h^*)},K_{(\fw',V',Y'\otimes k^*)}\rangle=(Y\otimes h^*)^*K(\fw,\fw')(V^{*}V')(Y'\otimes k^*)$$ and it is easy to see that the two expressions are equal.

Also, for $D\in \rho(M)'$ and $R\in B(H_{\rho})$, we have (where $\sigma$ is the map defined in \cite{BMV1}) $\sigma(D)K_{(\fw,V,Y)}\otimes h^*R=K_{(\fw,DV,Y)}\otimes (Rh)^*$ which is mapped to $K_{(\fw,DV,Y\otimes (R^*h)^*)}=K_{(\fw,DV,Y\otimes (Rh)^*)}=K_{(\fw,DV,(Y\otimes h^*)R)}=\varphi(D)K_{(\fw,V,Y\otimes h^*)}\cdot R$ where $\varphi$ is the left action on $E(K,\rho)$.

Thus, this induces a bimodule isometric and surjective map from $\mathcal{H}(K)\otimes (H_{\rho})^*$ onto $ E(K,\rho)$.
\end{proof}

\begin{definition}\label{multiplier}
A global [matricial] multiplier on $E(K,\Sigma)$ is a global [matricial] function $S: \cU \rightarrow \coprod_{\sigma \in \Sigma} B(H_{\sigma})$ such that, for every $f \in E(K,\Sigma)$, the function that maps $\fz \in \cU(\sigma)$ to $S(\fz)f(\fz)$ lies in $E(K,\Sigma)$. This function will be written $Sf$.

For a multiplier $S$, we write $M_S$ for the operator that takes $f$ to $Sf$.
\end{definition}

Note that $S$ being a matricial function means:
\begin{enumerate}
\item[(a)] For $\sigma\in \Sigma$ and $\fz \in \cU(\sigma)$, $S(\fz)\in B(H_{\sigma})$.
\item[(b)] If $\fz \in \cU(\sigma)$ and $\fw \in \cU(\tau)$ and $D$ intertwines $\fz$ and $\fw$ (that is, $D\in \cI(\sigma, \tau)$ and $D\fz = \fw (I_E\otimes D)$) then $S(\fw)D=DS(\fz)$.
\end{enumerate}

\begin{lemma}\label{adjointS}
For a multiplier $S$, the operator $M_S$ on $E(K)$ is adjointable and $$M_S^* K_{(\fw,A,C)}=K_{(\fw, A, S(\fw)^*C)}$$ 
where $S(\fw)^*$ is the adjoint of the operator $S(\fw)$.
\end{lemma}
\begin{proof}
Using (\ref{reproducing}) we compute
\begin{align*}
	\langle K_{(\fw,A,C)},SK_{(\fw',A',C')}\rangle = C^*S(\fw)K_{(\fw',A',C')}(\fw)(A^*)&=C^*S(\fw)K(\fw,\fw')(A^*A')C'\\&= \langle K_{(\fw,A,S(\fw)^*C)},K_{(\fw',A',C')} \rangle 
\end{align*}
\end{proof}

\begin{proposition}\label{SK}
For a global [matricial] function $S$ and a  global [matricial] kernel $K$, we can define a function $K_S$ by
\be
K_S(\fz,\fw)(B)=K(\fz,\fw)(B)-S(\fz)K(\fz,\fw)(B)S(\fw)^*
\ee
(where $\fz\in \cU(\sigma)$, $\fw\in \cU(\tau)$ and $B\in \cI(\tau,\sigma)$) then $K_S$ is a global [matricial] kernel  and the following are equivalent.
\begin{enumerate}
\item[(1)] $S$ is a multiplier on $E(K,\Sigma)$ with $\|M_S\|\leq 1$.
\item[(2)] The kernel $K_S$ is a cp kernel.
\end{enumerate}
\end{proposition}
\begin{proof}
A straightforward computation shows that
\be
\|K_{(\fw,A,C)}\|^2-\|K_{(\fw, A, S(\fw)^*C)}\|^2=C^*K_S(\fw,\fw)(A^*A) C.
\ee
The direction (1) implies (2) follows immediately and, for the converse, we can define a linear map $\Gamma$ on $E(K)^{(0)}$ by $\Gamma(K_{(\fw,A,C)})=K_{(\fw, A, S(\fw)^*C)}$ and linearity. Then we compute
\begingroup
\allowdisplaybreaks
\begin{align*}
	\big\|\sum_{j=1}^N K_{(\fw_j,A_j,C_j)}\big\|^2-\big\|\Gamma(\sum_{j=1}^N K_{(\fw_j,A_j,C_j)})\big\|^2&=\big\|\sum_{j=1}^N K_{(\fw_j,A_j,C_j)}\big\|^2-\big\|\sum_{j=1}^N K_{(\fw_j,A_j,S(\fw_j)^*C_j)})\big\|^2\\&=\sum_{i,j} C_i^*K_S(\fw_i,\fw_j)(A_i^*A_j) C_j. 
\end{align*}
\endgroup 
The last expression is $C^*K_S(\fw,\fw)(A^*A)C$ where $$ C=\left( \begin{array}{c} C_1 \\ \cdot\\ \cdot \\ C_N \end{array}\right) , \fw=\left(\begin{array}{cccc}\fw_1 & & & \\ & \cdot & & \\ & & \cdot & \\ & & & \fw_N \end{array}\right), A=\left( \begin{array}{cccc} A_1 & \cdot & \cdot & A_N \end{array}\right).$$
Since $K_S$ is assumed to be a cp kernel, we get
$$\big\|\Gamma(\sum_{j=1}^N K_{(\fw_j,A_j,C_j)})\big\| \leq \big\|\sum_{j=1}^N K_{(\fw_j,A_j,C_j)}\big\|.$$ 
Thus, $\Gamma$ is well defined and extends to a contraction on $E(K)$. A computation as in Lemma~\ref{adjointS} shows that $\Gamma=M_S^*$ and we conclude that $M_S$ is a well defined contraction and (1) follows.

\end{proof}
\begin{remark} Note that the definition of $E(K)$ depends on a choice of $\rho$ and, thus, it seems that the definition of a multiplier depends on $\rho$. But Proposition~\ref{SK} shows that it depends only on $K$ and not on $\rho$.
\end{remark}

\section{The relationship with nonmatricial kernels}\label{nonmatricial}
In \cite{G} and \cite{Meyer} the authors investigated reproducing kernel $W^*$-correspondences associated with non-matricial kernels.
The kernels they studied are maps $K:Z\times Z \rightarrow B_*(N,L)$ where $Z$ is a set and $N$, $L$ are $W^*$-algebras and $K$ is completely positive in the sense that given $n$ and points $z_1,\ldots,z_n$ in $Z$, the matrix (of maps) $(K(z_i,z_j))_{i,j=1}^{n}$ represents a (normal) completely positive map from $M_n(N)$ to $M_n(L)$. Such a maps-valued kernel gives rise to a positive operator-valued kernel $$K':(N\times Z)\times (N\times Z)\rightarrow L$$ defined by $$K'((a_1,z_1),(a_2,z_2))=K(z_1,z_2)(a_1a_2^*).$$
Associated with such a kernel one gets an $N-L$ $W^*$-correspondence $E_K$. They write the elements of $E_K$ as $L$-valued functions on $N\times Z$. It will be more convenient for us here to view these functions as functions $f:Z\rightarrow B(N,L)$.
The left action of $N$ on $f$ is defined by $(a\cdot f)(z)(b)=f(z)(ba)$.
Writing
\be
k_{(a,w)} (z)(b)=K(z,w)(ba^*)
\ee for $z,w \in Z$ and $a,b \in N$,we get, for $f\in E_K$, $$f(z)(a)=\langle k_{(a,z)},f\rangle $$ so these are the kernel functions that induce point evaluations. We also have
 $$ \langle k_{(b,z)},k_{(a,w)}\rangle = k_{(a,w)}(z)(b)= K(z,w)(ba^*) .$$
 The $W^*$-correspondence $E_K$ is generated by $\{k_{(a,w)}\cdot c: a\in N, w\in Z, c\in M\}$ and
 $$\langle k_{(b,z)}\cdot c,k_{(a,w)}\cdot d\rangle = c^*K(z,w)(ba^*)d \in L .$$

We shall now let $K$, $M$, $\cU$ and $E$ be as in Theorem~\ref{cpkerneltheorem}. Given a faithful representation $\rho$ of $M$, we get, by part (2) of that theorem, a $\rho(M)'-B(H_{\rho})$ $W^*$-correspondence $E^{\rho}(K)$.
Letting $$K^{\rho}:\cU(\rho)\times \cU(\rho) \rightarrow B_*(\rho(M)',B(H_{\rho}))$$ be the restriction of $K$ to $\cU(\rho)\times \cU(\rho)$, we get a (non matricial) completely positive kernel in the sense studied by Good and Meyer where $Z=\cU(\rho)$, $N=\rho(M)'$ and $L=B(H_{\rho})$. We, therefore, get a $\rho(M)'-B(H_{\rho})$ $W^*$-correspondence $E_{K^{\rho}}$. The following lemma clarifies the relationship between these two correspondences.

\begin{lemma}\label{lambda}
With the notation above, there is an embedding of correspondences
$$\Lambda_{\rho}:E_{K^{\rho}} \rightarrow E^{\rho}(K)$$
mapping $k_{(a,\fw)}\cdot c$ to $K_{(\fw,a^*,c)}$.
\end{lemma}
\begin{proof}
This follows from a straightforward computation. For $\fw,\fz \in \cU(\rho)$, $a,b \in \rho(M)'$ and $c,d \in B(H_{\rho})$, we have
$$\langle K_{\fw,a^*,c},K_{\fz,b^*,d}\rangle = c^*K(\fw,\fz)(ab^*)d=c^*K^{\rho}(\fw,\fz)(ab^*)d=\langle k_{(a,\fw)}\cdot c,k_{(b,\fz)}\cdot d \rangle .$$
\end{proof}
Recall (Proposition~\ref{propME}) that, given two faithful representations $\rho$ and $\pi$, the correspondences $E^{\rho}(K)$ and $E^{\pi}(K)$ are Morita equivalent. In fact, $\cI(\rho,\pi)\otimes_{\rho(M)'}E^{\rho}(K)\otimes_{B(H_{\rho})}B(H_{\rho},  H_{\pi})$ is isomorphic to $E^{\pi}(K)$. When $\rho \preceq \pi$ (that is, $\rho$ is equivalent to a subrepresentation of $\pi$), we have the following.

\begin{proposition}
If $\rho$ and $\pi$ are faithful representations of $M$ and $\rho \preceq \pi$ and $\Psi:\cI(\rho,\pi)\otimes_{\rho(M)'}E^{\rho}(K)\otimes_{B(H_{\rho})}B(H_{\rho},  H_{\pi})\rightarrow E^{\pi}(K)$ is the isomorphism mentioned above, then
\be\label{lambdarhosigma}
\Psi(\cI(\rho,\pi)\otimes_{\rho(M)'}\Lambda_{\rho}(E_{K^{\rho}})\otimes_{B(H_{\rho})}B(H_{\rho},  H_{\pi}))\subseteq \Lambda_{\pi}(E_{K^{\pi}}).
\ee

\end{proposition}
\begin{proof}
For $d\in \cI(\rho,\pi)$, $g\in B(H_{\pi},H_{\rho})$ and $K^{\rho}_{\fw,a,e}\in E^{\rho}(K)$ we have $$\Psi(d\otimes K_{\fw,a,e}\otimes g)=\varphi(d)K^{\rho}_{\fw,a,e}\cdot g=K^{\pi}_{\fw, da,eg}.$$ We now assume that $K^{\rho}_{\fw,a,e}\in \Lambda_{\rho}(E_{K^{\rho}})$. So, $\fw\in \cU(\rho)$, $a\in \rho(M)'$ and $e\in B(H_{\rho})$ and we should show that there are $\tilde{\fw}\in \cU(\pi)$, $t\in \pi(M)'$ and $c\in B(H_{\pi})$ such that
\be\label{Kpi}
K^{\pi}_{\fw, da,eg}= K^{\pi}_{\tilde{\fw},t,c}.
\ee
Since $\rho \preceq \pi$, there is some $D\in \cI(\rho,\pi)$ such that $D^*D=I_{H_{\rho}}$. We claim that $\tilde{\fw}=D\fw (I_{E}\otimes D^*)$, $t=daD^*$ and $c=Deg$ will do. It is easy to see that $\tilde{\fw}\in \cU(\pi)$, $t\in \pi(M)'$ and $c\in B(H_{\pi})$. Also note that $D\fw=\tilde{\fw}(I_E\otimes D)$. It is left to verify (\ref{Kpi}). For that, we fix a representation $\sigma$, $\fz\in \cU(\sigma)$, $b\in \cI(\pi,\sigma)$ and compute, using(\ref{kernelintertwining}),
\begin{align*}
	K^{\pi}_{\fw, da,eg}(\fz)(b)&=K(\fz,\fw)(bda)eg=K(\fz,\fw)(bda)D^*Deg\\&=
K(\fz,\tilde{\fw})(bdaD^*)Deg=
K(\fz,\tilde{\fw})(bt)c=K^{\pi}_{\tilde{\fw},t,c}(\fz)(b).
\end{align*}
This completes the proof.
\end{proof}
We do not know if, in general, we have equality in (\ref{lambdarhosigma})

The following result relates the multipliers of $E(K,\Sigma)$ with the multipliers of the correspondences $E(K^{\rho})$ ($\rho \in \Sigma$).

\begin{proposition}\label{Multipliers}
Given a multiplier $S:\cU \rightarrow \coprod_{\sigma \in \Sigma} B(H_{\sigma})$, the family $\{\phi_{\sigma}\}$, where $\phi_{\sigma}:\cU(\sigma)\rightarrow B(H_{\sigma})$ is the restriction of $S$ to $\cU(\sigma)$ , is a matricial family of functions and $\phi_{\sigma}$ is a multiplier of $E(K^{\sigma})$.
Conversely, every such family gives rise to a multiplier of $E(K,\Sigma)$.

We have $||S||=\sup\{||\phi_{\sigma}||\} $.

\end{proposition}
\begin{proof}
Given a multiplier $S$ of $E(K,\Sigma)$ with $||S||\leq 1$, write $\phi_{\sigma}:\cU(\sigma)\rightarrow B(H_{\sigma})$ for the restriction of $S$ to $\cU(\sigma)$. The fact that this is a matricial family follows from the fact that $S$ is matricial.
Note that, for $\fz\in \cU(\sigma)$ and $b\in \sigma(M)'$,
\be\label{KS}
K_S(\fz,\fz)(b)=K(\fz,\fz)(b)-S(\fz)K(\fz,\fz)(b)S(\fz)^*\ee
\[\qquad\qquad=K(\fz,\fz)(b)-\phi_{\sigma}(\fz)K(\fz,\fz)(b)\phi_{\sigma}(\fz)^* = K^{\sigma}_{\phi_{\sigma}}(\fz,\fz)(b).
\]
Since $S$ is a contractive multiplier, $K_S(\fz,\fz)$ is positive and, therefore, $K^{\sigma}_{\phi_{\sigma}}(\fz,\fz)$ is positive . Thus $\phi_{\sigma}$ is a contractive multiplier for $E(K^{\sigma})$ for every $\sigma$.
For the converse, we start with a matricial family $\{\phi_{\sigma}\}$ of contractive multipliers of $E(K^{\sigma})$. We let $S:\cU\rightarrow \coprod_{\sigma}B(H_{\sigma})$ be defined such that $S$, restricted to $\cU(\sigma)$, is $\phi_{\sigma}$ and then Equation (\ref{KS}) shows that $K_S$ is positive, thus $S$ is a contractive multiplier of $E(K,\Sigma)$.
\end{proof}

\section{Example: Weighted shift algebra}\label{example}
The kernel we study here can be viewed as a non commutative analogue of the unitarily invariant kernels discussed in \cite{CH}. To set up the notation, we start with a $W^*$-correspondence $E$ over $M$ and fix a sequence $\{R_k\}_{k=0}^{\infty}$ satisfying
\begin{enumerate}
\item[(1)] For each $k$, $R_k$ is a positive, invertible operator in $\varphi_k(M)^c$ ($:=\cL(E^{\otimes k})\cap \varphi_k(M)'$) and $R_0=I$.
\item[(2)] $\sup_k ||R_k^{-1}(I_E\otimes R_{k-1})||<\infty$.
\item[(3)] $\limsup_k ||R_k^2||^{\frac{1}{k}}<\infty$
\item[(4)] The unique solution $\{X_k\}_{k=1}^{\infty}$ (with $X_k\in \varphi_k(M)^c$) of the set of equations
\be\label{RX}
\Sigma_{l=1}^k X_l \otimes R_{k-l}^2=R_k^2 \;, \;\; k\geq 1
\ee
satisfies $X_k\geq 0$ (for all $k$).
\end{enumerate}

\begin{lemma}
Let $\{R_k\}$ be as above. Then equations (\ref{RX}) have a unique solution $\{X_k\}$ (with $X_k\in \varphi_k(M)^c$) and this solution satisfies:
$X_1$ is invertible and $\limsup_k ||X_k||^{\frac{1}{k}} <\infty$. Thus, this sequence is admissible in the sense of \cite[Definition 4.1]{MS16}.

\end{lemma}
\begin{proof}
The uniqueness is straightforward. In fact $X_1=R_1^2$ (using the equation with $k=1$) and, thus, $X_1$ is invertible. When $k=2$ we get $X_2=R_2^2-X_1\otimes R_1^2=R_2^2-R_1^2\otimes R_1^2$ and, continuing inductively, we can get, for $k\geq 1$,
\be
X_k=\Sigma_{l=1}^k (-1)^{l+1}\Sigma_{\alpha\in F(k,l)}\otimes_{i=1}^l R^2_{\alpha(i)}
\ee
where $F(k,l):=\{\alpha:\{1,2,\ldots,l\}\rightarrow \mathbb{N} : \; \Sigma_{i=1}^l \alpha(i)=k \}$.
We shall not need this explicit formula but the uniqueness should be clear.

The positivity requirement in \cite[Definition 4.1]{MS16} does not follow from this formula and has to be assumed (as we did above in the requirements of $\{R_k\}$). It is left to prove that $\limsup_k ||X_k||^{\frac{1}{k}} <\infty$. For that, note that we assume that $\limsup_k ||R_k^2||^{\frac{1}{k}} <\infty$ and this implies that there is some $C<\infty$ with $||R_k^2||\leq C^k $ for all $k$.
To prove $\limsup_k ||X_k||^{\frac{1}{k}} <\infty$ it suffices to show that, for every $k$, $||X_k||\leq (2C)^k$.
This will be shown by induction. For $k=1$, $||X_1||=||R_1^2||\leq C\leq 2C$. Now assume that $||X_l||\leq (2C)^l$ for all $l<k$.

Using (\ref{RX}), we have, for every $k$,
$ X_k=R_k^2-\Sigma_{l=1}^{k-1} X_l\otimes R^2_{k-l}$ and, thus,
\begin{align*}
	\|X_k\|\leq \|R_k^2\|+ \Sigma_{l=1}^{k-1} \|X_l\| C^{k-l}&\leq C^k+ \Sigma_{l=1}^{k-1} (2C)^l C^{k-l}\\&=(\Sigma_{l=0}^{k-1} 2^l)C^k=(2^k-1)C^k<(2C)^k.
\end{align*}
This completes the proof.
\end{proof}

\begin{definition}\label{def:DX}
\begin{enumerate}
\item[(1)] With $X=\{X_k\}$ as above (an admissible sequence) and $\sigma\in \widehat{M}$ we associate the following set
$$\mathbb{D}(X,\sigma)=\{ \fz\in \cI(\sigma^E\circ \varphi, \sigma) : \; ||\Sigma_{k=1}^{\infty} \fz^{(k)} (X_k\otimes I_{H_{\sigma}})\fz^{(k) *}||<1 \} .$$ (See also \cite[Definition 4.3]{MS16}). We also write
\be
\mathbb{D}_X=\{\mathbb{D}(X,\sigma): \sigma \in \widehat{M}\}.
\ee
\item[(2)] Define the sequence $Z=\{Z_k\}_{k\geq 0}$ by $$Z_k=R_k^{-1}(I_E\otimes R_{k-1})\;\;\;, k\geq 1$$ (and $Z_0=I$ ) and write $$Z^{(k)}:=Z_k(I_E\otimes Z_{k-1})\cdots (I_{E^{\otimes (k-1)}}\otimes Z_1)=R_k^{-1}.$$
\end{enumerate}
\end{definition}

\begin{proposition}\label{weighted Szego}
\begin{enumerate}
\item[(i)] The sequence $Z$ is a weight sequence associated with the admissible sequence $X$ in the sense of \cite[Definition 4.6]{MS16}. In fact, it is the canonical weight sequence associated with $X$.
\item[(ii)] Given $\sigma, \tau \in \widehat{M}$, for $\fz\in \mathbb{D}(X,\sigma)$ and $\fw \in \mathbb{D}(X,\tau)$, the map
\be \label{Kweighted}
K(\fz,\fw)(a)=\sum_{k=0}^{\infty} \fz^{(k)}(R_k^2 \otimes a)\fw^{(k) *}
 \ee
is well defined (in the sense of norm convergence) and maps $\cI(\tau,\sigma)$ into itself.
\item[(iii)] $\mathbb{D}_X=\{\mathbb{D}(X,\sigma)\}_{\sigma}$ is a matricial family of sets for $(M,E)$ and
$$ K:\mathbb{D}_X\times \mathbb{D}_X \rightarrow \coprod_{\sigma,\tau \in \widehat{M}} B_{*}(\cI(\tau,\sigma),\cI(\tau,\sigma))\subseteq \coprod_{\sigma,\tau \in \widehat{M}} B_{*}(\cI(\tau,\sigma),B(H_{\tau},H_{\sigma}))$$ is a cp matricial kernel.
\end{enumerate}

\end{proposition}
\begin{proof}
(i) To show that $Z$ is a weight sequence associated with $X$ we need to show that, for all $k$, $Z_k\in \varphi_k(M)^c$, $Z_k$ is invertible and $Z^{(k) *}Z^{(k)}=R_k^{-1}$. Also, $Z_0=I$ and $\sup||Z_k||<\infty$. The canonical weight sequence associated with $X$ is the unique one that satisfies also $Z^{(k)}\geq 0$ for all $k$. All these properties are obvious from the above.

(ii) For $\sigma_0\in \widehat{M}$ and $\fz_0\in \mathbb{D}(X,\sigma_0)$, it was proved in \cite[Theorem 4.5]{MS16} that $\sum_{k=0}^{\infty} \fz_0^{(k)}(R_k^2 \otimes a)\fz_0^{(k) *}$ converges in norm and maps $\sigma_0(M)'$ to itself. Setting $\sigma_0=\sigma \oplus \tau$, $\fz_0=\fz\oplus \fw$ and considering the off diagonal elements (in the resulting $2\times 2$ matrix) proves (ii).

(iii) It is straightforward to check that $\mathbb{D}_X=\{\mathbb{D}(X,\sigma)\}_{\sigma}$ is a matricial family of sets for $(M,E)$. The rest is also easy to check noting that, if $C\fz=\tilde{\fz}(I_E\otimes C)$ (as in Definition~\ref{kernel} (3)), then, for $k\geq 1$, $C\fz^{(k)}=\tilde{\fz}^{(k)}(I_{E^{\otimes k}}\otimes C)$.
\end{proof}

\begin{remark}
The condition $\limsup ||R_k^2||^{\frac{1}{k}}<\infty$ implies that the kernel $K$ could be defined on balls of some positive radius (instead of $\mathbb{D}(X,\sigma)$). For the case where $M=\mathbb{C}$ (and $\sigma $ is fixed, not an nc kernel) this was done in \cite{CH} (where the conditions on the coefficients are slightly different to make the domain the unit ball). But, as the results of \cite{MS16} and \cite{G} show, the natural domain associated with $K$ is $\mathbb{D}(X,\sigma)$. In fact, the results of \cite{G} also show that (for a fixed $\sigma$) the resulting kernel is a complete Nevanlinna-Pick kernel.
\end{remark}

\begin{proposition}\label{EK}
Assume $\Sigma \subseteq \widehat{M}$ is the subset of all faithful representations, $\rho\in \Sigma$ is fixed, $K$ is given in (\ref{Kweighted}) and $E(K,\Sigma)$ is the associated correspondence. Recall that $E^{\rho}$ is the $W^*$-corresponcence $\cI(\rho,\rho^E \circ \varphi)$.
Then there is an isomorphism of $W^*$-correspondences
$$V: E(K,\Sigma) \rightarrow \mathcal{F}(E^{\rho})\otimes_{\rho(M)'}B(H_{\rho}) $$
such that 
$$K_{\fw,a,cr}\mapsto\big(\Sigma_{k=0}^{\infty}(R_k \otimes a)(\fw^*)^{(k)}c\big)\otimes r$$
for all $\fw \in \mathbb{D}(X,\tau)\subseteq \cI(\tau^E\circ \varphi , \tau)$, $a\in \cI(\tau,\rho)$, $c\in \cI(\rho,\tau)$ and $r\in B(H_{\rho})$. 

\end{proposition}
\begin{proof}
Let $L=\text{span}\{K_{\fw,a,cr}: \fw \in \mathbb{D}(X,\tau), a\in \cI(\tau,\rho), c\in \cI(\rho,\tau), r\in B(H_{\rho})\}$ and let $V_0$ be the linear map from $L$ into $\mathcal{F}(E^{\rho})\otimes_{\rho(M)'}B(H_{\rho})$ that maps $K_{\fw,a,cr}$ to $(\Sigma_{k=0}^{\infty}(R_k \otimes a)(\fw^*)^{(k)}c)\otimes r$. It is straightforward to check that $V_0$ is a bimodule map.
We first compute
\begingroup
\allowdisplaybreaks
\begin{align*}
	\langle K_{\fw,a,cr},K_{\fw',a',c'r'}\rangle = r^*c^*K(\fw,\fw')(a^*a')c'r'&=\sum_{k=0}^{\infty} r^*c^*\fw^{(k)}(R_k^2 \otimes a^*a')\fw'^{(k) *}c'r'\\=\sum_{k=0}^{\infty} &((R_k \otimes a)\fw^{(k) *}cr)^*((R_k \otimes a')\fw'^{(k) *}c'r').
	\end{align*}
	\endgroup
	 For a fixed $k$, $(R_k \otimes a)\fw^{(k) *}c \in (E^{\rho})^{\otimes k}=\cI(\rho, \rho^{E^{\otimes k}}\circ \varphi )$ since $c\in \cI(\rho,\tau)$, $\fw^{(k) *}\in \cI(\tau^{E^{\otimes k}}\circ \varphi,\tau)$ and $R_k\in \varphi(M)'$. Thus $\sum_{k=0}^{\infty}(R_k \otimes a)\fw^{(k) *}c \in \mathcal{F}(E^{\rho})$. We, thus, have
$$\langle K_{\fw,a,cr},K_{\fw',a',c'r'}\rangle =\langle (\sum_{k=0}^{\infty}(R_k \otimes a)\fw^{(k) *}c)\otimes r,(\sum_{k=0}^{\infty}(R_k \otimes a')\fw'^{(k) *}c')\otimes r'\rangle.$$
So the map $V_0$, defined above on $L$ , preserves inner products.

It is left to show that $L$ is $w^*$-dense in $E(K,\Sigma)$ and the image of $V_0$ is $w^*$-dense in $\mathcal{F}(E^{\rho})\otimes_{\rho(M)'}B(H_{\rho})$.

The fact that $L$ is $w^*$-dense in $E(K,\Sigma)$ follows from the definition of $E(K,\Sigma)$ and the fact that, for $\fw\in \mathbb{D}(X,\tau)$, $a\in\cI(\tau,\rho)$ and $b\in B(H_{\rho},H_{\tau})$, $K_{\fw,a,b}$ lies in the $w^*$-closure of $L$ since $B(H_{\rho},H_{\tau})=\overline{\cI(\rho,\tau)B(H_{\rho})}^{w^*}$ (Lemma~\ref{corintertwiners}(1)).

To show that the image of $V_0$ is $w^*$-dense in $\mathcal{F}(E^{\rho})\otimes_{\rho(M)'}B(H_{\rho})$ it suffices to prove that $C:=\text{span}\{\Sigma_{k=0}^{\infty}(R_k \otimes a)(\fw^*)^{(k)}c:\fw \in \mathbb{D}(X,\tau), a\in \cI(\tau,\rho), c\in \cI(\rho,\tau)\}$ is $w^*$-dense in $\mathcal{F}(E^{\rho})$. But $\overline{C}^{w^*}$ is a sub-bimodule of $\mathcal{F}(E^{\rho})$ and is invariant under the action of the gauge group (just replace $\fw$ by $\lambda \fw$ for $|\lambda|=1$). Therefore it suffices to show that, for every $k\geq 0$, the set $C_k:=\{(R_k \otimes a)(\fw^*)^{(k)}c : \tau\in \Sigma, \fw \in \mathbb{D}(X,\tau), a\in \cI(\tau,\rho), c\in \cI(\rho,\tau)\}$ is $w^*$-dense in $(E^{\rho})^{\otimes k }$ which can be identified with $\cI(\rho,\rho^{E^{\otimes k}}\circ \varphi)$.

But, since $R_k\otimes I_{H_{\rho}}$ is invertible, it is enough to show that $C^0_k:=\{(I_{E^{\otimes k}} \otimes a)(\fw^*)^{(k)}c : \tau\in \Sigma, \fw \in \mathbb{D}(X,\tau), a\in \cI(\tau,\rho), c\in \cI(\rho,\tau)\}$ is $w^*$-dense in $(R_k^{-1}\otimes I)(E^{\rho})^{\otimes k }$. Since $R_k\otimes I_{H_{\rho}}$ is invertible and commutes with $\varphi_{E^{\otimes k}}(M)\otimes I_{H_{\rho}}=(\rho^{E^{\otimes k}}\circ \varphi)(M)$, and since we identify $(E^{\rho})^{\otimes k }$ with $\cI(\rho,\rho^{E^{\otimes k}}\circ \varphi)$, we have $(R_k^{-1}\otimes I)(E^{\rho})^{\otimes k }=(E^{\rho})^{\otimes k }$. Therefore, we should show that $C^0_k$ is $w^*$-dense in $(E^{\rho})^{\otimes k }$.

But $(E^{\rho})^{\otimes k }$ is generated by elements of the form $\eta_k \otimes \eta_{k-1}\otimes \cdots \otimes \eta_1$ for $\eta_i\in E^{\rho}$ (identified with $\cI(\rho,\rho^{E}\circ \varphi)$). Every such element is being identified with the element $(I_{E^{\otimes k-1}}\otimes \eta_k)\cdots (I_{E^{\otimes 2}}\otimes \eta_3)(I_{E}\otimes \eta_2)\eta_1 \in \cI(\rho,\rho^{E^{\otimes k}}\circ \varphi)$. So, it suffices to show that $(I_{E^{\otimes k-1}}\otimes \eta_k)\cdots (I_{E^{\otimes 2}}\otimes \eta_3)(I_{E}\otimes \eta_2)\eta_1 \in C^0_k$.

For this, let $\tau$ be $k\rho$ and note that $\cI(\tau,\tau^E\circ \varphi)$ is isomorphic to $M_k(\cI(\rho,\rho^E\circ \varphi))$ so that we can define $\fw^* \in \cI(\tau,\tau^E\circ \varphi)$ as a $k\times k$ matrix over $\cI(\rho,\rho^E\circ \varphi)$ by
\be
\fw^*=\left( \begin{array}{cccccc} 0 & \eta_k & 0 & 0 & 0 & 0 \\ 0 & 0 & \eta_{k-1} & 0 & 0 & 0 \\0&0&0&\cdot&0&0 \\0 & 0 & 0 & 0 & \cdot & 0  \\ 0 & 0 & 0 & 0 & 0 & \cdot \\ \eta_1 & 0&\cdot &\cdot &\cdot & 0 \end{array}\right) .
\ee
Now $(\fw^*)^{(k)}=(I_{E^{\otimes k-1}}\otimes \fw^*)\cdots (I_{E^{\otimes 2}}\otimes \fw^*)(I_{E}\otimes \fw^*)\fw^*$ which is a $k\times k$ matrix with $(I_{E^{\otimes k-1}}\otimes \eta_k)\cdots (I_{E^{\otimes 2}}\otimes \eta_3)(I_{E}\otimes \eta_2)\eta_1$ in the $(1,1)$ entry. Note that elements of  $\cI(\tau,\rho)=\cI(k\rho,\rho)$ can be written as rows of length $k$ over $\rho(M)'$ and elements of $\cI(\rho,\tau)=\cI(\rho,k\rho)$ can be written as columns of length $k$ over $\rho(M)'$. Thus we can set $a=(I,0,\ldots,0)\in \cI(\tau,\rho)$ and $c=a^t\in \cI(\rho,\tau)$ and we get $(I_{E^{\otimes k-1}}\otimes \eta_k)\cdots (I_{E^{\otimes 2}}\otimes \eta_3)(I_{E}\otimes \eta_2)\eta_1 =(I_{E^{\otimes k}} \otimes a)(\fw^*)^{(k)}c$. Take $t:=\limsup{\|X_k\|^{1/k}}$ and $r>\|\fw\|$. Consider 
 \[
\fw_0':=t^{-1/2}\frac{\fw}{r}.
 \]
So, $\|\fw_0'\|<t^{-1/2}$ and then, by \cite[Lemma 4.2]{MS16}, 
\[
\sum_{k=1}^\infty \|\fw_0'^{(k)}(X_k\otimes I_{\Hs})\fw_0'^{(k)*}\|<M<\infty
\]
for some $M>0$. If $M= 1$ we take $\fw_0=\fw_0'$. If $M>1$, we take $\fw_0=\frac{\fw_0'}{\sqrt{M}}$. For both cases,
\[
\fw_0^{(k)}=\frac{\fw_0'^{(k)}}{(\sqrt{M})^{k}}\quad\text{and}\quad \fw_0^{(k)*}=\frac{\fw_0'^{(k)*}}{(\sqrt{M})^{k}}.
\]
Therefore, 
\[
\sum_{k=1}^\infty \|\fw_0^{(k)}(X_k\otimes I_{\Hs})\fw_0^{(k)*}\|\leq\frac{1}{M}\sum_{k=1}^\infty \|\fw_0'^{(k)}(X_k\otimes I_{\Hs})\fw_0'^{(k)*}\|<1.
\]
Thus, $\fw_0=\frac{\fw_0'}{\sqrt{M}}=\frac{\fw}{r\sqrt{M t}}\in \bD(X,\tau)$. So, $(I_{E^{\otimes k}} \otimes a)(\fw_0^*)^{(k)}c\in C^0_k$. Therefore, $(I_{E^{\otimes k}} \otimes a)(\fw^*)^{(k)}c=(r\sqrt{M t})(I_{E^{\otimes k}} \otimes a)(\fw_0^*)^{(k)}c\in C^0_k$. Hence, $(I_{E^{\otimes k-1}}\otimes \eta_k)\cdots (I_{E^{\otimes 2}}\otimes \eta_3)(I_{E}\otimes \eta_2)\eta_1 \in C^0_k$.
\end{proof}

\begin{remark}\label{directsum}
Since $\mathcal{F}(E^{\rho})\otimes_{\rho(M)'}B(H_{\rho})$ can be written as a direct sum of correspondences
$$\mathcal{F}(E^{\rho})\otimes_{\rho(M)'}B(H_{\rho})=\Sigma^{\oplus}_{k\geq 0} (E^{\rho})^{\otimes k}\otimes_{\rho(M)'}B(H_{\rho}),$$
we can use $V$ (of Proposition~\ref{EK}) to write $E(K,\Sigma)$ as a direct sum
$$E(K,\Sigma)=\Sigma^{\oplus} V^*((E^{\rho})^{\otimes k}\otimes_{\rho(M)'}B(H_{\rho})).$$
In fact, if we write $Q_k$ for the projection onto the $k$th summand, we get
$$Q_k(K_{\fw,a,c})(\fz)(b)=\fz^{(k)}(R_k^2 \otimes ba)(\fw^*)^{(k)}c .$$
In the following we shall write $\tilde{V}$ for the unitary operator
$$\tilde{V}:=V\otimes I_{H_{\rho}}: E(K,\Sigma)\otimes_{B(H_{\rho})} H_{\rho} \rightarrow \mathcal{F}(E^{\rho})\otimes_{\rho(M)'}B(H_{\rho})\otimes_{B(H_{\rho})} H_{\rho}\cong \mathcal{F}(E^{\rho})\otimes_{\rho(M)'} H_{\rho} $$ and $\tilde{V}_k$ for its restriction to $(Q_k\otimes I)(E(K,\Sigma)\otimes_{B(H_{\rho})} H_{\rho}) $.

\end{remark}

\begin{corollary}
With $K,\Sigma$ and $\rho$ as in the proposition, $E(K,\Sigma)$ is isomorphic, as a $\rho(M)'-B(H_{\rho})$ $W^*$-correspondence, to $B(H_{\rho},\mathcal{F}(E)\otimes_{\rho}H_{\rho})$ where the structure of $B(H_{\rho},\mathcal{F}(E)\otimes_{\rho}H_{\rho})$ as a $W^*$-module over $B(H_{\rho})$ is the obvious one and the left action of $b\in \rho(M)'$ is given by $b\cdot T=(I_{\mathcal{F}(E)}\otimes b)\circ T$.

\end{corollary}
\begin{proof}
We will prove the isomorphism of $B(H_{\rho},\mathcal{F}(E)\otimes_{\rho}H_{\rho})$ and $\mathcal{F}(E^{\rho})\otimes_{\rho(M)'}B(H_{\rho}) $. Note that it follows from Lemma~\ref{corintertwiners}(1) that
$B(H_{\rho},\mathcal{F}(E)\otimes_{\rho}H_{\rho})$ is isomorphic to $\cI(\rho,\rho^{\cF(E)}\circ \varphi_E)\otimes_{\rho(M)'}B(H_{\rho})$. Since $\cI(\rho,\rho^{\cF(E)}\circ \varphi_E)=(\cF(E))^{\rho}\cong \cF(E^{\rho})$ (\cite[Lemma 3.8 (1)]{MS04}), we are done.
\end{proof}

\begin{lemma}\label{UVMSVU}
Let $U$ be the unitary map from $\cF(E^{\rho})\otimes_{\rho(M)'}H_{\rho}$ onto $\cF(E)\otimes_{\rho}H_{\rho}$ as in \cite[Theorem 3.9]{MS04} (see also \cite[Section 7]{MS16}). Let $C=\{C_m\}$, $H^{\infty}(E^{\rho},C)\subseteq \cL(\cF(E^{\rho}))$ and $W'_{\eta}$ (for $\eta\in E^{\rho}$) be as defined in \cite[Discussion following Theorem 7.6]{MS16}.
Suppose $S$ is in $M(K)$, $V$ is the unitary map of Proposition~\ref{EK} and $\tilde{V}$ is as in Remark~\ref{directsum}. Then $\tilde{V}(M_S\otimes I_{H_{\rho}})\tilde{V}^*$ commutes with $\iota^{\cF(E^{\rho})}(H^{\infty}(E^{\rho},C))$. Thus
$$U\tilde{V}(M_S\otimes I_{H_{\rho}})\tilde{V}^*U^* \in \rho^{\cF(E)}(H^{\infty}(E,Z)) .$$

\end{lemma}
\begin{proof}
Once we show that $\tilde{V}(M_S\otimes I_{H_{\rho}})\tilde{V}^*$ commutes with $\iota^{\cF(E^{\rho})}(H^{\infty}(E^{\rho},C))$, the last statement would follow from \cite[Theorem 7.7]{MS16}. Since $H^{\infty}(E,C)$ is generated by $\{W'_{\eta}: \eta \in E^{\rho}\}\cup \varphi_{\infty}(\rho(M)')$, it suffices to show that $\tilde{V}(M_S\otimes I_{H_{\rho}})\tilde{V}^*$ commutes with these generators. In fact, it is enough to prove, for every  $\eta \in E^{\rho}$ and $d\in \rho(M)'$,
\be\label{commuteeta}
(W'_{\eta})^*VM_S^*V^*=VM_S^*V^*(W'_{\eta})^*
\ee  and
\be\label{commuted}
\varphi_{\infty}(d)VM_S^*V^*=VM_S^*V^*\varphi_{\infty}(d).
\ee
By the definitions of $V$ and $M_S$, we have
$$(M_S^*\otimes I_{H_{\rho}})(V\otimes I_{H_{\rho}})^*((\Sigma_{k=0}^{\infty}(R_k \otimes a)(\fw^*)^{(k)}c)\otimes rh)=K_{\fw,a,S(\fw)^*cr}\otimes h.$$ To apply $V\otimes I_{H_{\rho}}$ to it we need to use Lemma~\ref{corintertwiners} and write $S(\fw)^*cr$ as the sum (in the $w^*$-topology) $\Sigma_i c_ir_i$ where $c_i\in \cI(\rho,\tau)$ and $r_i\in B(H_{\rho})$. Then we get
\be\label{VMV}
	(V\otimes I_{H_{\rho}})(M_S^*\otimes I_{H_{\rho}})(V\otimes I_{H_{\rho}})^*((\Sigma_{k=0}^{\infty}(R_k \otimes a)(\fw^*)^{(k)}c)\otimes rh)\ee
    \[\hspace{5cm}=\Sigma_i(\Sigma_{k=0}^{\infty}(R_k \otimes a)(\fw^*)^{(k)}c_i)\otimes r_ih).
\]
Therefore
\be\label{UVMV}
U(V\otimes I_{H_{\rho}})(M_S^*\otimes I_{H_{\rho}})(V\otimes I_{H_{\rho}})^*((\Sigma_{k=0}^{\infty}(R_k \otimes a)(\fw^*)^{(k)}c)\otimes rh)
\ee
\begingroup
\allowdisplaybreaks
\begin{align*}
\hspace{5cm}&=U\Sigma_i(\Sigma_{k=0}^{\infty}(R_k \otimes a)(\fw^*)^{(k)}c_i)\otimes r_ih)\\&=\Sigma_i(\Sigma_{k=0}^{\infty}(R_k \otimes a)(\fw^*)^{(k)}c_ir_ih)\\&=\Sigma_{k=0}^{\infty}(R_k \otimes a)(\fw^*)^{(k)}S(\fw)^*crh
\end{align*}
\endgroup
Applying $\varphi_{\infty}(d)$ to a map in $(E^{\rho})^{\otimes k}$ amounts to composing it with $I_{E^{\otimes k}}\otimes d$. Thus
\begingroup
\allowdisplaybreaks
\begin{align*}
	&(V\otimes I_{H_{\rho}})(M_S^*\otimes I_{H_{\rho}})(V\otimes I_{H_{\rho}})^*(\varphi_{\infty}(d)\otimes I_{H_{\rho}})((\Sigma_{k=0}^{\infty}(R_k \otimes a)(\fw^*)^{(k)}c)\otimes rh)\\&=(V\otimes I_{H_{\rho}})(M_S^*\otimes I_{H_{\rho}})(V\otimes I_{H_{\rho}})^*((\Sigma_{k=0}^{\infty}(R_k \otimes da)(\fw^*)^{(k)}c)\otimes rh)\\&=\Sigma_i(\Sigma_{k=0}^{\infty}(R_k \otimes da)(\fw^*)^{(k)}c_i)\otimes r_ih)\\&=(\varphi_{\infty}(d)\otimes I_{H_{\rho}})(V\otimes I_{H_{\rho}})(M_S^*\otimes I_{H_{\rho}})(V\otimes I_{H_{\rho}})^*((\Sigma_{k=0}^{\infty}(R_k \otimes a)(\fw^*)^{(k)}c)\otimes rh)
\end{align*} 
\endgroup
proving (\ref{commuted}).

To prove (\ref{commuteeta}) we need to compute $(W_{\eta}')^*((R_k \otimes a)(\fw^*)^{(k)}c)$ for $\eta \in E^{\rho}$, $k\geq 0$ and $g\in H_{\rho}$.
We use the notation set up in \cite[Theorem 4.5]{MS16}. In particular, for every $k\geq 1$, $R_k=(Z^{(k)})^{-1}$ is positive and invertible. We shall also need the notation of \cite[Lemma 7.4]{MS16}. Thus $Z'_k=Z^{(k)}(Z^{(k-1)}\otimes I_E)^{-1}=R_k^{-1}(R_{k-1}\otimes I_E)$ (thus, $(Z'_k)^*R_k=R_{k-1}\otimes I_E$) and $C_k\otimes I_H=U_k^*(Z'_k \otimes I_H)U_k$.
 We have
$(W_{\eta}')^*((R_k \otimes a)(\fw^*)^{(k)}c)=T_{\eta}^*C_k^*(((R_k \otimes a)(\fw^*)^{(k)}c)$.

We use $\eta$ for an element of $E^{\rho}$ but it will be convenient also to view $\eta$ as a map in $\cI(\rho,\rho^{E}\circ \varphi)$ and then, with $\fw$ as above and $a\in \cI(\tau,\rho)$,  $\eta^*(I_E\otimes a)\fw^*\in \cI(\tau,\rho)$.
We compute
\begingroup
\allowdisplaybreaks
\begin{align*}
	(W'_{\eta})^*\otimes &I_H)(\Sigma_k (R_k \otimes a)(\fw^*)^{(k)}c\otimes rh)\\&=\Sigma_k (I_{k-1}\otimes \eta)^*((Z'_k)^*\otimes I_H)(R_k \otimes a)(\fw^*)^{(k)}c\otimes rh \\&=\Sigma_k (I_{k-1}\otimes \eta)^*(R_{k-1}\otimes I_E\otimes a)(\fw^*)^{(k)}c\otimes rh\\&=\Sigma_k (R_{k-1}\otimes \eta^*(I_E\otimes a)(\fw)^*)(\fw^*)^{(k-1)}c\otimes rh.
\end{align*}
\endgroup
Thus
$$(W'_{\eta})^*\otimes I_H)(V\otimes I_H)(K_{\fw,a,cr}\otimes h)=(V\otimes I_H)(K_{\fw,\eta^*(I_E\otimes a)\fw^*,cr}\otimes h)$$ 
and 
\begingroup
\allowdisplaybreaks
\begin{align*}
	(W'_{\eta})^*\otimes &I_H)(V\otimes I_H)(M_S^*\otimes I_H)(K_{\fw,a,cr}\otimes h)\\&=(W'_{\eta})^*\otimes I_H)(V\otimes I_H)(K_{\fw,a,S(\fw)^*cr}\otimes h)\\&=(V\otimes I_H)(K_{\fw,\eta^*(I_E\otimes a)\fw^*,S(\fw)^*cr}\otimes h)\\&=(V\otimes I_H)(M_S^*\otimes I_H)(K_{\fw,\eta^*(I_E\otimes a)\fw^*,cr}\otimes h)\\&=(V\otimes I_H)(M_S^*\otimes I_H)(V^*\otimes I_H)((W'_{\eta})^*\otimes I_H)(V\otimes I_H)(K_{\fw,a,cr}\otimes h).
\end{align*}
\endgroup
This completes the proof of (\ref{commuteeta}).
\end{proof}

In the following we write $\mathcal{M}(K)$ for  $\{M_S : S\in M(K)\}$.

\begin{theorem}\label{isomalg}
Let $U$ be the unitary map from $\cF(E^{\rho})\otimes_{\rho(M)'}H_{\rho}$ onto $\cF(E)\otimes_{\rho}H_{\rho}$ as in \cite[Theorem 3.9]{MS04} (see also \cite[Section 7]{MS16}), let $V$ is the unitary map of Proposition~\ref{EK} and $\tilde{V}$ is as in Remark~\ref{directsum}.
Then
\begin{enumerate}
\item[(1)] For every $X\in H^{\infty}(E,Z)$, the map $\widehat{X}$ lies in $M(K)$ and
 \be\label{X}
 U\tilde{V}(M_{\widehat{X}}\otimes I_{H_{\rho}})\tilde{V}^*U^*=X\otimes I_{H_{\rho}}.
 \ee
In particular, given $\xi\in E$, $S_{\xi}(\fw):=\fw L_{\xi}$ ($=\widehat{W_{\xi}}(\fw)$) lies in $M(K)$ (where $L_{\xi}h=\xi\otimes h$ ) and we have
 \be\label{xi}
 U\tilde{V}(M_{S_{\xi}}\otimes I_{H_{\rho}})\tilde{V}^*U^*=W_{\xi}\otimes I_{H_{\rho}}.
 \ee
\item[(2)] $\mathcal{M}(K)$ is completely isometrically isomorphic to $H^{\infty}(E,Z)$ via an isomorphism that is $w^*$-continuous. In fact,
    \be\label{theisomorphism}
   U\tilde{V}(\mathcal{M}(K)\otimes I_{H_{\rho}})\tilde{V}^*U^*=H^{\infty}(E,Z)\otimes I_{H_{\rho}}.
    \ee
\end{enumerate}

\end{theorem}
\begin{proof}
We start by proving (\ref{xi}). In the following computation we use the equality $Z_k^*R_k=I_E\otimes R_{k-1}=R_kZ_k$ (\cite[Equation (4.5)]{MS16}) and, thus, $R_kZ_kZ_k^*R_k=I_k\otimes R_{k-1}^2$. Note that $S_{\xi}$ satisfies conditions (a) and (b) following Definition~\ref{multiplier} and one can define $M_S^*$ at least on a dense subspace of $E(K,\Sigma)$ by $M_{S_{\xi}}^*K_{\fw,a,c}=K_{\fw,a, S_{\xi}(\fw)^*c}$. To show that $M_{S_{\xi}}$ is bounded we need to show that there is a constant $C$ such that, for $\sigma \in \Sigma$, $\fz\in \mathbb{D}(X,\sigma)$ and a positive element $a^2\in \sigma(M)'$, $S_{\xi}(\fz)K(\fz,\fz)(a^2)S_{\xi}(\fz)^*\leq CK(\fz,\fz)(a^2)$. Now compute
\begingroup
\allowdisplaybreaks
\begin{align*}
	S_{\xi}(\fz)K(\fz,\fz)(a^2)S_{\xi}(\fz)^*&=\Sigma_k \fz L_{\xi}\fz^{(k)}(R_k^2\otimes a^2)\fz^{(k) *}L_{\xi}^*\fz^*=\Sigma_k \fz^{(k+1)}L_{\xi}(R_k^2\otimes a^2)L_{\xi}^*\fz^{(k+1) *}\\&=\Sigma_k \fz^{(k+1)}(I_E\otimes R_k \otimes a)L_{\xi}L_{\xi}^*(I_E \otimes R_k \otimes a)\fz^{(k+1) *}\\&\leq ||\xi||^2 \Sigma_k \fz^{(k+1)}(I_E\otimes R_k^2 \otimes a^2)(\fz^{(k+1)})^*\\&=||\xi||^2\Sigma_k \fz^{(k+1)}( R_{k+1}Z_{k+1}Z_{k+1}^*R_{k+1}\otimes a^2)(\fz^{(k+1)})^*\\&\leq ||\xi||^2\sup_k||Z_k||^2 (\Sigma_k \fz^{(k+1)}( R_{k+1}^2\otimes a)(\fz^{(k+1)})^*)\\&\leq ||\xi||^2\sup_k ||Z_k||^2 K(\fz,\fz)(a^2)
\end{align*} 
\endgroup
completing the proof that $M_{S_{\xi}}$ is a bounded multiplier. In fact, $||M_{S_{\xi}}||\leq ||\xi||\sup||Z_k||$.

For $\xi\in E$, $\fw\in \mathbb{D}(X,\tau)$, $a\in \cI(\tau,\rho)$, $c\in \cI(\rho,\tau)$ , $r\in B(H_{\rho})$ and $h\in H_{\rho}$, we have
\begingroup
\allowdisplaybreaks
\begin{align*}
	(W_{\xi}^*& \otimes I_{H_{\rho}})U((R_k\otimes a)(\fw^*)^{(k)}c\otimes rh)=(W_{\xi}^* \otimes I_{H_{\rho}})(R_k\otimes a)(\fw^*)^{(k)}crh)\\&=(T_{\xi}^*\otimes I_{H_{\rho}})(Z_k^*R_k\otimes a)(\fw^*)^{(k)}crh=(T_{\xi}^*\otimes I_{H_{\rho}})(I_E\otimes R_{k-1}\otimes a)(\fw^*)^{(k)}crh \\&=( R_{k-1}\otimes a)(\fw^*)^{(k-1)}L_{\xi}^*\fw^*crh =( R_{k-1}\otimes a)(\fw^*)^{(k-1)}S_{\xi}(\fw)^*crh.
\end{align*}
\endgroup
Thus \begin{align*}
	(W_{\xi}^* \otimes I_{H_{\rho}})U(\Sigma_k(R_k\otimes a)(\fw^*)^{(k)}c\otimes rh)&=\Sigma_k( R_{k-1}\otimes a)(\fw^*)^{(k-1)}S_{\xi}(\fw)^*crh \\&=\Sigma_k( R_{k}\otimes a)(\fw^*)^{(k)}S_{\xi}(\fw)^*crh .
\end{align*}
Using (\ref{UVMV}), this is equal to $$U(V\otimes I_{H_{\rho}})(M_{S_{\xi}}^*\otimes I_{H_{\rho}})(V\otimes I_{H_{\rho}})^*((\Sigma_{k=0}^{\infty}(R_k \otimes a)(\fw^*)^{(k)}c)\otimes rh)$$ and (\ref{xi}) follows.

We now turn to prove (\ref{X}) for $Y=\varphi_{\infty}(b)$ for $b\in M$. In this case $\tilde{Y}(\fz)=\sigma(b)$ (for $\fz\in \mathbb(X,\sigma)$. It is straightforward to check that $S(\fz)=\sigma(b)$ is a bounded multiplier (in fact, $||M_S||\leq ||b||$ in this case). To prove (\ref{X}) for this case, we use (\ref{UVMV}) and the fact that $S(\fw)^*cr=\tau(b)^*cr=c\rho(b)^*r$ (as $c\in \cI(\rho,\tau)$) to compute
\begingroup
\allowdisplaybreaks
\begin{align*}
	U(V\otimes I_{H_{\rho}})&(M_S^*\otimes I_{H_{\rho}})(V\otimes I_{H_{\rho}})^*((\Sigma_{k=0}^{\infty}(R_k \otimes a)(\fw^*)^{(k)}c)\otimes rh)
 \\&=\Sigma_{k=0}^{\infty}(R_k \otimes a)(\fw^*)^{(k)}S(\fw)^*crh=\Sigma_{k=0}^{\infty}(R_k \otimes a)(\fw^*)^{(k)}\tau(b)^*cr\\&=\Sigma_{k=0}^{\infty}(R_k \otimes a)(\varphi_k(b)^*\otimes I_{H_{\tau}})(\fw^*)^{(k)}cr\\&=\Sigma_{k=0}^{\infty}(\varphi_k(b)^*\otimes I_{H_{\tau}})(R_k \otimes a)(\fw^*)^{(k)}cr\\&=(\varphi_{\infty}(b)^*\otimes I_{H_{\tau}})\Sigma_{k=0}^{\infty}(R_k \otimes a)(\fw^*)^{(k)}cr\\&=(\varphi_{\infty}(b)^*\otimes I_{H_{\tau}})U((\Sigma_{k=0}^{\infty}(R_k \otimes a)(\fw^*)^{(k)}c)\otimes rh)
\end{align*}
\endgroup
proving (\ref{X}) for this case. Therefore, (\ref{X}) holds for every $Y$ in the algebra $H_0^{\infty}(E,Z)$ generated by $\{W'_{\eta}: \eta \in E^{\rho}\}\cup \varphi_{\infty}(M)$ which is $w^*$-dense in $H^{\infty}(E,Z)$.

 Given a $Y\in H^{\infty}(E,Z)$ there is a convegent net $Y_{\alpha}\rightarrow Y$ in the $w^*$-topology where each $Y_{\alpha}$ lies in $H_0^{\infty}(E,Z)$. By the Krein-Smulian Theorem we can assume that the net is bounded. To prove (\ref{X}) for $Y$ it suffices to show that $M_{\widehat{Y_{\alpha}}}$ converges weakly to $M_{\widehat{Y}}$. Since the net is bounded it is enough to check that $$\langle (M_{\widehat{Y_{\alpha}}})^*K_{\fw,a,c},K_{\fw',a',c'}\rangle \rightarrow^{w^*} \langle (M_{\widehat{Y}})^*K_{\fw,a,c},K_{\fw',a',c'}\rangle$$ for every $\fw,\fw',a,a',c,c'$. But $$\langle (M_{\widehat{Y_{\alpha}}})^*K_{\fw,a,c},K_{\fw',a',c'}\rangle=\langle K_{\fw,a,(\widehat{Y_{\alpha}}(\fw))^*c},K_{\fw',a',c'}\rangle= c^*\widehat{Y_{\alpha}}(\fw)K(\fw,\fw')(a^*a')c'$$
   and this converges to $c^*\widehat{Y}(\fw)K(\fw,\fw')(a^*a')c'$ in the $w^*$-topology. This completes the proof of (\ref{X}).

Combining this with Lemma~\ref{UVMSVU} completes the proof of part (2).
\end{proof}

Note that this generalizes part of \cite[Corollary 3.6]{SSS} and a version of this theorem, in the case of a fixed $\sigma$,  was proved by Good (\cite[Theorem 5.5]{G}).

\begin{corollary}\label{Sbounded}
For every $S\in M(K)$,
\begin{enumerate}
\item[(1)] there is $Y\in H^{\infty}(E,Z)$ such that, for every $\tau\in \Sigma$ and $\fw\in \mathbb{D}(X,\tau)$, $S(\fw)=\widehat{Y}(\fw)$ and
\item[(2)] $\sup\{||S(\fw)||: \fw\in \mathbb{D}(X,\tau),\; \tau\in \Sigma \}<\infty .$
\end{enumerate}
\end{corollary}
\begin{proof}
Clearly, (2) follows from (1). To prove (1) note that it follows from Theorem~\ref{isomalg} that there is some $Y\in H^{\infty}(E,Z)$ such that $M_S=M_{\widehat{Y}}$. So it is left to show that, if $S\in M(K)$ and $M_S=0$ then $S=0$. But if $M_S=0$ then, for every $\tau \in \Sigma$, $\fw\in \cU(\tau)$, $a\in \cI(\tau,\rho)$ and $c\in B(H_{\rho},H_{\tau})$, $K_{\fw,a,S(\fw)^*c}=0$. It suffices to show now that, if $K_{\fw,a,c}=0$ for all $\fw, a$, then $c=0$. Using (\ref{innerproduct}), we get $c^*K(\fw,\fw)(a^*a)c=0$ for all such $\fw,a$. Using Lemma~\ref{intertwiners}, we get $c^*K(\fw,\fw)(I)c=0$ for all $\fw$. For our kernel (\ref{Kweighted}) we have $I\leq K(\fw,\fw)(I)$ (as $I$ is the term for $k=0$) and, thus, $c=0$.

\end{proof}

\begin{definition}
We write $H^{\infty}(D_X)$ for the algebra of all matricial functions $f:D_X \rightarrow \mathcal{B}_H$ that are bounded.

\end{definition}

\begin{remark}\label{boundedmatricial}
Given $Y\in H^{\infty}(E,Z)$, we saw above that we can identify it with some $f\in H^{\infty}(D_X)$. In \cite[Theorem 3.1]{SSS} it is shown that, when $M=\mathbb{C}$ and $Z_k=I$ for all $k$, the converse also holds (i.e. every $f\in H^{\infty}(D_X)$ defines a bounded multiplier and, thus, some $Y\in H^{\infty}(E,Z)$). In our case here, to prove it, we need to assume that $\Sigma$ contains a universal generator $\sigma_0$ (see below). We may be able to get away with assuming that a natural sequence of subrepresentations of $\sigma_0$ is contained in $\Sigma$ (see the argument in \cite[Theorem 3.1]{SSS} ) but it does not seem to be worth the trouble.

\end{remark}

Write $\pi$ for an infinite multiple of a faithful normal representation of $M$ and let $\sigma_0=\pi^{\cF(E)}\circ \varphi_{\infty}$ (a normal representation of $M$ on $\cF(E)\otimes_{\pi}H$). 
Note that every $\sigma \in \widehat{M}$ is a subrepresentation of $\pi$. 
We refer to $\pi$ as the universal generator. 
Define ${\fs}_0$ by the formula ${\fs}_0(\xi\otimes h)=W_{\xi}h$ for $\xi\in E$ and $h\in \cF(E)\otimes_{\pi}H_{\pi}$. 
We can view ${\fs}_0$ as a map on $\Sigma_{k=1}^{\infty}\oplus E^{\otimes k}\otimes_{\pi}H_{\pi}$ defined by
\be\label{s0}
{\fs}_0=\Sigma_{k=1}^{\infty}\oplus Z_k\otimes I_{H_{\pi}}
\ee and
\be\label{s0txi}
{\fs}_0(T_{\xi}\otimes I_{H_{\pi}})=W_{\xi}\otimes I_{H_{\pi}}
\ee for $\xi\in E$.

Similarly, $${\fs}_0^{(m)}=\Sigma_{k=m}^{\infty}\oplus Z_k(I_E\otimes Z_{k-1})\cdots(I_{E^{\otimes m-1}}\otimes Z_{k-m+1}) \otimes I_{H_{\pi}}$$  on $\Sigma_{k=m}^{\infty} \oplus E^{\otimes k}\otimes_{\pi}H_{\pi}$.

The computation in \cite[Equation (4.8)]{MS16} shows that
 \be\label{eq:DX}
 \Sigma_{k=1}^{\infty} {\fs}_0^{(k)}(X_k\otimes I_{H_{\pi}}){\fs}_0^{(k) *} =I_{\cF(E)\otimes_{\pi}H_{\pi} }.
 \ee
Thus ${\fs}_0$ does not lie in $\mathbb{D}(X,\sigma_0)$ but it does define a representation of $H^{\infty}(E,Z)$. The representation it defines is $T\mapsto T\otimes_{\pi}I \in B(\cF(E)\otimes H_{\pi})$.

For $0<r<1$ it follows from (\ref{def:DX}) that $r{\fs}_0\in \mathbb{D}(X,\sigma_0)$.

\begin{theorem}\label{boundedmultipliers}
Suppose $\sigma_0\in \Sigma$. Then every $f\in H^{\infty}(D_X)$ lies in $M(K)$ (i.e. it is a bounded multiplier). Thus, we can identify $H^{\infty}(D_X)$ with $M(K)$.
\end{theorem}
\begin{proof}
Fix $f\in H^{\infty}(D_X)$ and $0<r<1$.
Suppose $D\in (\pi^{\cF(E)}(H^{\infty}(E,Z)))'$ . Then $D\in \sigma_0(M)'$ and $D$ commutes with $W_{\xi}\otimes I_{H_{\pi}} $ for all $\xi\in E$. Since ${\fs}_0(T_{\xi}\otimes I_{H_{\pi}})=W_{\xi}\otimes I_{H_{\pi}} $, $D$ commutes with ${\fs}_0(T_{\xi}\otimes I_{H_{\pi}})$ . Thus $D{\fs}_0(T_{\xi}\otimes_{\pi}I_{ H_{\pi}})={\fs}_0(T_{\xi}\otimes_{\pi}I_{ H_{\pi}})D={\fs}_0(I_E\otimes D)(T_{\xi}\otimes_{\pi}I_{ H_{\pi}})$ for all $\xi\in E$. It follows that $D{\fs}_0={\fs}_0(I_E\otimes D)$ and, thus, $D(r{\fs}_0)=(r{\fs}_0)(I_E\otimes D)$.

Since $f$ is matricial, we get $f(r{\fs}_0)D=Df(r{\fs}_0)$ for every $D\in (\pi^{\cF(E)}(H^{\infty}(E,Z)))'$. Thus (\cite{G}) $f(r{\fs}_0)\in \pi^{\cF(E)}(H^{\infty}(E,Z))$. Thus, there is some $Y_r\in H^{\infty}(E,Z)$ such that $f(r{\fs}_0)= \pi^{\cF(E)}(Y_r)$. It follows that $||Y_r||=||\pi^{\cF(E)}(Y_r)||\leq \sup||f(\fz)||<\infty$ so that $\{Y_r : 0<r<1 \}$ is bounded.

Now fix $\fz\in r\mathbb{D}(X,\sigma)$ so that $\frac{1}{r}\fz \in \mathbb{D}(X,\sigma)$. From \cite[Lemmas 5.4 and 5.5]{MS16} we see that there is an isometry $V_r: H\rightarrow \cF(E)\otimes_{\sigma}H_{\sigma}$ (written $K(\frac{1}{r}\fz)$ there) such that, for $\xi\in E$,
\be\label{Vr}
V_r (\frac{1}{r}L_{\xi}^*\fz ^*)=(W_{\xi}^*\otimes I_{H_{\sigma}})V_r
\ee  and, for $R\in H^{\infty}(E,Z)$,
\be
V_r \widehat{R}(\frac{1}{r}\fz)V_r^*=\sigma^{\cF(E)}(R) .
\ee
Since $\sigma$ is a subrepresentation of $\pi$, there is an isometry $w: H_{\sigma}\to  H_{\pi}$ such that
\[\sigma(a)=w^*\pi(a)w\] for $a\in M$. It follows that $v:\cF(E)\otimes_{\sigma}H_{\sigma} \rightarrow \cF(E)\otimes_{\pi}H_{\pi}$ defined by $v(x\otimes h)=x \otimes wh$ for $x\in \cF(E)$ and $h\in H_{\sigma}$, is an isometry and we have
$$v(W_{\xi}^*\otimes I_{\Hs})(x\otimes h)=v(T_{\xi}^*x\otimes h)=T_{\xi}^*x \otimes wh=(W_{\xi}^*\otimes I_{H_{\pi}})v(x\otimes h)$$ and, thus
\begin{equation}\label{v}
v(W_{\xi}^*\otimes I_{\Hs})=(W_{\xi}^*\otimes I_{H_{\pi}})v\quad(\text{for all} \,\xi\in E).
\end{equation}
Write $v_r=vV_r$. Then $v_r: H_{\sigma}\to \cF(E)\otimes H_{\pi}$ and, applying (\ref{Vr}) and (\ref{v}), 
\[
v_r(\frac{1}{r}T_{\xi}^*\fz^*)=(W_{\xi}^*\otimes I_{H_{\pi}})v_r\quad(\forall\,\xi\in E).
\] 
Taking adjoint in the both sides,
\begin{equation}\label{vr}
	\fz T_{\xi}v_r^*=rv_r^*(W_{\xi}\otimes I_{H_{\pi}})\quad(\forall\,\xi\in E)
\end{equation}
Observe that, for all $\eta\in\cF(E)$ and $h_{\pi}\in H_{\pi}$,
\[
T_{\xi}v_r^*(\eta\otimes h_\pi)=\xi\otimes v_r^*(\eta\otimes h_\pi)=(I_E\otimes v_r^*)(\xi\otimes\eta\otimes h_{\pi})=(I_E\otimes v_r^*)(T_{\xi}\otimes I_{H_{\pi}})(\eta\otimes h_\pi),
\]
that is,
\[
T_{\xi}v_r^*=(I_E\otimes v_r^*)(T_{\xi}\otimes I_{H_{\pi}}).
\]
Therefore, it follows from (\ref{vr}) that,
\be
	\fz(I_E \otimes v_r^*) (T_{\xi}\otimes I_{H_{\pi}})=r v_r^*(W_{\xi}\otimes I_{H_{\pi}})=r v_r^*\fs_0(T_{\xi}\otimes I_{H_{\pi}})\quad(\forall\,\xi\in E)
\ee
Hence, for all $0<r< 1$,
\[
\fz(I_E \otimes v_r^*)=v_r^*(r\fs_0).
\]
Since $f$ is matricial, 
\[
f(\fz) v_r^*=v_r^*f(r\fs_0)
\]
Thus, for $\fz\in r\mathbb{D}(X,\sigma)$,
$$f(\fz)=v_r^*f(r{\fs}_0)v_r=v_r^*\pi^{\cF(E)}(Y_r)v_r=V_r^*\sigma^{\cF(E)}(Y_r)V_r=\widehat{Y}_r(\fz).$$
Therefore, $f=\widehat{Y}_r$ on $r\mathbb{D}(X,\sigma)$ for $0<r<1$.
As shown above, the set $\{Y_r :0<r<1\}$ is bounded and we can fix a sequence $\{r_n\}$ increasing to $1$ and find a subnet $\{Y_{r_{\alpha}}\}$ of it such that $Y_{r_{\alpha}}\rightarrow Y$ in the $w^*$-topology for some $Y\in H^{\infty}(E,Z)$.

Now fix $\fz\in \mathbb{D}(X,\sigma)$. For some $0<r_0<1$, we have $\fz\in r\mathbb{D}(X,\sigma)$ for all $r>r_0$. Thus $f(\fz)=\widehat{Y}_r(\fz)$ for every $r>r_0$ and, in the limit, $f(\fz)=\widehat{Y}(\fz)$. Thus $f=\widehat{Y}$ and, by Theorem\ref{isomalg} (1), $f\in M(K)$.
\end{proof}


\section{Kernels defined on the matricial unit ball}\label{kernel on matricial ball}
We now discuss the same kernel as in the previous section except that it is defined on the matricial unit ball
\[ \mathbb{D}(E)=\{\mathbb{D}(E,\sigma) : \sigma\in \Sigma \} \]
where
\[ \mathbb{D}(E,\sigma)=\{\fz\in \cI(\sigma^E\circ \varphi, \sigma): ||\fz||<1 \} \]
where $E$ is a fixed $W^*$-correspondence over $M$.

The kernel we study here can be viewed as a non commutative analogue of the unitarily invariant kernels discussed in \cite{CH}.
Again, to define it, we start with  sequence $\{R_k\}_{k=0}^{\infty}$ satisfying
\begin{enumerate}
\item[(1)] For each $k$, $R_k$ is a positive, invertible operator in $\varphi_k(M)^c$ ($:=\cL(E^{\otimes k})\cap \varphi_k(M)'$) and $R_0=I$.
\item[(2)] $\sup_k ||R_k^{-1}(I_E\otimes R_{k-1})||<\infty$.
\item[(3)] $\limsup_k ||R_k^2||^{\frac{1}{k}}\leq 1$

\end{enumerate}

\begin{definition}\label{DX1}
Define the sequence $Z=\{Z_k\}_{k\geq 0}$ by $$Z_k=R_k^{-1}(I_E\otimes R_{k-1})\;\;\;, k\geq 1$$ (and $Z_0=I$ ) and write $$Z^{(k)}:=Z_k(I_E\otimes Z_{k-1})\cdots (I_{E^{\otimes (k-1)}}\otimes Z_1)=R_k^{-1}.$$

\end{definition}

\begin{proposition}\label{kernel1}
\begin{enumerate}
\item[(i)] Given $\sigma, \tau \in \widehat{M}$, for $\fz\in \mathbb{D}(E,\sigma)$ and $\fw \in \mathbb{D}(E,\tau)$, the map
\be \label{Kweighted1}
K(\fz,\fw)(a)=\sum_{k=0}^{\infty} \fz^{(k)}(R_k^2 \otimes a)\fw^{(k) *}
 \ee
is well defined (in the sense of norm convergence) and maps $\cI(\tau,\sigma)$ into itself.
\item[(ii)] $\mathbb{D}=\{\mathbb{D}(E,\sigma)\}_{\sigma}$ is a matricial family of sets for $(M,E)$ and
$$ K:\mathbb{D}\times \mathbb{D} \rightarrow \coprod_{\sigma,\tau \in \widehat{M}} B_{*}(\cI(\tau,\sigma),\cI(\tau,\sigma))$$  $$\subseteq \coprod_{\sigma,\tau \in \widehat{M}} B_{*}(\cI(\tau,\sigma),B(H_{\tau},H_{\sigma}))$$ is a cp matricial kernel.
\item[(iii)]  For all $k$, $Z_k\in \varphi_k(M)^c$, $Z_k$ is invertible and $Z^{(k) *}Z^{(k)}=R_k^{-2}$. Also, $Z_0=I$ and $\sup||Z_k||<\infty$.
\end{enumerate}

\end{proposition}
\begin{proof}
Part (i) follows from the assumption that $\limsup_k ||R_k^2||^{\frac{1}{k}}\leq 1$ and, as $||\fz||<q$ and $||\fw||<q$ for some $q<1$, $||\fz^{(k)}||||\fw^{(k)}||<q^{2k}$.
 
(ii) It is straightforward to check that $\mathbb{D}=\{\mathbb{D}(E,\sigma)\}_{\sigma}$ is a matricial family of sets for $(M,E)$. The rest is also easy to check noting that, if $C\fz=\tilde{\fz}(I_E\otimes C)$ (as in Definition~\ref{kernel} (3)), then, for $k\geq 1$, $C\fz^{(k)}=\tilde{\fz}^{(k)}(I_{E^{\otimes k}}\otimes C)$.
\end{proof}

\begin{proposition}\label{EK1}
Assume $\Sigma \subseteq \widehat{M}$ is the subset of all faithful representations, $\rho\in \Sigma$ is fixed, $K$ is given in (\ref{Kweighted1}) and $E(K,\Sigma)$ is the associated correspondence. Recall that $E^{\rho}$ is the $W^*$-corresponcence $\cI(\rho,\rho^E \circ \varphi)$.
Then there is an isomorphism of $W^*$-correspondences
$$V: E(K,\Sigma) \rightarrow \mathcal{F}(E^{\rho})\otimes_{\rho(M)'}B(H_{\rho}) $$
such that $K_{\fw,a,cr}$ (for $\fw \in \mathbb{D}(E,\tau)\subseteq \cI(\tau^E\circ \varphi , \tau)$, $a\in \cI(\tau,\rho)$, $c\in \cI(\rho,\tau)$ and $r\in B(H_{\rho})$) is mapped into $(\Sigma_{k=0}^{\infty}(R_k \otimes a)(\fw^*)^{(k)}c)\otimes r$.

\end{proposition}
\begin{proof}
Let $L=\text{span}\{K_{\fw,a,cr}: \fw \in \mathbb{D}(E,\tau), a\in \cI(\tau,\rho), c\in \cI(\rho,\tau), r\in B(H_{\rho})\}$ and let $V_0$ be the linear map from $L$ into $\mathcal{F}(E^{\rho})\otimes_{\rho(M)'}B(H_{\rho})$ that maps $K_{\fw,a,cr}$ to $(\Sigma_{k=0}^{\infty}(R_k \otimes a)(\fw^*)^{(k)}c)\otimes r$. It is straightforward to check that $U_0$ is a bimodule map.
We first compute
\begingroup
\allowdisplaybreaks
\begin{align*}
	\langle K_{\fw,a,cr},K_{\fw',a',c'r'}\rangle &= r^*c^*K(\fw,\fw')(a^*a')c'r'=\sum_{k=0}^{\infty} r^*c^*\fw^{(k)}(R_k^2 \otimes a^*a')\fw'^{(k) *}c'r'\\&=\sum_{k=0}^{\infty} ((R_k \otimes a)\fw^{(k) *}cr)^*((R_k \otimes a')\fw'^{(k) *}c'r').
\end{align*}
\endgroup
For a fixed $k$, $(R_k \otimes a)\fw^{(k) *}c \in (E^{\rho})^{\otimes k}=\cI(\rho, \rho^{E^{\otimes k}}\circ \varphi )$ since $c\in \cI(\rho,\tau)$, $\fw^{(k) *}\in \cI(\tau^{E^{\otimes k}}\circ \varphi,\tau)$ and $R_k\in \varphi(M)'$. Thus $\sum_{k=0}^{\infty}(R_k \otimes a)\fw^{(k) *}c \in \mathcal{F}(E^{\rho})$. We, thus, have
$$\langle K_{\fw,a,cr},K_{\fw',a',c'r'}\rangle =\langle (\sum_{k=0}^{\infty}(R_k \otimes a)\fw^{(k) *}c)\otimes r,(\sum_{k=0}^{\infty}(R_k \otimes a')\fw'^{(k) *}c')\otimes r'\rangle.$$
So the map $V_0$, defined above on $L$ , preserves inner products.

It is left to show that $L$ is $w^*$-dense in $E(K,\Sigma)$ and the image of $U_0$ is $w^*$-dense in $\mathcal{F}(E^{\rho})\otimes_{\rho(M)'}B(H_{\rho})$.

The fact that $L$ is $w^*$-dense in $E(K,\Sigma)$ follows from the definition of $E(K,\Sigma)$ and the fact that, for $\fw\in \mathbb{D}(E,\tau)$, $a\in\cI(\tau,\rho)$ and $b\in B(H_{\rho},H_{\tau})$, $K_{\fw,a,b}$ lies in the $w^*$-closure of $L$ since $B(H_{\rho},H_{\tau})=\overline{\cI(\rho,\tau)B(H_{\rho})}^{w^*}$ (Lemma~\ref{corintertwiners}(1)).

To show that the image of $V_0$ is $w^*$-dense in $\mathcal{F}(E^{\rho})\otimes_{\rho(M)'}B(H_{\rho})$ it suffices to prove that $C:=\text{span}\{\Sigma_{k=0}^{\infty}(R_k \otimes a)(\fw^*)^{(k)}c\}$ is $w^*$-dense in $\mathcal{F}(E^{\rho})$. But $\overline{C}^{w^*}$ is a sub-bimodule of $\mathcal{F}(E^{\rho})$ and is invariant under the action of the gauge group (just replace $\fw$ by $\lambda \fw$ for $|\lambda|=1$). Therefore it suffices to show that, for every $k\geq 0$, the set $C_k:=\text{span}\{(R_k \otimes a)(\fw^*)^{(k)}c : \tau\in \Sigma, \fw\in \cI(\tau^{E}\circ \varphi,\tau),||\fw||<1, a\in \cI(\tau,\rho), c\in \cI(\rho,\tau)\}$ is $w^*$-dense in $(E^{\rho})^{\otimes k }$ which can be identified with $\cI(\rho,\rho^{E^{\otimes k}}\circ \varphi)$.

But, since $R_k\otimes I_{H_{\rho}}$ is invertible, it is enough to show that $C^0_k:=\text{span}\{(I_{E^{\otimes k}} \otimes a)(\fw^*)^{(k)}c : \tau\in \Sigma, \fw\in \cI(\tau^{E}\circ \varphi,\tau), a\in \cI(\tau,\rho), c\in \cI(\rho,\tau)\}$ is $w^*$-dense in $(R_k^{-1}\otimes I)(E^{\rho})^{\otimes k }$. Since $R_k\otimes I_{H_{\rho}}$ is invertible and commutes with $\varphi_{E^{\otimes k}}(M)\otimes I_{H_{\rho}}=(\rho^{E^{\otimes k}}\circ \varphi)(M)$, and since we identify $(E^{\rho})^{\otimes k }$ with $\cI(\rho,\rho^{E^{\otimes k}}\circ \varphi)$, we have $(R_k^{-1}\otimes I)(E^{\rho})^{\otimes k }=(E^{\rho})^{\otimes k }$. Therefore, we should show that $C^0_k$ is $w^*$-dense in $(E^{\rho})^{\otimes k }$.

But $(E^{\rho})^{\otimes k }$ is generated by elements of the form $\eta_k \otimes \eta_{k-1}\otimes \cdots \otimes \eta_1$ for $\eta_i\in E^{\rho}$ (identified with $\cI(\rho,\rho^{E}\circ \varphi)$). Every such element is being identified with the element $(I_{E^{\otimes k-1}}\otimes \eta_k)\cdots (I_{E^{\otimes 2}}\otimes \eta_3)(I_{E}\otimes \eta_2)\eta_1 \in \cI(\rho,\rho^{E^{\otimes k}}\circ \varphi)$. So, it suffices to show that $(I_{E^{\otimes k-1}}\otimes \eta_k)\cdots (I_{E^{\otimes 2}}\otimes \eta_3)(I_{E}\otimes \eta_2)\eta_1 \in C^0_k$.

For this, let $\tau$ be $k\rho$ and note that $\cI(\tau,\tau^E\circ \varphi)$ is isomorphic to $M_k(\cI(\rho,\rho^E\circ \varphi))$ so that we can define $\fw^* \in \cI(\tau,\tau^E\circ \varphi)$ as a $k\times k$ matrix over $\cI(\rho,\rho^E\circ \varphi)$ by
\be
\fw^*=\left( \begin{array}{cccccc} 0 & \eta_k & 0 & 0 & 0 & 0 \\ 0 & 0 & \eta_{k-1} & 0 & 0 & 0 \\0&0&0&\cdot&0&0 \\0 & 0 & 0 & 0 & \cdot & 0  \\ 0 & 0 & 0 & 0 & 0 & \cdot \\ \eta_1 & 0&\cdot &\cdot &\cdot & 0 \end{array}\right) .
\ee
Now $(\fw^*)^{(k)}=(I_{E^{\otimes k-1}}\otimes \fw^*)\cdots (I_{E^{\otimes 2}}\otimes \fw^*)(I_{E}\otimes \fw^*)\fw^*$ which is a $k\times k$ matrix with $(I_{E^{\otimes k-1}}\otimes \eta_k)\cdots (I_{E^{\otimes 2}}\otimes \eta_3)(I_{E}\otimes \eta_2)\eta_1$ in the $(1,1)$ entry. Note that elements of  $\cI(\tau,\rho)=\cI(k\rho,\rho)$ can be written as rows of length $k$ over $\rho(M)'$ and elements of $\cI(\rho,\tau)=\cI(\rho,k\rho)$ can be written as columns of length $k$ over $\rho(M)'$. Thus we can set $a=(I,0,\ldots,0)\in \cI(\tau,\rho)$ and $c=a^t\in \cI(\rho,\tau)$ and we get $(I_{E^{\otimes k-1}}\otimes \eta_k)\cdots (I_{E^{\otimes 2}}\otimes \eta_3)(I_{E}\otimes \eta_2)\eta_1 =(I_{E^{\otimes k}} \otimes a)(\fw^*)^{(k)}c$. Write $t=||\fw||$,  fix $r>t$ and write $\fw_0=\frac{1}{r}\fw$ and, thus, $(\fw^*)^{(k)}=r^k(\fw_0^*)^{(k)}$ and $(I_{E^{\otimes k}} \otimes a)(\fw^*)^{(k)}c=r^k(I_{E^{\otimes k}} \otimes a)(\fw_0^*)^{(k)}c\in C^0_k$.
\end{proof}

\section{Complete Pick property of matricial Szeg\"o kernel}\label{sec:CPSzego}
Let $M$ be a $W^*$-algebra. We denote the set of all faithful normal representations of $M$ as $F\widehat{M}$. 

  Let $K$ be a matricial cp kernel defined on a matricial family  $\cU=\{\cU(\sigma)\}_{\sigma\in FN Rep(M)}$. Suppose $\cE$ be a finite dimensional Hilbert space, let $\Omega\subseteq\cU$ and let
    \[
    f_0:\Omega\to\coprod_{\sigma\in F\widehat{M}}\cB(\Hs)\otimes\cB(\cE)\quad \left(\cong\coprod_{\sigma\in F\widehat{M}}\cB(\Hs\otimes\cE)\right)
    \]
be a matricial function. Then the complete Pick property of a kernel $K$ is the equivalence of the following two assertions:
\begin{itemize}
    \item[(I)] If $f_0$ extends to a matricial function defined on the full matricial  envelope $\Omega'$ of $\Omega$, then there exists $f\in \text{Mult}(E(K)\otimes\cE)$ such that $f|_{\Omega}=f_0$ and  $\|M_f\|_{\text{Mult}(E(K)\otimes\cE)}\leq 1$

    \item[(II)] The kernel \[
K_{f_0}(\fz,\fw)(a)=K(\fz,\fw)(a)\otimes I_{\cE}-f_0(\fz)\big[K(\fz,\fw)(a)\otimes I_{\cE}\big]f_0(\fw)^*
\]
is a cp matricial kernel on $\Omega$.
\end{itemize} 

\noindent\underline{Note:} $E(K)\otimes\cE$ is the $\rho(M)'-\cB(\Hr\otimes\cE)$ reproducing kernel $W^*$-correspondence consisting of matricial functions
\[
g:\cU\to\coprod_{\sigma\in F\widehat{M}}\cB_*\big(\cI(\rho,\sigma), \cB(\Hr\otimes\cE,\Hs\otimes\cE)\big).
\]
\begin{remark}
    For any contractive multiplier $f\in \text{Mult}(E(K)\otimes\cE)$ if $f_0=f|_{\Omega}$ then by Proposition \ref{SK} $K_{f_0}$ is a cp matricial kernel.
\end{remark}

The above remark suggests that, to prove that a matricial kernel $K$ is a complete Pick kernel it is enough to prove that (II) implies (I).

In this section we study the complete Pick property of the kernel defined in Proposition  \ref{weighted Szego} but without a non trivial weight sequence, that is, $R_k=I_{E^{\otimes k}}$ for all $k\in\bZ_+$. In this case, $X_1=I_E$ and $X_k=0$ for all $k\geq 2$ consequently $\bD(X,\sigma)$ turns into the unit ball $\bD(E,\sigma)$ as in Section \ref{kernel on matricial ball}. More precisely,  the kernel we are studying here is a matricial cp kernel on the matricial unit ball $\bD:=\{\bD(E,\sigma)\}_{\sigma\in F\widehat{M}}$ defined by 
\be \label{Szego}
\bK(\fz,\fw)(a)=\sum_{k=0}^{\infty} \fz^{(k)}(I_{E^{\otimes k}} \otimes a)\fw^{(k) *}\quad(a\in\cI(\tau,\sigma))
 \ee
where $\fz\in \mathbb{D}(E,\sigma)$ and $\fw \in \mathbb{D}(E,\tau)$ for $\sigma, \tau \in \widehat{M}$. We call this kernel as matricial Szeg\"o kernel.

Our goal now is to prove that matricial Szeg\"o kernel is a complete Pick Kernel, Theorem \ref{complete Pick}. Before going into the proof which includes several lemmas, we recall a few facts which will be used in the sequel. 

We connect this discussion with the discussion in \cite{MS13} that defines $H^{\infty}(E)$. The \emph{algebraic tensor algebra of $E$} is defined to be \[
\cT_0(E):=\text{alg}\{\phi_{\infty}(M), L_\eta:\eta\in E\}\subseteq\cL(\cF(E)),
\]
where, recall, $L_{\eta}$ is the \emph{creation operator} defined on $\cF(E)$ by the formula \begin{equation*}
    L_{\eta}\xi:=\eta \otimes \xi, \qquad \xi\in \cF(E).
\end{equation*} Thus $\cT_0(E)$ may be viewed as the subalgebra of $\cL(\cF(E))$ freely generated by $M$ and $E$.  The \emph{Hardy algebra of $E$}, $H^{\infty}(E)$, then is defined to be the ultraweak closure of $\cT_0(E)$ in $\cL(\cF(E))$. In symbols,
\[
H^{\infty}(E):=\overline{\text{alg}\{\phi_{\infty}(M), L_\eta:\eta\in E\}}^{uw}\subseteq\cL(\cF(E)).
\]
Each $\sigma\in \widehat{M}$ gives an induced representation, $\sigma^{\cL(F(E))}$, of $\cL(F(E))$ on the \emph{Hilbert space} $\cF(E)\otimes_M H_{\sigma}$. It is defined by the formula \be\sigma^{\cF(E)}(T)=T\otimes I_{H_{\sigma}},\qquad T\in \cL(F(E)).\ee (cf. \eqref{eq:inducedrep})
Since $H^{\infty}(E)\subseteq\cL(\cF(E))$, we may restrict $\sigma^{(\cF(E))}$ to $H^{\infty}$ and refer to it as the representation of $H^{\infty}$ on $\cL(E)\otimes_M H_{\sigma}$ induced by $\sigma^{\cF(E)}$. In symbols, 
\[
\sigma^{\cF(E)}:H^{\infty}(E)\to\cB(\cF(E)\otimes\Hs).
\]
The representation $\sigma^{\cF(E)} $leads to an important class of representations of $H^{\infty}(E)$ on the spaces $H_{\sigma}$. They come about in the following way.   If $\fz_{\sigma}\in B(E\otimes H_{\sigma},H_{\sigma})$ and $\|\fz_\sigma\|<1$ then there is a representation $\sigma\times\fz_\sigma: H^{\infty}(E)\to\cB(\Hs)$ defined on generators by
\begin{enumerate}
    \item $\sigma\times\fz_\sigma (\phi_{\infty}(a)):=\sigma(a),\qquad a\in M$
    \item $\sigma\times\fz_\sigma (L_\eta)(h):=\fz_\sigma(\eta\otimes h),\qquad \eta\in E, h\in\Hs$.
\end{enumerate}
Recall from \cite{MS13} that for $F\in H^{\infty}(E)$, the $\sigma$-Berezin transform $\h{F}_\sigma: \bD(E,\sigma)\to\cB(\Hs)$ is defined by
\[
\h{F}_\sigma(\zs):=\sigma\times\zs(F).
\]
We need the following lemmas and proposition to prove the complete Pick property of the matricial Szeg\"o kernel $\bK$. Recall that $\cE$ is a finite dimensional Hilbert space that is fixed for our discussion.
\begin{lemma}\label{lem:W-in-alg}
   Suppose $\fz\in\bD(E,\sigma)$, $W\in\cB(\Hs\otimes\cE)$ and $\Omega=\{\fz\}$. Then $W\in\{\h{F}_\sigma(\fz):F\in\cT_0(E)\}\otimes\cB(\cE)$ if and only if the map $\fz\mapsto W$ extends to a matricial function on $\Omega'$, the matricial envelope of $\Omega$.
\end{lemma}
\begin{proof}
   Let $W=\sum_i\h{F}_i(\fz)\otimes A_i$ for some $F_i\in\cT_0(E)$ and $A_i\in\cB(\cE)$. Since $\fw\mapsto\sum_i\h{F}_i(\fw)\otimes A_i$ is a matricial function on $\bD$, restricting it to the set $\Omega'$ gives us the required function which extends the map $\fz\mapsto W$. This proves the necessary part.
    
    To prove the sufficient part, we assume that $f$ is the matricial function on $\Omega'$ that extends the map $\fz\mapsto W$.  We first claim the following:
    
    \noindent\underline{Claim:} If $\cQ\subseteq \Hs\otimes\cE$ is an invariant subspace for $\{\h{F}_\sigma(\fz):F\in\cT_0(E)\}\otimes\cB(\cE)$, then it is also an invariant subspace for $W$.

    \noindent\underline{Proof of the claim:} We denote
    \[
    \cT_{00}(E):=\{\phi_{\infty}(M), L_\eta:\eta\in E\}.
    \]
    The space $\cQ$ is of the form $\cQ_\sigma\otimes\cE$ where $\cQ_\sigma$ is an invariant subspace for $\big\{\h{F}_\sigma(\fz) : F\in\cT_{00}(E)\big\}$, that is,
    \[
    \sigma(a)\cQ_\sigma=\sigma\times\fz(\phi_{\infty}(a))\cQ_\sigma\subseteq\cQ_\sigma\quad(\forall\,a\in M)
    \]
    and
    \[
    \big(\h{L}_\eta\big)_\sigma(\fz)(\cQ_\sigma)\subseteq\cQ_\sigma\quad(\forall\,\eta\in E).
    \]
This implies, for each $\eta\in E$,
\[
\fz(\eta\otimes h)\in\cQ_{\sigma},\quad\forall h\in\cQ_\sigma.
\]
So, $\fz$ on $E\otimes\cQ_\sigma$ can be written in the form,
\[
\fz=\begin{bmatrix}
    \underline{\fz} & *\\
    0 & *
\end{bmatrix}
\]
with respect to the decomposition $(E\otimes\cQ_\sigma)\oplus(E\otimes[\Hs\ominus\cQ_\sigma])$. Then $\fz \bI=\bI\underline{\fz}$ where $\bI=\begin{bmatrix}
    I\\ 0
\end{bmatrix}$. Therefore, $\underline{\fz}\in\Omega'$. Thus by assumption, $f$ extends to $\underline{\fz}$. Since $f$ is a matricial function it respects intertwiners, so
\[
W\bI=f(\fz)\bI=\bI f(\underline{\fz}).
\]
Therefore, 
\[
W=\begin{bmatrix}
    f(\underline{\fz}) & *\\
    0 & *
\end{bmatrix}
\]
on $\Hs\otimes\cE=(\cQ_\sigma\otimes\cE)\oplus([\Hs\ominus\cQ_\sigma]\otimes\cE)$. Hence, $\cQ=\cQ_\sigma\otimes\cE$ is invariant under $W$. This proves the claim.

Now, if possible assume that, $W\notin \{\h{F}_\sigma(\fz):F\in\cT_0(E)\}\otimes\cB(\cE)$. Then there exists a linear functional $\Psi:\cB(\Hs)\otimes\cB(\cE)\to\bC$ such that $\Psi\big(\h{F}_\sigma(\fz)\otimes D\big)=0$ for all $F\in\cT_0(E)$, $D\in\cB(\cE)$ and $\Psi(W)\neq 0$. We can find some $n\in\bN$ and $u,v\in(\Hs\otimes\cE)^n$ such that 
\[
\Psi(X)=\langle X^{(n)}u,v\rangle\quad \forall X\in\cB(\Hs)\otimes\cB(\cE),
\]
where $X^{(n)}$ denotes the $n$-fold ampliation of $X$, that is, $X^{(n)}=X\otimes\cdots\otimes X$ ($n$-times). Denote,
\[
Q=\big[\big\{\h{F}_\sigma(\fz)^{(n)}:F\in\cT_0(E)\big\}\otimes\cB(\cE)\big]u
\]
Then for all $F\in\cT_0(E)$ and $D\in\cB(\cE)$,
\[
\left\langle \big[\h{F}_\sigma(\fz)^{(n)}\otimes D\big]u,v\right\rangle=\Psi(\h{F}_\sigma(\fz)\otimes D)=0.
\]
This implies $Q$ is a non-trivial invariant subspace for $\big\{\h{F}_\sigma(\fz)^{(n)}:F\in\cT_0(E)\big\}\otimes\cB(\cE)$ and $v\in Q^{\perp}$. Since $\fz\in\bD(E,\sigma)$, $\fz^{(n)}\in\bD(n\sigma)$ and $\h{F}_\sigma(\fz)^{(n)}=\h{F}_{n\sigma}\big(\fz^{(n)}\big)$. Apply the above claim for $\fz^{(n)}$, $W^{(n)}$ and for the representation $n\sigma$, we get that $Q$ is invariant under $W^{(n)}$. So,
\[
\Psi(W)=\langle W^{(n)}u,v\rangle=0
\]
which is a contradiction. This completes the proof of the lemma.
\end{proof}
For an ideal $\fJ\subseteq H^{\infty}(E)$ we denote \[
\h{\fJ}:=\{\h{F}:F\in\fJ\}.
\]
We need the following Davidson-Pitt \cite{DP} type distance formula that we prove using the distance formula that appears in the thesis of M. Gurevich \cite{Max} and J.R. Meyer \cite{Meyer}.
\begin{proposition}\label{prop:basic-distance-formula}
Let $\fJ\subseteq H^{\infty}(E)$ be a closed ideal and $\cG=\overline{\h\fJ E\big(\bK\big)\otimes_{\cB(\Hr)}\Hr}$. For all $X\in H^{\infty}(E)$ we have\[
\text{dist}\,(X,\fJ)=\left\|P_{\cG^{\perp}}\big(M_{\h{X}}\otimes I_{\Hr}\big)P_{\cG^{\perp}}\right\|.
\] 
\end{proposition}
\begin{proof}
    Let $\rho^{\cF(E)}:H^{\infty}(E)\to\cB(\cF(E)\otimes\Hr)$ be the induced representation. Denote $\tilde{G}=\overline{\rho^{\cF(E)}(\fJ)\big(\cF(E)\otimes\Hr\big)}$. By the distance formula in Maxim's thesis \cite[Theorem 4.1.9]{Max},
    \[
    \text{dist}(X,\fJ)=\left\|P_{\tilde{G}^{\perp}}\big(\rho^{\cF(E)}X\big)\right\|=\|P_{\tilde{G}^{\perp}}(X\otimes I_{\Hr})\|.
    \]
    By Theorem \ref{isomalg},
    \begin{align*}
        \|P_{\tilde{G}^{\perp}}(X\otimes I_{\Hr})\|=&\big\|P_{\tilde{G}^{\perp}}U\tilde{V}\big(M_{\tilde{X}}\otimes I_{\Hr}\big)\tilde{V}^*U^*\big\|
    \end{align*}
    where $U:\cF(E^{\rho})\otimes_{\rho(M)'}\Hr\to\cF(E)\otimes_{\rho}\Hr$ is the unitary map as in \cite[Theorem 3.9]{MS04} and $\tilde{V}: E(\bK)\otimes_{\cB(\Hr)}\Hr\to \cF(E^{\rho})\otimes_{\rho(M)'}\Hr$ is the unitary map as in Remark \ref{directsum}. So, $U\tilde{V}:E(\bK)\otimes_{\cB(\Hr)}\Hr\to \cF(E)\otimes_{\rho}\Hr$ is a unitary operator. Observe that $\tilde{V}^*U^*(\cF(E)\otimes_{\rho}\Hr)=E(\bK)\otimes_{\cB(\Hr)}\Hr$ and by Theorem \ref{isomalg} $$\tilde{V}^*U^*\big[\rho^{\cF(E)}\fJ\big](\cF(E)\otimes_{\rho}\Hr)=\h\fJ E(\bK)\otimes_{\cB(\Hr)}\Hr.$$
    Therefore,
    \[
\tilde{V}^*U^*(\tilde{G})=\cG.
    \]
    Now,
    \begingroup
\allowdisplaybreaks
    \begin{align*}
\big\|P_{\tilde{G}^{\perp}}U\tilde{V}\big(M_{\h{X}}\otimes I_{\Hr}\big)\tilde{V}^*U^*\big\|=& \big\|U\tilde{V}P_{\tilde{V}^*U^*(\tilde{G})^{\perp}}\big(M_{\h{X}}\otimes I_{\Hr}\big)\big(P_{\tilde{V}^*U^*(\tilde{G})}+P_{\tilde{V}^*U^*(\tilde{G})^{\perp}}\big)\big\|\\
=&\big\|U\tilde{V}P_{\cG^{\perp}}\big(M_{\h{X}}\otimes I_{\Hr}\big)P_{\cG^{\perp}}\big\|\\
=&\big\|P_{\cG^{\perp}}\big(M_{\h{X}}\otimes I_{\Hr}\big)P_{\cG^{\perp}}\big\|
\end{align*}
\endgroup
The second last equality follows because $\cG$ is invariant under $M_{\h{X}}\otimes I_{\Hr}$.
\end{proof}

Now we prove the following generalization of the above distance formula.
\begin{lemma}\label{distance formula}
    Let $\fJ\subseteq H^{\infty}(E)$ be a closed ideal and $\cG=\overline{\h\fJ E\big(\bK\big)\otimes\Hr}$. For a finite dimensional Hilbert space $\cE$ and for all $\bX=\sum_i X_i\otimes A_i\in H^{\infty}(E)\otimes\cB(\cE)$ we have\[
\text{dist}\,(\bX,\fJ\otimes\cB(\cE))=\Big\|\big(P_{\cG^{\perp}}\otimes I_{\cE}\big)\Big(\sum_i M_{\h{X}_i}\otimes I_{\Hr}\otimes A_i\Big)\big(P_{\cG^{\perp}}\otimes I_{\cE}\big)\Big\|.
\] 
\end{lemma}
\begin{proof}
    Since $\cE$ is finite dimensional, by \cite[Lemma 4.1.12]{Max}
    \[
    \text{dist}\,(\bX,\fJ)=\Big\|\big(P_{\tilde{G}^{\perp}}\otimes I_{\cE}\big)\Big(\sum_i X_i\otimes I_{\Hr}\otimes A_i\Big)\Big\|,
    \]
    where, $\tilde{G}=\overline{\rho^{\cF(E)}\cJ\big(\cF(E)\otimes\Hr\big)}$ as in the Proposition above. Following the computation as in the above proposition, we have
    \begingroup
\allowdisplaybreaks 
    \begin{align*} 
    \Big\|\big(P_{\tilde{G}^{\perp}}\otimes I_{\cE}\big)&\Big(\sum_i X_i\otimes I_{\Hr}\otimes A_i\Big) \Big\| \\
    =&\Big\|\sum_i P_{\tilde{G}^{\perp}}\big(X_i\otimes I_{\Hr}\big)\otimes A_i\Big\|\\
    =&\Big\|\sum_i P_{\cG^{\perp}}\big(M_{\h{X}_i}\otimes I_{\Hr}\big)P_{\cG^{\perp}}\otimes A_i\Big\|\quad(\text{using computation of Proposition \ref{prop:basic-distance-formula}})\\
    =&\Big\| \big(P_{\cG^{\perp}}\otimes I_{\cE}\big)\Big(\sum_i M_{\h{X}_i}\otimes I_{\Hr}\otimes A_i\Big)\big(P_{\cG^{\perp}}\otimes I_{\cE} \big)\Big\|
    \end{align*}
    \endgroup
    This completes the proof.
\end{proof}
For a matricial kernel $K$ and a set $\Omega\subseteq \bD$, we let
\[
E(K)_\Omega=\overline{\text{span}}\{K_{\fw,a,c}:\fz\in\Omega(\tau), a\in\cI(\tau,\rho), c\in\cB(\Hr,\Ht),\tau\in F\widehat{M}\}.
\]
\begin{lemma}
	For the matricial Szeg\"o kernel $\bK$ and for $\Omega\subseteq \bD$, $E(\bK)_{\Omega}^{\perp}\otimes\Hr$ is invariant under $\text{Mult}(E(\bK))\otimes B(\Hr)$ and $[\text{Mult}(E(\bK))\otimes B(\Hr)]'$.
\end{lemma}
\begin{proof}
	First of all, if $T\in \bB(E(\bK)\otimes \Hr)$ such that $T(M_{\h{X}}\otimes I_{\Hr})=(M_{\h{X}}\otimes I_{\Hr})T$ then using the identity 
    \[
    U\tilde{V}(M_{\h{X}}\otimes I_{\Hr})\tilde{V}^*U^*=X\otimes I_{\Hr}
    \]
    we have 
    \[
    U\tilde{V}T\tilde{V}^*U^*(X\otimes I_{\Hr})=(X\otimes I_{\Hr}) U\tilde{V}T\tilde{V}^*U^*
    \]
    Therefore, 
    \[
     U\tilde{V}T\tilde{V}^*U^*\in \gamma\big(H^{\infty}(E^{\rho})\big)
    \]
    where, $\gamma(\cdot):=U\iota^{\cF(F^{\sigma})}(\cdot)U^*$ is the representation of the algebra $H^{\infty}(E^{\rho})$ on the space $\cF(E)\otimes_{\rho}\Hr$ as in \cite[Theorem 3.9]{MS04}. Conversely, in a similar manner one can see, for any $S\in \gamma\big(H^{\infty}(E^{\rho})\big)$, $\tilde{V}^*U^*SU\tilde{V}\in[\text{Mult}(E(\bK)\otimes \Hr)]'$. So, $\gamma\big(H^{\infty}(E^{\rho})\big)$ is unitary equivalent to $[\text{Mult}(E(\bK)\otimes \Hr)]'$ via the unitary map $\tilde{V}^*U^*$. Recall that,
    \[
E(\bK)_\Omega=\overline{\text{span}}\{\bK_{\fw,a,c}:\fw\in\Omega(\tau), a\in\cI(\tau,\rho), c\in\cB(\Hr,\Ht),\tau\in F\widehat{M}\}.
\]
A function $f\in E(\bK)_{\Omega}^{\perp}$ if and only if
\[
\langle f, \bK_{\fw,a,c}\rangle=c^*f(\fw)(a^*)=0
\]
for all $\fw\in\Omega,\,a\in\cI(\tau,\rho),\, c\in\cB(\Hr,\Ht)$, for each $\tau\in F\widehat{M}$. This happens if and only if $f(\fw)=0$ for all $\fw\in \Omega$. So,
\[
E(\bK)_{\Omega}^{\perp}=\{f\in E(\bK): f(\fw)=0\}.
\]
 Clearly, $E(\bK)_{\Omega}^{\perp}\otimes \Hr$ is invariant under $[\text{Mult}(E(\bK)\otimes \Hr)]$. To show $E(\bK)_{\Omega}^{\perp}\otimes \Hr$ is invariant under $[\text{Mult}(E(\bK)\otimes \Hr)]'$ we proceed as follows. Let, $T\in [\text{Mult}(E(\bK)\otimes \Hr)]'$. Then by a previous argument $U\tilde{V}T\tilde{V}^*U^*\in\gamma(H^{\infty}(E^{\rho}))$ which implies $\tilde{V}T\tilde{V}^*\in \iota^{\cF(E^{\rho})}(H^{\infty}(E^{\rho}))$. So, $\tilde{V}T\tilde{V}^*=P\otimes I_{H_{\rho}}$ where $P\in H^{\infty}(E^{\rho})\subseteq \cL(\cF(E^{\rho}))$. Therefore, $\tilde{V}T=(P\otimes I_{H_{\rho}})\tilde{V}$. Recall that (see Remark \ref{directsum}), $\tilde{V}:=V\otimes I_{\Hr}:E(\bK)\otimes\Hr\to \cF(E^{\rho})\otimes_{\rho(M)'}\Hr$ is defined by,
 \[
 \tilde{V}(\bK_{\fw,a,cr}\otimes h)=\Big(\sum_{k=0}^{\infty}(I_k\otimes a)(\fw^*)^{(k)}c\Big)\otimes rh.
 \]
 For all $\eta\in E^{\rho}$,
 \[
 (L_{\eta}\otimes I_{\Hr})\tilde{V}(\bK_{\fw,a,cr}\otimes h)=\Big(\sum_{k=0}^{\infty}\eta\otimes (I_k\otimes a)(\fw^*)^{(k)}c\Big)\otimes rh\cong \Big(\sum_{k=0}^{\infty} (I_k\otimes \eta a)(\fw^*)^{(k)}c\Big)\otimes rh
 \]
 For all $d\in \rho(M)'$,
 \begin{align*}
 \big(\phi_{\infty}(d)\otimes I_{\Hr}\big)\tilde{V}(\bK_{\fw,a,cr}\otimes h)=&\Big(\sum_{k=0}^{\infty}\big(I_{k}\otimes da\big)(\fw^*)^{(k)}c\Big)\otimes rh \end{align*}

 \noindent\underline{Claim:} If $\bK_{\fw,a,cr}\otimes h\perp E(\bK)_{\Omega}\otimes \Hr$ then $$\tilde{V}^*(L_{\eta}\otimes I_{\Hr})\tilde{V}(\bK_{\fw,a,cr}\otimes h)\perp E(\bK)_{\Omega}\otimes \Hr$$
 and
 $$\tilde{V}^*\big(\phi_{\infty}(d)\otimes I_{\Hr}\big)\tilde{V}(\bK_{\fw,a,cr}\otimes h)\perp E(\bK)_{\Omega}\otimes \Hr.$$ 
 
 \noindent\underline{Proof of the claim:} It is enough to prove that, if  $\bK_{\fw,a,cr}\perp E(\bK)_{\Omega}$ then $$V^*L_{\eta}V(\bK_{\fw,a,cr})\perp E(\bK)_{\Omega}\quad \text{and}\quad V^*\phi_{\infty}(d)V(\bK_{\fw,a,cr})\perp E(\bK)_{\Omega}.$$ First note that, for all $\fw\in\bD(E,\sigma), \fw'\in\bD(\tau), a\in\cI(\sigma,\rho), c\in\cI(\Hr,\Hs),r\in B(\Hr), a'\in\cI(\tau,\rho), c'\in\cI(\Hr,\Ht)$ and $r'\in B(\Hr)$, 
 \begingroup
\allowdisplaybreaks
 \begin{align*}
 	\langle\bK_{\fw,a,cr},\bK_{\fw',a',c'r'}\rangle=&\Big\langle \Big(\sum_{k=0}^{\infty}(I_k\otimes a)(\fw^*)^{(k)}c\Big)\otimes r, \Big(\sum_{k=0}^{\infty}(I_k\otimes a')(\fw'^*)^{(k)}c'\Big)\otimes r' \Big\rangle\\
 	=&\sum_{k=0}^{\infty} r^*c^*\fw^{(k)}(I_k\otimes a^*a')\fw'^{*(k)}c'r'\\
 	=& r^*c^*\bK(\fw,\fw')(a^*a')c'r'.
 \end{align*}
 \endgroup
 On the other hand,
 \begingroup
\allowdisplaybreaks
 \begin{align*}
 	\langle V^*L_{\eta}V\bK_{\fw,a,cr},&\bK_{\fw',a',c'r'}\rangle\\
    =&\Big\langle \Big(\sum_{k=0}^{\infty}(I_k\otimes \eta a)(\fw^*)^{(k)}c\Big)\otimes r, \Big(\sum_{k=0}^{\infty}(I_k\otimes a')(\fw'^*)^{(k)}c'\Big)\otimes r' \Big\rangle\\
 	=&\sum_{k=1}^{\infty} r^*c^*\fw^{(k-1)}\Big((I_{k-1}\otimes \eta a)^*(I_k\otimes a')\Big)\fw'^{*(k)}c'r'\\
 	=&\sum_{k=1}^{\infty} r^*c^*\fw^{(k-1)}\Big((I_{k-1}\otimes a^*\eta^*)\big(I_{k-1}\otimes(I\otimes a')\fw'^*\big)\Big)\fw'^{*(k-1)}c'r'\\
 	=& r^*c^*\bK(\fw,\fw')[a^*\eta^*(I\otimes a')\fw'^*]c'r'.
 \end{align*}
 \endgroup
 Here observe that, $\eta^*(I\otimes a')\fw'^*\in\cI(\tau,\rho)$ as for all $m\in M$,
 \begingroup
\allowdisplaybreaks 
 \begin{align*}
 \eta^*(I\otimes a')\fw'^*\tau(m)=&\eta^*(I\otimes a')\big(\phi(m)\otimes I_{\Ht}\big)\fw'^*\\
 =& \eta^*\big(\phi(m)\otimes I_{\Hr}\big)(I\otimes a')\fw'^*\\
 =&\rho(m)\eta^*(I\otimes a')\fw'^*.
 \end{align*}
 \endgroup
 Also, for all $d\in\rho(M)'$,
 \begingroup
\allowdisplaybreaks
 \begin{align*}
 	\langle V^*\phi_{\infty}(d)&V\bK_{\fw,a,cr},\bK_{\fw',a',c'r'}\rangle\\
 	=&\Big \langle \Big(\sum_{k=0}^{\infty} \big(I_{k}\otimes da\big)(\fw^*)^{(k)}c\Big)\otimes r, \Big(\sum_{k=0}^{\infty}(I_k\otimes a')(\fw'^*)^{(k)}c'\Big)\otimes r' \Big \rangle\\
 	=& r^*c^*\bK(\fw,\fw')[(da)^* a')]c'r'.
 \end{align*}
 \endgroup
 Here observe that,  $da\in \cI(\sigma,\rho)$ for all $d\in \rho(M)'$ as $a\in \cI(\sigma,\rho)$.

 Now, for all $a,c,r,a',c',r'$ lie in respective spaces, $$\langle\bK_{\fw,a,cr},\bK_{\fw',a',c'r'}\rangle=0	 \implies 
 r^*c^*\bK(\fw,\fw')(a^*a')c'r'=0$$ which again implies $$r^*c^*\bK(\fw,\fw')[a^*\eta^*(I_E\otimes a')\fw'^*]c'r'=0\quad\text{and}\quad r^*c^*\bK(\fw,\fw')[(da)^* a')]c'r'=0.$$ Therefore, if
 $\langle\bK_{\fw,a,cr},\bK_{\fw',a',c'r'}\rangle=0$,  we can conclude $$\langle V^*L_{\eta}V\bK_{\fw,a,cr},\bK_{\fw',a',c'r'}\rangle=0=\langle V^*\phi_{\infty}(d)V\bK_{\fw,a,cr},\bK_{\fw',a',c'r'}\rangle.$$
 This proves the claim when we vary $\fw'\in\Omega$..
 
 Using the claim, via a limiting argument, we conclude that, if $f\perp E(\bK)_{\Omega}\otimes \Hr$ then 
 $$\tilde{V}^*(L_{\eta}\otimes I_{\Hr})\tilde{V}f \perp E(\bK)_{\Omega}\otimes \Hr\,\,(\forall\eta\in E^{\rho})\,\,\text{and}\,\, \tilde{V}^*(\varphi_{\infty}(d)\otimes I_{\Hr})\tilde{V}f\perp E(\bK)_{\Omega}\otimes \Hr.$$ Therefore, 
 $$\tilde{V}^*(F\otimes I_{\Hr})\tilde{V}f\perp E(\bK)_{\Omega}\otimes \Hr\quad (\forall\,F\in H^{\infty}(E^{\rho})).$$  Hence, $E(\bK)_{\Omega}^{\perp}\otimes \Hr$ is invariant under $[\text{Mult}(E(\bK)\otimes \Hr)]'$.
\end{proof}
We denote
\[
\fJ_{\Omega}:=\{F\in H^{\infty}(E): \h{F}(\fz)=0,\,\forall\,\fz\in\Omega\}
\]
that is 
\[
\fJ_{\Omega}=\{F\in H^{\infty}(E): \h{F}_\sigma(\zs)=0_{op}\in\cB(\Hs),\,\forall\,\zs\in\Omega(\sigma)\}.
\]
Using the fact that $\sigma\times\zs:H^{\infty}(E)\to\cB(\Hs)$ is a homomorphism one can see that $\fJ_\Omega$ is a two sided ideal in $H^{\infty}(E)$. Since $\fJ_\Omega$ is an ideal in $H^{\infty}(E)$, $\h{\fJ}_\Omega$ is an ideal in Mult$(E(K))$. The following lemma is crucial for our purpose.
\begin{lemma}\label{Gperp}
    Let $\Omega\subseteq\bD$ and $\fJ_\Omega$ be the ideal in $H^{\infty}(E)$, defined as above. Then for the matricial Szeg\"o kernel $\bK$, 
    \[
    \overline{\h{\fJ}_{\Omega}E(\bK)\otimes\Hr}=E(\bK)_{\Omega}^{\perp}\otimes \Hr.
    \]
\end{lemma}
\begin{proof}
    By \cite[Corollary 4.1.6]{Max} there is a complete lattice isomorphism $\mu$ from the lattice of ultraweakly closed two-sided ideals of $H^{\infty}(E)$ to the lattice of subspaces of $\cF(E)\otimes_{\rho}\Hr$ which are invariant under $H^{\infty}(E)$ and $H^{\infty}(E)'$. Indeed, for an ideal $\fJ\subseteq H^{\infty}(E)$, $$\mu(\fJ)=\overline{\fJ(\cF(E)\otimes_{\rho}\Hr)}.$$ 
    For $\cN\subseteq \cF(E)\otimes_{\rho}\Hr$ that is invariant under $H^{\infty}(E)$ and $H^{\infty}(E)'$,
    \[
    \mu^{-1}(\cN)=\big\{S\in H^{\infty}(E):S\big(\cF(E)\otimes_{\rho}\Hr\big)\subseteq \cN\big\}.
    \]
    By Theorem \ref{isomalg}, there is a unitary map $\tilde{V}^*U^*: \cF(E)\otimes_{\rho}\Hr\to E(\bK)\otimes_{\cB(\Hr)}\Hr$ (we have already used this map in Proposition \ref{prop:basic-distance-formula}) such that
    \[
    \tilde{V}^*U^*\overline{\big[\fJ\big(\cF(E)\otimes_{\rho}\Hr\big)\big]}=\overline{\h{\fJ}E(\bK)\otimes\Hr}.
    \]
    On the lattice of two-sided ideals in $H^{\infty}(E)$, define
    \[
    \tilde{\mu}:= \tilde{V}^*U^*\mu.
    \]
    Since, for a two-sided ideal $\fJ\subseteq H^{\infty}(E)$, $\overline{\h{\fJ}E(\bK)\otimes\Hr}$ is invariant under left and right multiplications, $\tilde{\mu}$ is a map from the lattice of two-sided ideals in $H^{\infty}(E)$ to the lattice of subspaces of $E(\bK)\otimes_{\cB(\Hr)}\Hr$ that are invariant under $\text{Mult}(E(\bK)\otimes \Hr)$ and $[\text{Mult}(E(\bK)\otimes \Hr)]'$.
        
    Recall the fact that $\tilde{V}^*U^*$ is a unitary and $H^{\infty}(E)$ and $\gamma\big(H^{\infty}(E^{\rho})\big)$ is unitary equivalent to $\text{Mult}(E(\bK)\otimes \Hr)$ and $[\text{Mult}(E(\bK)\otimes \Hr)]'$, respectively, via the unitary map $\tilde{V}^*U^*$. So, $\tilde{V}^*U^*$ is, in fact, a lattice isomorphism between the lattice of closed subspaces of $\cF(E)\otimes_{\rho}\Hr$ that are invariant under $H^{\infty}(E)$ and $H^{\infty}(E)'$ and the lattice of closed subspaces of $E(\bK)\otimes_{\cB(\Hr)}\Hr$ that are invariant under $\text{Mult}(E(\bK)\otimes \Hr)$ and $[\text{Mult}(E(\bK)\otimes \Hr)]'$.    
    Therefore, $\tilde{\mu}$ is a lattice isomorphism, where
    \[
    \tilde{\mu}(\fJ)=\overline{\h{\fJ}E(\bK)\otimes\Hr}
    \]
    and for a closed invariant (under left multiplication operators and its commutators) subspace $M$ of $E(\bK)\otimes\Hr$,
    \begin{align*}
    \tilde{\mu}^{-1}(M)=&\big\{S\in H^{\infty}(E):\tilde{V}^*U^*\big[S\big(\cF(E)\otimes_{\rho}\Hr\big)\big]\subseteq M\big\}\\
    =&\big\{S\in H^{\infty}(E):M_{\h{S}}E(\bK)\otimes_{\rho}\Hr\subseteq M\big\}.
    \end{align*}
    Let $\cM:=E(\bK)_{\Omega}^{\perp}\otimes \Hr$.  
  Due to the set theoretic representation of $E(\bK)_{\Omega}^{\perp}$, which is
 \[
E(\bK)_{\Omega}^{\perp}=\{f\in E(\bK): f(\fw)=0\}.
\]
we deduce,
\begin{align*}
    \tilde{\mu}^{-1}(\cM)=&\big\{S\in H^{\infty}(E):M_{\h{S}}f=0\,\,\text{on}\,\Omega,\,\forall f\in E(\bK)\big\}\\
    =&\big\{S\in H^{\infty}(E):\h{S}=0\,\,\text{on}\,\Omega\big\}\\
    =&\fJ_{\Omega}.
\end{align*}
This implies,
\[
\cM=\tilde{\mu}\big[\tilde{\mu}^{-1}(\cM)\big]=\overline{\h{\fJ}_{\Omega}E(\bK)\otimes\Hr}
\]
and so 
\[
E(\bK)_{\Omega}^{\perp}\otimes \Hr=\overline{\h{\fJ}_{\Omega}E(\bK)\otimes\Hr}.
\]
Hence, the result follows.
\end{proof}
\begin{theorem}\label{complete Pick}
$\bK$ is a complete Pick kernel.
\end{theorem}
\begin{proof}
    As we have discussed earlier, we need only to prove that item (II) of the complete Pick property implies item (I). We will prove this for a finite set $\Omega$. The result for an arbitrary set follows from Theorem \ref{isomalg} and a compactness argument. 

    Let $\Omega=\{\fz_1,\ldots,\fz_n\}$ and $f_0(\fz_i)=W_i$, for all $i=1,\ldots,n$. By taking direct sums and using the assumption that $\fz\to W$ extends to a matricial function on $\Omega'$, we may assume that $\Omega=\{\fz\}$ and $f_0(\fz)=W\in\cB(\Hs\otimes \cE)$. Since $f_0$ extends to a matricial function on $\Omega'$, by Lemma \ref{lem:W-in-alg}, there is some $P\in\cT_0(E)\otimes\cB(\cE)$ such that $W=\h{P}(\fz)$.

    Assume that $\bK_{f_0}$, as in the item (II) in the complete Pick property, is a cp kernel in $\Omega$. We need to find a $f\in H^{\infty}(E)\otimes \cB(\cE)$ such that $\|M_{\h{f}}\|\leq 1$ and $\h{f}(\fz)=W$. Recall that, $\fJ_\Omega$ as defined above is an ideal in $H^{\infty}(E)$. If we can find $G\in\fJ_\Omega\otimes\cB(\cE)$ such that $\|P+G\|\leq 1$ then putting $f=P+G$ we will be done. Since $\fJ_\Omega$ is weak-$*$ closed, it suffices to prove that dist$\big(P, \fJ_\Omega\otimes\cB(\cE)\big)\leq 1$.

    Let $\cG=\overline{\h\fJ_\Omega E\big(\bK\big)\otimes\Hr}$. The space $\cG$ is invariant under Mult$(E(\bK)\otimes\Hr)$. By Lemma \ref{distance formula} 
    \[
    \text{dist}(P,\fJ_\Omega\otimes\cB(\cE))=\Big\|(P_{\cG^{\perp}}\otimes I_{\cE})\Big(\sum_i M_{\h{P}_i}\otimes I_{\Hr}\otimes A_i\Big)(P_{\cG^{\perp}}\otimes I_{\cE})\Big\|,
    \]
    where $P=\sum_i P_i\otimes A_i$ for some $P_i\in H^{\infty}(E)$ and $A_i\in\cB(\cE)$. We will compute the norm in the right-hand side of the above. By Lemma \ref{Gperp}, $\cG^{\perp}=E(\bK)_\Omega\otimes\Hr$.  Let $g=\sum_{i}K_{\fz, a_i, c_i\otimes h_i\otimes\epsilon_i}=\sum_{i}K_{\fz, a_i, c_i}\otimes h_i\otimes\epsilon_i$ be an element of the correspondence $E(\bK)_\Omega\otimes\Hr\otimes\cE$. Now, using $\h{P}(\fz)=f_0(\fz)$, we calculate
    \begingroup
\allowdisplaybreaks
    \begin{align*}
        \|g\|^2-&\Big\|\Big(\sum_i M^*_{\h{P}_i}\otimes I_{\Hr}\otimes A_i\Big)g\Big\|^2\\
        =&\sum_{i,j}\Big[(c_i\otimes\epsilon_i\otimes h_i)^*\big(\bK(\fz,\fz)(a_ja_i^*)\otimes I_{\cE}\otimes I_{\Hr}\big) (c_j\otimes \epsilon_j\otimes h_j)\\
        & \quad\quad-(c_i\otimes\epsilon_i\otimes h_i)^*\Big(f_0(\fz)\big(\bK(\fz,\fz)(a_ja_i^*)\otimes I_{\cE}\big)f_0(\fz)^*\otimes I_{\Hr}\Big)(c_j\otimes \epsilon_j\otimes h_j)\Big]\\
        \geq& 0 \quad(\text{by the assumption}).
    \end{align*}
    \endgroup
    (In the above, we have used the isomorphism between the Hilbert spaces $\Hr\otimes\cE$ and $\cE\otimes\Hr$ that takes $h_{\rho}\otimes\epsilon$ to $\epsilon\otimes h_{\rho}$). This implies, $$\Big\|\Big(\sum_i M^*_{\h{P}_i}\otimes I_{\Hr}\otimes A_i\Big)\Big|_{\cG^{\perp}\otimes\cE}\Big\|\leq 1.$$
    But since
    \[
    \Big\|\Big(\sum_i M^*_{\h{P}_i}\otimes I_{\Hr}\otimes A_i\Big)\Big|_{\cG^{\perp}\otimes\cE}\Big\|=\Big\|(P_{\cG^{\perp}}\otimes I_{\cE})\Big(\sum_i M^*_{\h{P}_i}\otimes I_{\Hr}\otimes A_i\Big)(P_{\cG^{\perp}}\otimes I_{\cE})\Big\|
    \]
    it follows that 
    \[
    \text{dist}(P,\fJ_\Omega\otimes\cB(\cE))=\Big\|(P_{\cG^{\perp}}\otimes I_{\cE})\Big(\sum_i M_{\h{P}_i}\otimes I_{\Hr}\otimes A_i\Big)(P_{\cG^{\perp}}\otimes I_{\cE})\Big\|\leq 1.
    \]
    This completes the proof.
\end{proof}

\section{Application of Complete Pick property}\label{sec:AppCompPick}
Now we discuss some applications of the complete Pick property in the context of matricial variety. Let $\cF$ be a family of matricial functions on $\bD$, corresponding to some kernel $K$. Define, 
\[
V(\cF)=\{\fw\in\bD:f(\fw)=0,\,\forall\,f\in \cF\}.
\]
We say, a set $\cV\subseteq \bD$ is a matricial variety corresponding to some kernel $K$  if it is the joint zero set of a family of multipliers in $\text{Mult}\,E(K)$, that is, $\cV=V(\cF)$ for some $\cF\subseteq \text{Mult}\,E(K)$. The following lemma is useful for our purpose which is the converse of Lemma \ref{adjointS}.
\begin{lemma}\label{mult-action}
	Let $S:\bD\to \coprod_{\sigma\in F\widehat{M}}\cB(\Hs)$ be a matricial function and $T$ be an adjointable operator on $E(K)$ such that 
	\[
	TK_{\fw,a,c}=K_{\fw,a,S(\fw)^*c}
	\]
	for all $\fw,a,c$. Then,  $T=M_S^*$.
\end{lemma}
\begin{proof}
	For all $\fw,a,c$, we compute
	\begin{align*}
		\langle K_{\fw,a,c}, T^*h\rangle=\langle TK_{\fw,a,c}, h\rangle=\langle K_{\fw,a,S(\fw)^*c}, h\rangle=c^*S(\fw)h(\fw)(a^*)=\langle K_{\fw,a,c}, S\cdot h\rangle
	\end{align*}
	for all $h\in E(K)$. Hence the result follows.
\end{proof}
We also need the following lemma.
\begin{lemma}\label{restriction}
	Let $K$ be a matricial kernel on $\bD$ and $\Omega\subseteq \bD$. Then for a multiplier $S$ of $E(K)$, the restriction $S|_{\Omega}\in \text{Mult}\,E(K)_{\Omega}$. Moreover, 
	\[
	M_{S|_{\Omega}}=P_{\Omega}M_S=P_{\Omega}M_SP_{\Omega},
	\]
	where, $P_{\Omega}:E(K)\to E(K)_{\Omega}$ is the orthogonal projection. 
\end{lemma}
\begin{proof}
	By Lemma \ref{adjointS}, for all $\fw\in\Omega$, $M_S^*K_{\fw,a,c}=K_{\fw,a,S(\fw)^*c}$. Therefore, $E(K)_{\Omega}$ is invariant under $M_S^*$. Using Lemma \ref{mult-action} we conclude that $S|_{\Omega}$ is a multiplier and $M^*_{S|_{\Omega}}=M^*_S|_{E(K)_{\Omega}}$.
\end{proof}
For $\Omega\subseteq \bD$, recall the following notation,
\[
\fJ_{\Omega}=\{F\in H^{\infty}(E): \h{F}_\sigma(\zs)=0_{op}\in\cB(\Hs),\,\forall\,\zs\in\Omega(\sigma)\}.
\]
Considering the matricial Szeg\"o kernel $\bK$, the following theorem is an application of complete Pick property.
\begin{theorem}
	Let $\cV\subseteq\bD$ be a matricial variety corresponding to $E(\bK)$. Then the map $S\mapsto S|_{\cV}$ is a completely contractive and surjective homomorphism  $\text{Mult}\,E(\bK)\to \text{Mult}\,E(\bK)_{\cV}$, which induces a completely isometric isomorphism $\text{Mult}\,E(\bK)/\h{\fJ}_{\cV}\to \\\text{Mult}\,E(\bK)_{\cV}$.
\end{theorem}

\begin{proof}
	By Lemma \ref{restriction} the map $S\mapsto S|_{\cV}$ is completely contractive and it is easy to verify that the kernel of this map is $\h{\fJ}_{\cV}$. To show the remaining parts we use the complete Pick property. 
	
	Let $\cE$ be a finite dimensional Hilbert space and $S_0\in \text{Mult}(E(\bK)_{\Omega}\otimes\cE)$ such that \\$\|M_{S_0}\|_{\text{Mult}(E(\bK)_{\Omega}\otimes\cE)}\leq 1$. By Proposition \ref{SK} the kernel 
	\[
\bK_{S_0}(\fz,\fw)(a)=\bK(\fz,\fw)(a)\otimes I_{\cE}-S_0(\fz)\big[\bK(\fz,\fw)(a)\otimes I_{\cE}\big]S_0(\fw)^*
\]
is a cp matricial kernel on $\cV$. By Theorem \ref{complete Pick}, $S_0$ extends to a multiplier $S\in \text{Mult}(E(\bK)\otimes\cE)$ such that $\|M_S\|_{\text{Mult}(E(\bK)\otimes\cE)}\leq 1$. In particular the restriction map $S\mapsto S|_{\cV}$ is surjective and it induces a completely contractive linear isomorphism $\Phi: \text{Mult}\,E(\bK)/\h{\fJ}_{\cV}\to \text{Mult}\,E(\bK)_{\cV}$. In fact, in the above one can find the multiplier $S$ with $$\|M_S\|_{\text{Mult}(E(\bK)\otimes\cE)}=\|M_{S_0}\|_{\text{Mult}(E(\bK)_{\Omega}\otimes\cE)}.$$
 Indeed one can consider $S_0'=\frac{S_0}{\|M_{S_0}\|}$. Since $\|M_{S_0'}\|=1$, using complete Pick property, we get an extension $S'$ with $\|M_{S'}\|=1$. Take $S=\|M_{S_0}\|S'$ and observe that the above norm equality follows. Since, $\cE$ is a finite dimensional Hilbert space with arbitrary dimension the above norm equality give us that $\Phi$ is in fact a complete isometry.
\end{proof}
Below is an immediate remark which follows from above.
\begin{remark}\label{remark on multipiler}
	For a matricial variety $\cV\subseteq\bD$, $\text{Mult}\,E(\bK)_{\cV}=\text{Mult}\,E(\bK)|_{\cV}$ and for every $T\in \text{Mult}\,E(\bK)_{\cV}$ there exists $S\in \text{Mult}\,E(\bK)$ such that $T=S|_{\cV}$ and $\|T\|_{\text{Mult}\,E(\bK)_{\cV}}=\|S\|_{\text{Mult}\,E(\bK)}$.
\end{remark}

\begin{theorem}\label{HinfinityV}
	Suppose $\cV$ be a matricial variety with respect to some ideal $\fJ\subseteq H^{\infty}(E)$, that is,
	\[
	\cV=\{\fw\in\bD:\h{F}(\fw)=0,\,\forall\,F\in \fJ\}.
	\]
	Then, 
	\[
	\text{Mult}\,E(\bK)_{\cV}=H^{\infty}(\cV),
	\]
	where $H^{\infty}(\cV)$ is the set of all bounded holomorphic \normalcolor  matricial functions on $\cV$.
\end{theorem}
\begin{proof}
	Let $\pi$ be an infinite multiple of a faithful normal representation of $M$. We refer to $\pi$ as a universal generator and every $\sigma\in\,NRep\, (M)$ is a sub representation of $\pi$. 
	
	Let $\mu$ be the map as in the proof of Lemma \ref{Gperp}. The space $\mu(\fJ)\subseteq \cF(E)\otimes H_{\pi}$ is invariant under $\pi^{\cF(E)}(H^{\infty}(E))$. Let $P: \cF(E)\otimes H_{\pi}\to \mu(\fJ)^{\perp}$ be the orthogonal projection. Denote $\cH:= \cF(E)\otimes H_{\pi}$ and write $t_1$ for the map 
	$$t_1:E \rightarrow B(\cH)$$ defined by 
	\[ t_1(\xi)=P(T_{\xi}\otimes I_{H_{\pi}})P .\]
	We also define 
	\[ \sigma_P:M \rightarrow B(\cH) \] by
	\[\sigma_P(a)=P(\varphi_{\infty}(a)\otimes I_{H_{\pi}})P .\]
	
	In order to show that $(\sigma_P,t_1)$ is a representation of $(M,E)$, we note first that $P$ commutes with $\varphi_{\infty}(M) \otimes I_{H_{\pi}}$.
	
	Clearly, $(\varphi_{\infty}(M) \otimes I_{H_{\pi}})P^{\perp}=P^{\perp}(\varphi_{\infty}(M) \otimes I_{H_{\pi}})P^{\perp}$ and ,by taking adjoints, we find that $\varphi_{\infty}(M) \otimes I_{H_{\pi}}$ commutes with $P^{\perp}$, and,thus, with $P$.
	
	This shows that $\sigma_P$ is a representation of $M$ and we have ,for $a,b\in M$,
	$$ t_1 (\varphi(a)\xi b)=P(T_{\varphi(a)\xi b}\otimes I_{H_{\pi}})P=P(\varphi_{\infty}(a)\otimes I_{H_{\pi}})(T_{\xi}\otimes I_{H_{\pi}})(\varphi_{\infty}(b)\otimes I_{H_{\pi}})=\sigma_P(a) t_1(\xi) \sigma_P(b) $$ showing  that $(\sigma_P,t_1)$ is a representation of $(M,E)$.

	Now, we set $t_0=\tilde{t_1}$. Thus $t_0:E\otimes_{\sigma_P}\cH\to\cH$ by
	\[
	t_0(\xi\otimes h)=P(T_{\xi}\otimes I_{H_{\pi}})Ph\quad (\xi\in E,h\in \cH).
	\] 
	It is clear that, $\|t_0\|\leq 1$. Therefore, $rt_0\in \bD(\sigma_P)$ for all $0\leq r<1$.

	Observe that,
	\[
	t_0(T_{\xi}\otimes I_{H_{\pi}})=P(T_{\xi}\otimes I_{H_{\pi}})P,
	\] 
	for all $\xi\in E$. Also,  
	\begin{equation}\label{Emult}
	\h{T}_{\xi}(t_0)h=t_0(\xi\otimes h)=P(T_{\xi}\otimes I_{H_{\pi}})Ph. \end{equation}
	 and 
	\begin{equation}\label{Mmult} \widehat{\varphi_{\infty}(a)}(t_0)=\sigma_P(a)=P(\varphi_{\infty}(a)\otimes I_{H_{\pi}})P.\end{equation}

	The projection $P^{\perp}$ is invariant under $H^{\infty}(E)$ and, thus, $P$ is coinvariant and it follows that the map $F\mapsto P(T_{\xi}\otimes I_{H_{\pi}})P$ is multiplicative on $H^{\infty}(E)$. It, thus , follows from (\ref{Mmult}) and (\ref{Emult}) that, for every $F\in H^{\infty}(E)$,
	\[\h{F}(t_0)=P(F\otimes I_{H_{\pi}})P.\] and, for $0\leq r<1$
	\begin{equation}\label{Fmult}
	\h{F}(rt_0)=rP(F\otimes I_{H_{\pi}})P.\end{equation}

Suppose $0<r<1$ and $D\in \big( P\pi^{\cF(E)}(H^{\infty}(E))P\big)'$. So $D$ commutes with $P(T_{\xi}\otimes I_{H_{\pi}})P$. Then, for all $\xi\in E$ and $h\in H_{\pi}$,
\begin{align*}
	Drt_0(\xi\otimes h)=rDP(T_{\xi}\otimes I_{H_{\pi}})Ph=rP(T_{\xi}\otimes I_{H_{\pi}})PDh&=rt_0(\xi\otimes Dh)\\
	&=rt_0(I_E\otimes D)(\xi\otimes h).
\end{align*}
Thus, $D(rt_0)=(rt_0)(I_E\otimes D)$. For $f\in H^{\infty}(\cV)$ since $f$ is matricial, $f(rt_0)D=Df(rt_0)$ and that implies $f(rt_0)\in P\pi^{\cF(E)}(H^{\infty}(E))P$. Therefore, there exists some $Y_r\in H^{\infty}(E)$ such that $f(rt_0)=P\pi^{\cF(E)}(Y_r)P$, which is isomorphic to some $\tilde{Y_r}\in H^{\infty}(E)/\fJ$ by \cite[Corollary 4.1.8 and Theorem 4.1.9]{Max}. Then it follows that 
\[
\|\tilde{Y}_r\|=\big\|P\pi^{\cF(E)}(Y_r)P\big\|\leq \sup_{\fz\in\cV}\|f(\fz)\|<\infty.
\] 
Hence $\{\tilde{Y}_r:0<r<1\}$ is bounded.

Denote, $\cV(\sigma):=\{\fw\in\bD(E,\sigma):\h{F}(\fw)=0,\,\forall F\in\fJ \}$. 

\noindent\underline{Claim:} For all  $0< r\leq 1$ and $\fz\in r\bD(E,\sigma)\cap r\cV(\sigma)$, there exist isometries $V_r:H_{\sigma}\to\cF(E)\otimes H_{\sigma}$ such that
\begin{equation}\label{dilation of creation}
V_r(\frac{1}{r}L_{\xi}^*\fz^*)=(L_{\xi}^*\otimes I_{H_{\sigma}})PV_r\quad(\forall\,\xi\in E)
\end{equation}
and consequently, for all $F\in H^{\infty}(E)$, 
\begin{equation}\label{dilation of F}
\h{F}(\frac{1}{r}\fz)=V_r^*P(F\otimes I_{H_{\sigma}})PV_r=V_r^*P\sigma^{\cF(E)}(F)PV_r
\end{equation}

\noindent\underline{Proof of the claim:} To prove the claim it is enough to prove that for all $\fz\in \bD(E,\sigma)\cap \cV(\sigma)$ the isometry $K(\fz):H_{\sigma}\to\cF(E)\otimes H_{\sigma}$, defined in \cite{MS16}, by
\[
K(\fz)h= \sum_{k=0}^{\infty} \big(I_{E^{\otimes k}}\otimes \Delta_*(\fz)\big)\fz^{(k)*}h
\]
such that 
\[
K(\fz)L_{\xi}^*\fz^*=(L_{\xi}^*\otimes I_{H_{\sigma}})PK(\fz)\quad(\forall\,\xi\in E).
\]
Since, $\fz\in \bD(E,\sigma)$, from \cite{MS16} we know that,
\[
K(\fz)L_{\xi}^*\fz^*=(L_{\xi}^*\otimes I_{H_{\sigma}})K(\fz) \quad(\forall\,\xi\in E).
\]
Using the fact that $\fz\in \cV$ we will show
\[
(L_{\xi}^*\otimes I_{H_{\sigma}})K(\fz)=(L_{\xi}^*\otimes I_{H_{\sigma}})PK(\fz),
\]
that is to show, 
\begin{equation}\label{eKz}
(L_{\xi}^*\otimes I_{H_{\sigma}})P^{\perp}K(\fz)=(L_{\xi}^*\otimes I_{H_{\sigma}})P_{\mu({\fJ})}K(\fz)=0
\end{equation}
for all $\xi\in E$.

To prove it we use  \cite[Theorem 13]{MS08}. There we have 
\begin{equation}\label{Poisson} 
	K(\fz)\h{X}(\fz)^*=(X^*\otimes I_{H_{\sigma}})K(\fz).
	\end{equation}
	(Note that $\eta$ there should be translated into $\fz^*$ here.)
	
	Thus, if $X\in \fJ$, it follows from the assumption that $\fz\in \cV(\sigma)$ that $\h{X}(\fz)=0$ and, thus, 
	$$(X^*\otimes I_{H_{\sigma}})K(\fz)=0$$
	and, for $h,g \in H_{\sigma}$,
	$$\langle K(\fz)h, (X\otimes I_{H_{\sigma}})g\rangle= \langle (X^*\otimes I)K(\fz)h,g\rangle =0$$
	proving (\ref{eKz}) since $P_{\mu({\fJ})}K(\fz)=0$.

Hence the claim is proved by denoting $K(\frac{1}{r}\fz)$ by $V_r$.

Since $\sigma$ is a subrepresentation of $\pi$, there is an isometry $w: H_{\sigma}\to  H_{\pi}$ such that
\[\sigma(a)=w^*\pi(a)w\] for $a\in M$. It follows that $v: x\otimes h\mapsto x \otimes wh$ is an isometry for $x\in \cF(E)$ and $h\in H_{\sigma}$. In fact, noting that $P=P_0\otimes I$ for a projection $P_0$ on $\cF(E)$, we have
$v(T_{\xi}^*\otimes I)P(x\otimes h)=v(T_{\xi}^*P_0x\otimes h)=T_{\xi}^*P_0x \otimes wh=(T_{\xi}^*\otimes I)Pv(x\otimes h)$ and, thus
\begin{equation}\label{v1}
v(T_{\xi}^*\otimes I_{H_{\sigma}})P=(T_{\xi}^*\otimes I_{H_{\sigma}})Pv\quad(\text{for all} \,\xi\in E).
\end{equation}
Write $v_r=vV_r$. Then $v_r: H_{\sigma}\to \cF(E)\otimes H_{\pi}$ and, applying (\ref{dilation of creation}) and (\ref{v1}), 
\[
v_r(\frac{1}{r}T_{\xi}^*\fz^*)=(T_{\xi}^*\otimes I_{H_{\pi}})Pv_r\quad(\forall\,\xi\in E).
\] 
Since $\mu(\fJ)$ is coinvariant for $H^{\infty}(E)$ and, thus, invariant for $T_{\xi}^*$, we have 
\[
	v_r(\frac{1}{r}T_{\xi}^*\fz^*)=P(T_{\xi}^*\otimes I_{H_{\pi}})Pv_r
\]
Taking adjoint in the both sides,
\begin{equation}\label{v_r}
	\fz T_{\xi}v_r^*=rv_r^*P(T_{\xi}\otimes I_{H_{\pi}})P\quad(\forall\,\xi\in E)
\end{equation}
Observe that, for all $\eta\in\cF(E)$ and $h_{\pi}\in H_{\pi}$,
\[
T_{\xi}v_r^*(\eta\otimes h_\pi)=\xi\otimes v_r^*(\eta\otimes h_\pi)=(I_E\otimes v_r^*)(\xi\otimes\eta\otimes h_{\pi})=(I_E\otimes v_r^*)(T_{\xi}\otimes I_{H_{\pi}})(\eta\otimes h_\pi),
\]
that is,
\[
T_{\xi}v_r^*=(I_E\otimes v_r^*)(T_{\xi}\otimes I_{H_{\pi}}).
\]
Therefore, it follows from (\ref{v_r}) that,
\begingroup
\allowdisplaybreaks
\begin{align*}
	& \fz(I_E \otimes v_r^*) (T_{\xi}\otimes I_{H_{\pi}})=r v_r^*t_0(T_{\xi}\otimes I_{H_{\pi}})\quad(\forall\,\xi\in E)\\&\Rightarrow \fz(I_E \otimes v_r^*) (T_{\xi}\otimes I_{H_{\pi}})h=r v_r^*t_0(T_{\xi}\otimes I_{H_{\pi}})h\quad(\forall\,\xi\in E,\,h\in\cH)\\
	&\Rightarrow \fz(I_E \otimes v_r^*)(\xi \otimes h)=r v_r^*t_0(\xi\otimes h)\quad(\forall\,\xi\in E,\,h\in\cH).
\end{align*}
\endgroup
Hence, for all $0<r\leq 1$,
\[
\fz(I_E \otimes v_r^*)=v_r^*(rt_0).
\]
Since $f$ is matricial, for all $0<r< 1$, we get  
\[
f(\fz) v_r^*=v_r^*f(rt_0)
\]
and that gives 
\begin{equation}\label{f-Y_r}
	f(\fz)=v_r^*f(rt_0) v_r=V_r^*v^*f(rt_0)vV_r=V_r^*v^*P\pi^{\cF(E)}(Y_r)PvV_r
\end{equation}
  
 In above we have seen $\{\tilde{Y}_r:0<r<1\}$ is a bounded set in $H^{\infty}(E)/\fJ$. So, we can fix a sequence $\{r_n\}$ increasing to $1$ and find a subnet $\{\tilde{Y}_{r_\alpha}\}$ such that $\{\tilde{Y}_{r_\alpha}\}\to \tilde{Y}$ in the $w^*$-topology for some $\tilde{Y}\in H^{\infty}/\fJ$. Therefore, $P\pi^{\cF(E)}(Y_{r_\alpha})P\to P\pi^{\cF(E)}(Y)P$.
 
 Now fix $\fz\in \bD(E,\sigma)\cap\cV(\sigma)$. Then $r\fz\in r\bD(E,\sigma)\cap r\cV(\sigma)$. From the relation (\ref{f-Y_r}) we have 
 \begin{equation}\label{f-Yr}
	f(r\fz)=V_r^*v^*P\pi^{\cF(E)}(Y_r)PvV_r
\end{equation}

Being an element in $H^{\infty}(\cV)$, $f$ is holomorphic and hence continuous. Using relation (\ref{f-Yr}), via a limiting argument  we get,
\begingroup
\allowdisplaybreaks
\begin{align*}
	f(\fz)=\lim_{r_\alpha\to 1} f(r_{\alpha} \fz)&=\lim_{r_\alpha\to 1}V_{r_\alpha}^*v^*P\pi^{\cF(E)}(Y_{r_\alpha})PvV_{r_\alpha}\\
	&= V_1^*P\sigma^{\cF(E)}(Y)PV_1\\
	&= \h{Y}(\fz)\quad(\text{by}\,(\ref{dilation of F})).
\end{align*}
\endgroup
Thus, $f\in\text{Mult}\,E(\bK)|_\cV$. Since by Remark \ref{remark on multipiler} $\text{Mult}\,E(\bK)|_\cV=\text{Mult}\,E(\bK)_\cV$, we conclude that $f\in\text{Mult}\,E(\bK)_\cV$. Hence we prove that $H^{\infty}(\cV)\subseteq \text{Mult}\,E(\bK)_\cV$.

For the reverse inclusion, let $f\in \text{Mult}\,E(\bK)_\cV=\text{Mult}\,E(\bK)|_\cV$. Then there exists $S\in\text{Mult}\,E(\bK)$ such that $f=S|_{\cV}$. By Corollary \ref{Sbounded}, $S$ is bounded on $\bD$ with respect to the supremum norm. Therefore,  $f$ is bounded on $\cV$ with respect to the supremum norm and hence $f\in H^{\infty}(\cV)$. Thus, $\text{Mult}\,E(\bK)_\cV\subseteq H^{\infty}(\cV)$. This completes the proof of the theorem.
	 \end{proof}

\end{document}